# NONSINUSOIDAL PERIODIC FOURIER SERIES

Arpad Török[1], Stoian Petrescu[2], Michel Feidt[3]

[1] PhD student, The Polytechnic University of Bucharest, Department of Engineering Thermodynamics, 313, Splaiul Independentei, 060042 Bucharest, Romania, e-mail: arpi_torok@yahoo.com
[2] Prof. Dr. Eng., Polytechnic University of Bucharest, Department of Engineering Thermodynamics, România
[3] Prof. Dr. Eng., L.E.M.T.A., U.R.A. C.N.R.S. 7563, Université de Loraine Nancy 12, avenue de la Foret de Haye, 54516 Vandeuvre-lès-Nancy, France

**Abstract.** According to harmonic analysis (Fourier analysis), any function $f(x)$, periodic over the interval $[-L, L]$, which satisfies the Dirichlet conditions, can be developed into an infinite sum (known in the literature as the **trigonometric series,** and for which, for reasons which will become evident in the course of this work, we will use the name of **sinusoidal series**), consisting of the weighted components of a complete biortogonal base, formed of the unitary function *1*, of the fundamental harmonics $sin(\pi x/L)$, even and $cos(\pi x/L)$, odd ($2L$-periodic functions) and of the secondary harmonics $sin(n\pi x/L)$ and $cos(n\pi x/L)$ (periodic functions, with period $2L/n$, where $n \in \mathbf{Z}^+$, positive integers). The coefficients of these expansions (Fourier coefficients) can be calculated using Euler formulas. We will generalize this statement and show that the function $f(x)$ can also be developed into **non-sinusoidal periodic series**, formed from the sum of the weighted components of a complete, non-orthogonal base: the unit function *1*, the fundamental quasi-harmonics $g(x)$, even and $h(x)$, odd ($2L$-periodic functions, with zero mean value over the definition interval) and the secondary quasi-harmonics $g_n(x)$ and $h_n(x)$ ($2L/n$-periodic functions), where $n \in \mathbf{Z}^+$. The fundamental quasi-harmonics $g(x)$ and $h(x)$ are any functions which admit expansions in sinusoidal series (satisfy Dirichlet conditions, or belong to $L^2$ space). The coefficients of these expansions are obtained with the help of certain algebraic relationships between the Fourier coefficients of the expansions of the functions $f(x)$, $g(x)$ and $h(x)$. In addition to their obvious theoretical importance, these types of expansions can have practical importance in the approximation of functions and in the numerical and analytical resolution of certain classes of differential equations.

## 1. Introduction

Let $g(x)$ be any function of R, defined on a real interval $I$ (open, closed or semi-open), delimited by $x_1$ and $x_{m+1}$, introduced by a single expression $g(x)=g_1(x)$ on any $I$ ($m=1$), or by $m$ different expressions: $g(x)=g_i(x)$, $i=1, 2, ..., m$, each valid on one of the adjacent disjoint subintervals bounded by $x_i$ and $x_{i+1}$, $1\leq i \leq m$. To define simultaneously the function $g(x)$, as well as its definition domain, we will introduce a system to simultaneously mark the limits $x_i$ of (sub)intervals of definition and also the expressions of the function valid on these (sub)intervals:

$g(x) = G[x_1^+ > g_1 < x_2^-]$, or $g(x) = G[x_1 \geq g_1 \leq x_2]$, or $g(x) = G[x_1^+ * (g_1) * x_2^-]$, or
  $g(x) = G[x_1^+ (g_1) x_2^-]$, for $I$ closed,

$g(x) = G[x_1^+ > g_1 < x_2 \cup x_2^+ > ... < x_i \cup x_i^+ > g_i < x_{i+1} \cup x_{i+1}^+ > ... < x_m \cup x_m^+ > g_m < x_{m+1}^-]$, or

  $g(x) = G[x_1 \geq g_1 < x_2 \cup x_2 \geq ... < x_i \cup x_i \geq g_i < x_{i+1} \cup x_{i+1} \geq ... < x_m \cup x_m \geq g_m \leq x_{m+1}]$, or

  $g(x) = G[x_1^+(g_1)x_2 \cup x_2^+(g_2)...(g_{i-1})x_i \cup x_i^+(g_i)x_{i+1} \cup x_{i+1}^+(g_{i+1})...(g_{m-1})x_m \cup x_m^+(g_m)x_{m+1}^-]$, etc.

  for $I= [x_1, x_2) \cup [x_2, x_3) \cup ... \cup [x_i, x_{i+1}) \cup ... \cup [x_m, x_{m+1}]$

The upper index attached to the limits $x_i$ of these intervals mean:
- $+$ : $g(x_i) = \lim_{x \to x_i^+} g(x)$
- $-$ : $g(x_i) = \lim_{x \to x_i^-} g(x)$
- without index: the function $g(x)$ is undefined in $x_i$, or it has a fixed value $g(x_i)=a \in \mathbf{R}$.

Therefore:
- for $g(x)$ $):(x_1, x_2)$, we will use one of the notations:
$g(x)=G[x_1>g<x_2]$, or $g(x)=G[x_1*(g)*x_2]$, or $g(x)=G[x_1(g)x_2]$
- for $g(x):[x_1, x_2]$, with $g(x)=g_1(x)$ if $x \in (x_1, x_2)$, $g(x_1)=a$, $g(x_2)=b$, we will use:
$g(x)=G[a<x_1>g_1<x_2>b]$, $g(x)= G[(a)*x_1*(g_1)*x_2*(b)]$, or $g(x)= G[(a)x_1(g_1)x_2(b)]$
- for $g(x):(x_1, x_2) \cup (x_2, x_3)$, with $g(x)=g_1(x)$ if $x \in (x_1, x_2)$ and $g(x)=g_2(x)$ if $x \in (x_2, x_3)$ :
$g(x)=G[x_1>g_1<x_2>g_2<x_3]$, etc.
- for $g(x):[x_1, x_2) \cup (x_3, x_4]$, with $g(x)=g_1(x)$ if $x \in (x_1, x_2)$ and $g(x)=g_2(x)$ if $x \in (x_3, x_4)$:
$g(x)=G[x_1^+>g_1<x_2> \cup <x_3>g_2<x_4^-]$, or $g(x)=G[x_1 \geq g_1<x_2> \cup <x_3>g_2 \leq x_4]$
- for $g_a(x):[x_1, x_{m+1}]$, with $g_a(x)=g_i(x)$ if $x \in (x_i, x_{i+1})$, and $g_a(x_i)=a_i \neq \infty$, $i=1, 2, ..., m+1$:
$g_a(x)=G[a_1<x_1>g_1<x_2>a_2<x_2>g_2< x_3> a_3 ... a_m<x_m>g_m<x_{m+1}>a_{m+1}]$  (a)
- for $g_b(x): \bigcup_{i=1}^{m}(x_i, x_{i+1})$, with $g_b(x)=g_i(x)$ if $x \in (x_i, x_{i+1})$, $i=1, 2, ..., m+1$, and $g_b(x_i)$ undefined:
$g_b(x)= G[x_1>g_1<x_2>g_2< x_3>... <x_m>g_m<x_{m+1}]$.  (b)

For these last two examples, if $g_i(x)$ are the same, for $i=1, 2, ..., m+1$, the two functions (as well as other $g_r(x)$ functions which fulfill this condition, and moreover the condition $g_r(x_i)=r_i \neq \infty$), are *equals almost everywhere*. If they are square-integrable, their Fourier expansions $\overline{g}_a(x)$, $\overline{g}_b(x)$ and $\overline{g}_r(x)$ have the same expression $\overline{g}(x)$ which, at the points of discontinuity, converges towards:

$$\overline{g}(x_1)= \overline{g}(x_{m+1})=\frac{1}{2}\left[\lim_{x \to x_1^+} f(x)+ \lim_{x \to x_{m+1}^-} f(x)\right] \text{ and } \overline{g}(x_i)=\frac{1}{2}\left[\lim_{x \to x_i^+} f(x)+ \lim_{x \to x_i^-} f(x)\right], \text{ for } i \neq 1, m+1$$

We will note $g_F(x)$ a function of type $g_r(x)$, for which $g_r(x_i)= \overline{g}(x_i)$, $i=1, 2, ..., m+1$:
$g_F(x)= G[ \overline{g}(x_1)<x_1>g_1<x_2> \overline{g}(x_2)<x_2>g_2< x_3> \overline{g}(x_3) ... \overline{g}(x_m)<x_m>g_m<x_{m+1}> \overline{g}(x_1)]$  (c)

We can note $g_a(x) \stackrel{F}{=} g_b(x) \stackrel{F}{=} g_r(x) \stackrel{F}{=} g_F(x) \stackrel{F}{=} \lim_{N \to \infty} \overline{g}(x)$, equality almost everywhere. The Fourier expansion of the discontinuous function $g(x)$ is a continuous function which approaches as much as we want of $g_F(x)$. We can assign to the continuous function $\overline{g}(x)$ (which is an approximation as much as we want close to the function $g_F(x)$)), the designation of **Fourier-function**, or F-function. In many situations, including the majority of this paper, the values $r_i$ of the function $g_r(x)$ at the points of discontinuity are not relevant. Consequently, when we analyze functions of this type, without losing the character of generality, but for the sake of simplification of the exposure, we will always consider (except the cases expressly specified) that it is a function of type $g_F(x)$, and we will use the simplest equivalent notation, that of $g_b(x)$, the relation (b):
$g_F(x)=G[x_1>g_1<x_2>g_2< x_3>... <x_m>g_m<x_{m+1}]$

Other examples:
 for the Heaviside function: $H(x)=G[-\infty>0<0>1/2<0>1<\infty]$
 for the Dirac function: $\delta(x) = \Delta[-\infty>0<0>\infty<0>0<\infty]$
 For the Haar function: $w_{hy}(x) = \Psi[0^+>1<1/2^+>-1<1]$

For any real function $g(x)$, where $x \in I = [x_1, x_2]$, which has finites reals values in the codomenium $g(I)$ and for which $g(x_1)=g(x_2)=a$, we can construct, by successive translations, a **periodic extension** on **R**: $g_p(x_R) = \sum_k g_{pk}(x_R, k)$ where, for $\forall k \in \mathbf{Z}$ and $\forall x \in [x_1, x_2]$

$g_{ap}(a_r k)=g_{pk}(x+kT)=g(x_R-kT)=g(x)$ if $x_R \in [x_1+kT, x_2+kT]$ and (1)
$g_{pk}(x_R)=0$, if $x_R \notin [x_1+kT, x_2+kT]$. Here $T=x_2-x_1$,
For each $k$, there is on the real axis, an interval $[x_1+kT, x_2+kT]$, with $T=x_2-x_1$, for which relation (1) is true. If $k=0$, we obtain for $x \in [x_1, x_2]$: $g_{p0}(x, 0)=g(x)$.

For a definition in which the value of $k$ (dependent on $x_R$), appears implicitly, we can use the *floor function*: $E(x)=\lfloor x \rfloor$ = the biggest integer less than or equal to $x$: $E(x) \leq x < E(x)+1$.
For every $x_R \in \mathbf{R}$, we define the function $K(x_R)=E((x_R-x_1)/T)$ (1a)
So, for $\forall x_R \in R$, $\exists x \in [x_1, x_2]$, $x_R=x+KT$. By definition: $g(x_R)=g(x+KT)=g(x)$

For the function $g_p(x)=sin(x)$, defined on the **R** axis, the relation (1) is true, in the form $sin(x_R)=sin(x)$, for $\forall x_R=x+2\pi k$, **implicitly, simultaneously**, for any $k \in \mathbf{Z}$ and for all intervals $[(2k-1)\pi, (2k+1)\pi]$ that correspond to them. For a certain function $g(x):[x_1, x_2]$, outside this interval $g_p(x)$ must be **explicitly** specified, by **successive** translations, for all validity intervals $x_1+kT<x<x_1+(k+1)T$: $g_p(x_R)=g(x_R-kT)=g(x)$, for all $k \in \mathbf{Z}$, or **implicitly, simultaneously**:
$g_p(x_R)= g[x_R-T \cdot E((x_R-x_1)/T)]=g(x)$. (1b)

According to the theory developed by Fourier, the *2L-periodic function* $f_p(x):(-\infty, \infty)$, the extension on the real axis of the square-integrable function $f(x):[-L, L]$, can be decomposed into a sum: $\bar{f}(x)= f_0 + \sum_{n=1}^{\infty}[a_n \cos(\omega_n x)+b_n \sin(\omega_n x)]$ where, for all $n \in \mathbf{N}$, $\omega_n = n\frac{\pi}{L}$. Here, $f_0$ is the average value of the function $f(x)$, on the interval $[-L, L]$, $cos(\omega_n x)$ and $sin(\omega_n x)$ are continuous functions (called unitary even secondary harmonics, respectively odd unitary secondary harmonics), which come from the continuous functions $cos(\omega_0 x)$ and $sin(\omega_0 x)$, for $\omega_0=\pi/L$ (called unitary even fundamental harmonic, respectively unitary odd fundamental harmonic) by multiplying their arguments with a positive natural integer $n \in \mathbf{N}^+$. The zero mean value functions $sin(n\omega_0 x)$ and $cos(n\omega_0 x)$ take over the interval $[-L/n, L/n]$, the same values as the fundamental harmonics takes over the interval $[-L, L]$ and they implicitly satisfy: $sin(n\omega_0 x)=sin[n(\omega_0 x+2\pi k)]$ and $cos(n\omega_0 x)=cos[n(\omega_0 x+2\pi k)]$, for all intervals $[(2k-1)\omega_0/n, (2k+1)\omega_0/n]$, where $k \in \mathbf{Z}$.

Let be the function $g(x)$ of real variable $x \in [-L, L]$, which has finite real values in the codomenium **Im**$(g)$ and check for equality $g(-L)=g(L)$. Similar to sinusoids, from function $g_p(x) \to (-\infty, \infty)$, which is the *2L-periodic* extension on the real axis of the function $g(x)$, we can get for each $n \in \mathbf{N}^+$, by dilation, a function $2L/n$-periodic: $g_n(x)= g_p(nx) = \sum_{k=-\infty}^{\infty} g_{nk}(nx,k):(-\infty, \infty)$, where for any integer $k$, $g_{nk}(nx,k)$ is a function defined over the interval $[(2k-1)L/n, (2k+1)L/n]$. In this interval, $g_{nk}(nx,k)$ takes the same values as the ones what takes $g_{n0}(nx,0)=g(nx)$ over the interval $[-L/n, L/n]$ and $g(x)$ over the interval $[-L, L]$. Explicitly and successively: $g_{nk}(nx,k)=g(nx+2kL)=g_{n0}(nx,0)=g(nx)$ for $x \in [(2k-1)L/n, (2k+1)L/n]$ and
$g_{nk}(nx,k)=0$, for $x \notin [(2k-1)L/n, (2k+1)L/n]$.
The periodicity relationship becomes: $g_n(x)=g_n(x+2kL/n)$, for all $k \in Z$, or implicitly:
$g_n(x)=g_n(x+2L \cdot E(n(x+L)/2L))$. The function $g(x)$ and the functions $g_n(x)$ have over the interval $[-L,L]$, the same mean value $g_0$.

We will call the function $g_n(x)$, restricted to the interval $[-L, L]$, the **g-harmonic of order $n$** of the function $g(x)$ and the function $g_1(x)=g(x)$, the **fundamental g-harmonic**. We are also going to introduce a reduced notation for the g-harmonic of order $n$:
$g_n(x)=G[-L/n<g(nx)>L/n]_n$, $n \in \mathbf{N}^+$. (1c)

These translation and expansion operations are similar to those used to create the wavelet functions $\psi_{nk}(x)$, from a mother function $\Psi(x)$ [3]: $\psi_{nk}(x)=k\Psi[(x-b)/a]$, for $b=2k/n$ and $a=1/n$.

If the function $g(x)$ has on the interval $[x_1, x_{m+1}]$, a finite number $m$ of discontinuities, the function $g_n(x)$ (the g-harmonic of order $n$) will have a number $m \cdot n$ of such discontinuities, which tends to infinity if $n \to \infty$. For this reason, the function $g(x)$ is unsuitable for generating a base for a subspace of functions. But, if the function $g(x)$ is square-integrable (belong to the space $L^2[x_1,x_{m+1}]$), or if it satisfies the Dirichlet conditions, it can be developed into a series [1, 9]:

$$\bar{g}(x) = g_0 + \sum_{n=1}^{\infty}[a_n \cos(\omega_n x) + b_n \sin(\omega_n x)], \text{ where } \omega_n = n\pi/L, \ \forall n \in \mathbf{N}.$$

Here, $a_n = \dfrac{1}{L}\int_{x_1}^{x_{m+1}} g(x)\cos(\omega_n x)dx$ and $b_n = \dfrac{1}{L}\int_{x_1}^{x_{m+1}} g(x)\sin(\omega_n x)dx$; $g_0 = \dfrac{1}{L}\int_{x_1}^{x_{m+1}} g(x)\cos dx$

Because $\bar{g}(x)$ is a convergent series of continuous functions, it is a continuous function (a Fourier-function $\bar{g}_F(x)$) and can be taken into account to generate a basis for the functions of the space $L^2[x_1, x_{m+1}]$. Obviously, all continuous functions are F-functions. At all points of continuity, $\bar{g}(x) \to g_F(x)$, $\bar{g}(x) \to g_F(x)$, while in the vicinity of a point of discontinuity, $\bar{g}(x-h) \to \lim_{x \to x_i^-} g_F(x)$, and $\bar{g}(x+h) \to \lim_{x \to x_i^+} g_F(x)$, if $h \to 0$. On the interval $[x_i-h, x_i+h]$, for $h \to 0$, the function $\bar{g}(x)$ approach the line $g_F(x)=x[g_F(x_i+h)+g_F(x_i-h)]/2h$ and $\bar{g}(x_i)$ approach the value $[g_F(x_i+h)+g_F(x_i-h)]/2$. Consequently, all the g-harmonics $g_n(x)$, fundamental or secondary are continuous functions throughout the interval $I$.

An F-function can be constructed by definition: let be the function $g(x)$, defined in the interval $[x_1, x_2]$, with a jump discontinuity at the point $x_i$. The corresponding F-function is:

$$g_F(x) = \lim_{h \to 0} G\left[ x_1 * (g) * x_d - h * \left( x \frac{g(x_d - h) + g(x_d + h)}{2h} \right) * x_d + h * (g) * x_2 \right], h \text{ real}. \qquad (1d)$$

As we have already mentioned, in this paper, when we analyze the g-harmonics $g_n(x)$, continuous by pieces, we will always consider (except the cases expressly specified) that they are the Fourier functions $\bar{g}_F(x)$.

For the phenomena of Nature, described by the evolution of certain functions, at least for energy considerations, the discontinuous functions give way to the functions which approach the Fourier functions.
It's obvious that the functions $\bar{g}_n(x)$, $n=1, 2, ..., \infty$ are, two by two, **independents**. Consequently, they form a generating base of a subspace of $L^2$. We will call this base: **the base generated by $g(x)$** or, more simply, **the base $g(x)$**, denoted $B_g$.

## 2. Non-sinusoidal periodic Fourier series

In the previous section, we noted the existence of some formal analogy between the real finite functions $\cos(\omega_0 x)$, $\sin(\omega_0 x)$, defined on the interval $[-L, L]$, and the other real finite functions $g(x)$, defined on the same interval. In this section, we will try to discover those categories of functions $g(x)$ which accentuate this analogy, so that it becomes a functional analogy, useful for creating complete bases of independent functions.

We will use the notations $\bar{f}$, $\hat{f}$ and $\widetilde{f}$, for the expansions in sinusoidal Fourier series, in non-sinusoidal Fourier series and respectively, in orthogonal non-sinusoidal Fourier series.

For the formulas for the expansions in Fourier series and their properties, we have consulted renowned works [4-13].

### 2.1. Non-sinusoidal periodic Fourier series of even functions

**Theorem 1.** *The base $B_g$ of a pair function $g(x)$ defined on the interval $[−L, L]$ of the $L^2$ space (denoted $L^2[−L, L]$), having the mean value $g_0$ null on this interval, constitutes a complete basis for the $F_E$ system of all the even functions $f_e(x)$, real, of $L^2$-space, periodic of period $2L$, having the mean value zero on this interval.*

The proof of this theorem also includes, how to calculate the coefficients $A_n$ of the expansion in non-sinusoidal Fourier series of the even function $f_e(x)$ de $L^2$:

$$\hat{f}_e(x) = \sum_{n=1}^{\infty} A_n \bar{g}_n(x), \quad \text{where } \bar{g}_n(x) \text{ are Fourier series} \tag{2}$$

The function $f_e(x)$ which is, by definition, of zero mean value over the interval $[−L, L]$, can be developed, according to Fourier's thesis, unequivocally, into an infinite sum of even cosine functions:

$$\bar{f}_e(x) = \sum_{n=1}^{\infty} a_n \cos(\omega_n x), \quad \text{where } \omega_n = n\frac{\pi}{L} = n\omega_0. \tag{2.1}$$

At the same time, all quasi-harmonics $g_n(x)$ can be written as a linear combination of the function $\cos\omega_n x$ and the other cosines of higher rank:

$$\bar{g}_1(x) = c_1 \cos\omega_0 x + c_2 \cos 2\omega_0 x + c_3 \cos 3\omega_0 x + c_4 \cos 4\omega_0 x + ...$$
$$\bar{g}_2(x) = c_1 \cos 2\omega_0 x + c_2 \cos 4\omega_0 x + c_3 \cos 6\omega_0 x + c_4 \cos 8\omega_0 x + ...$$
$$\bar{g}_3(x) = c_1 \cos 3\omega_0 x + c_2 \cos 6\omega_0 x + c_3 \cos 9\omega_0 x + c_4 \cos 12\omega_0 x + ...$$
........................
$$\bar{g}_n(x) = c_1 \cos n\omega_0 x + c_2 \cos 2n\omega_0 x + c_3 \cos 3n\omega_0 x + c_4 \cos nN\omega_0 x + ...$$
........................

From these relations, for the general case, we obtain, for $c_1 \neq 0$:

$$\cos\omega_n x = (\bar{g}_n - c_2 \cos 2\omega_n x - c_3 \cos 3\omega_n x - c_4 \cos 4\omega_n x - ....)/c_1, \text{ for } n=1, 2, ..., \infty$$

Here, all $\bar{g}_n(x)$ functions are F-functions (therefore, continuous)

$$\bar{f}_e(x) = a_1 \cos\omega_0 x + a_2 \cos 2\omega_0 x + a_3 \cos 3\omega_0 x + a_4 \cos 4\omega_0 x + ... + a_i \cos i\omega_0 x + ... + a_n \cos n\omega_0 x + ... =$$

$$= \frac{a_1}{c_1}(\bar{g}_1 - c_2 \cos 2\omega_0 x - c_3 \cos 3\omega_0 x - ...) + \frac{a_2}{c_1}(\bar{g}_2 - c_2 \cos 4\omega_0 x - c_3 \cos 6\omega_0 x - ...) +$$

$$+ \frac{a_3}{c_1}(\bar{g}_3 - c_2 \cos 6\omega_0 x - c_3 \cos 9\omega_0 x - ...) + \frac{a_4}{c_1}(\bar{g}_4 - c_2 \cos 8\omega_0 x - c_3 \cos 12\omega_0 x - ...) + ... +$$

$$+ \frac{a_n}{c_1}(\bar{g}_n - c_2 \cos 2n\omega_0 x - c_3 \cos 3n\omega_0 x - ...) =$$

$$= \hat{f}_e(x) = A_1 \bar{g}_1(x) + A_2 \bar{g}_2(x) + A_3 \bar{g}_3(x) + A_4 \bar{g}_4(x) + A_5 \bar{g}_5(x) + ... = \sum_{n=1}^{\infty} A_n \bar{g}_n(x)$$

The equality of $\bar{f}_e(x)$ with $\hat{f}_e(x)$ is unequivocal, which had to be proven. So: (2.2)

$$A_1 = \frac{a_1}{c_1} = K_1, \quad A_2 = K_1\left(\frac{a_2}{a_1} - \frac{c_2}{c_1}\right), \quad A_3 = K_1\left(\frac{a_3}{a_1} - \frac{c_3}{c_1}\right), \quad A_4 = K_1\left(\frac{a_4}{a_1} - \frac{a_2}{a_1}\frac{c_2}{c_1} - \frac{c_4}{c_1} + \frac{c_2^2}{c_1^2}\right), \quad A_5 = K_1\left(\frac{a_5}{a_1} - \frac{c_5}{c_1}\right)$$

$$A_6 = K_1\left(\frac{a_6}{a_1} - \frac{a_2}{a_1}\frac{c_3}{c_1} - \frac{a_3}{a_1}\frac{c_2}{c_1} - \frac{c_6}{c_1} + 2\frac{c_2 c_3}{c_1^2}\right), \quad A_7 = K_1\left(\frac{a_7}{a_1} - \frac{c_7}{c_1}\right),$$

$$A_8 = K_1\left(\frac{a_8}{a_1} - \frac{a_2}{a_1}\frac{c_4}{c_1} - \frac{a_4}{a_1}\frac{c_2}{c_1} - \frac{c_8}{c_1} + 2\frac{c_2c_4}{c_1^2} + \frac{a_2}{a_1}\frac{c_2^2}{c_1^2} - \frac{c_2^3}{c_1^3}\right), \quad A_9 = K_1\left(\frac{a_9}{a_1} - \frac{a_3}{a_1}\frac{c_3}{c_1} - \frac{c_9}{c_1} + \frac{c_3^2}{c_1^2}\right),$$

$$A_{10} = K_1\left(\frac{a_{10}}{a_1} - \frac{a_2}{a_1}\frac{c_5}{c_1} - \frac{a_5}{a_1}\frac{c_2}{c_1} - \frac{c_{10}}{c_1} + \frac{2c_2c_5}{c_1^2}\right), \quad A_{11} = \frac{a_{11}c_1 - a_1c_{11}}{c_1^2} = K_1\left(\frac{a_{11}}{a_1} - \frac{c_{11}}{c_1}\right),$$

$$A_{12} = \frac{a_{12}c_1 - a_2c_6 - a_3c_4 - a_4c_3 - a_6c_2 - a_1c_{12}}{c_1^2} + \frac{2a_1c_3c_4 + 2a_1c_2c_6 + 2a_2c_2c_3 + a_3c_2^2}{c_1^3} - \frac{3a_1c_2^2c_3}{c_1^4} =$$

$$= K_1\left(\frac{a_{12}}{a_1} - \frac{a_2}{a_1}\frac{c_6}{c_1} - \frac{a_3}{a_1}\frac{c_4}{c_1} - \frac{a_4}{a_1}\frac{c_3}{c_1} - \frac{a_6}{a_1}\frac{c_2}{c_1} - \frac{c_{12}}{c_1} + 2\frac{c_3c_4 + c_2c_6}{c_1^2} + 2\frac{a_2}{a_1}\frac{c_2c_3}{c_1^2} + \frac{a_3}{a_1}\frac{c_2^2}{c_1^2} - 3\frac{c_2^2c_3}{c_1^3}\right), \ldots$$

In conclusion, because any even function $f_e(x)$ of the subspace $L^2[-L, L]$, can be developed in a sinusoidal Fourier series (2.1), it can also be developed in a non-sinusoidal Fourier series (2). To calculate the coefficients of this expansion, it is necessary to know the coefficients $a_n$ of the Fourier expansion of the function $f_e(x)$ as well as $c_n$, the coefficients of the function $g(x)$, which implies the calculation of the integrals $\int_{-L}^{L} f_e(x)\cos\omega_n x\, dx$, respectively $\int_{-L}^{L} g(x)\cos\omega_n x\, dx$. Of course, for another even function $f(x)=f_0(x)+f_e(x)$, $f_0(x)\neq 0$,

$\hat{f}(x) = f_0(x) + \sum_{n=1}^{\infty} A_n \bar{g}_n(x)$. The approximation of order $N$ is written:

$\hat{f}_N(x) = f_0(x) + \sum_{n=1}^{N} A_n \bar{g}_n(x)$, where $\bar{g}_n(x) = \sum_{m=1}^{N} c_m \cos(m\omega_n x)$. If $N\to\infty$, $\bar{g}_n(x) \to g(x)$

To illustrate the calculation method, let be the function $f_2(x)=G[-1^+>x^2<1^-]$, which is a second degree polynomial, with no discontinuities, and which has the Fourier expansion:

$$\bar{f}_2(x) = f_0 + \sum_{n=1}^{\infty} a_n \cos n\pi x = \frac{1}{3} + \sum_{n=1}^{\infty} \frac{4(-1)^n}{n^2\pi^2} \cos n\pi x \tag{2a}$$

We want to develop it into a base generated by the even function (rectangular pulses of zero mean value) $g_e = g^{dr} = G_e[-1>-1<-1/2>1<1/2>-1<1]$:

$$\hat{f}_2(x) = \frac{1}{3} + \sum_{n=1}^{\infty} A_n \bar{g}_n^{dr}(x), \tag{2b}$$

The expansion in trigonometric series of the function $g_e(x)$ is:

$$\bar{g}_e(x) = \sum_{n=1}^{\infty} c_n \cos(2n-1)\pi x = -\frac{2}{\pi}\sum_{n=1}^{\infty}\frac{(-1)^n \cos(2n-1)\pi x}{2n-1} \tag{2c}$$

The relations (2a) and (2c), provide the following coefficients:
$a_1=-4/\pi^2$, $a_2=1/\pi^2$, $a_3=-4/9\pi^2$, $a_4=1/4\pi^2$, $a_5=-4/25\pi^2$, $a_6=1/9\pi^2$, $a_7=-4/49\pi^2$, $a_8=1/16\pi^2$, $a_9=-4/81\pi^2$, $a_{10}=1/25\pi_2$, $a_{11}=-4/121\pi^2$, $a_{12}=1/36\,\pi^2$, ... and
$c_1=2/\pi$, $c_2=0$, $c_3=-2/3\pi$, $c_4=0$, $c_5=2/5\pi$, $c_6=0$, $c_7=-2/7\pi$, $c_8=0$, $c_9=2/9\pi$, $c_{10}=0$, $c_{11}=-2/11\pi$, $c_{12}=0$, ...

According to (2.2), the expansion coefficients (2b) are:
$A_1=-2/\pi$, $A_2=1/2\pi$, $A_3=-8/9\pi$, $A_4=1/8\pi$, $A_5=8/25\pi$, $A_6=2/9\pi$, $A_7=-16/49\pi$, $A_8=1/32\pi$, $A_9=-8/81\pi$, $A_{10}=-2/25\pi$, $A_{11}=-24/121\pi$, $A_{12}=1/18\,\pi$, ...

The representation of the corresponding quasi-harmonics and the resulting partial sums is given in figure 1. Here, we have represented the functions $g_{en}(x)$ instead of the functions $\bar{g}_{en}(x)$. Since the function $g_e(x)$ has two points of discontinuity, the partial sums $S_N$ of the non-sinusoidal expansion, have jump points in increasing number, as the rank $N$ increases. One can notice a low speed of convergence, compared to the traditional method of approximation of the same curve, by horizontal line segments

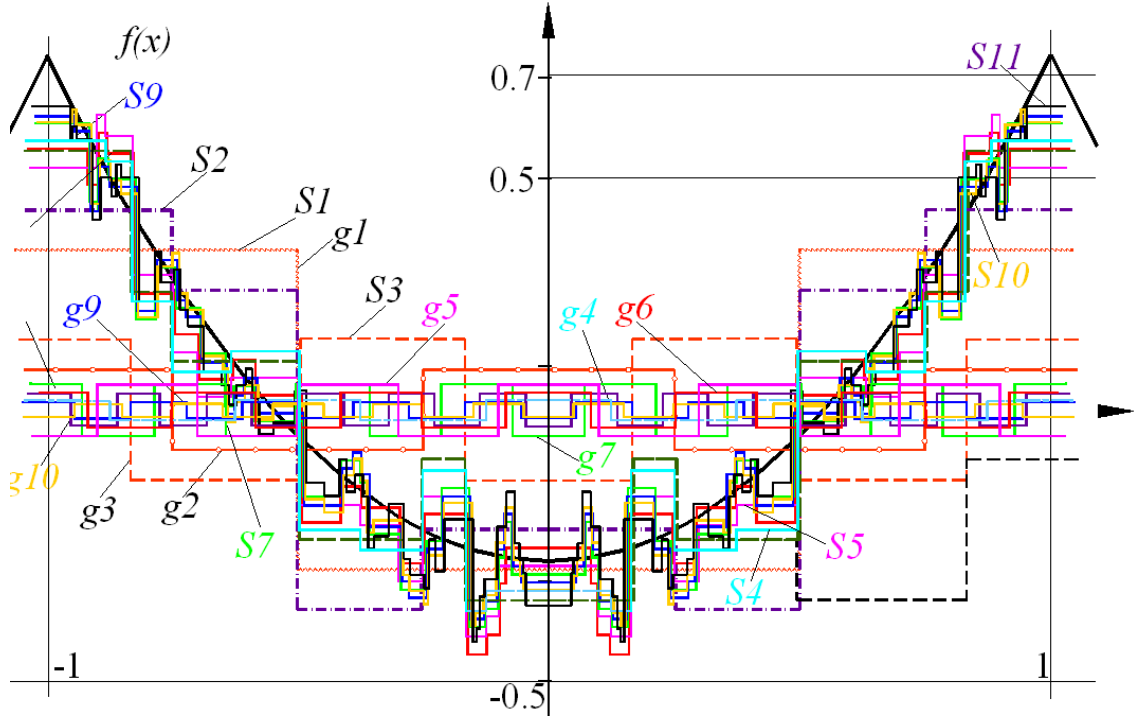

*Fig. 1. The approximation of the function $x^2-1/3$ by a sum of rectangular pulses*
*$g_i$= component i of partial sums, $S_i$= partial sum of order i*

Let us examine, for a comparison, the expansion of the same function by a sum of the triangle waves of zero mean value. Let be the same function: $f_2(x)=F_2[-1^+>x^2<1^-]$ and the function $g(x)=G[-1^+>x-1/2<0>x-1/2<1^-]$.

$\bar{f}_2(x) = \frac{1}{3} + \sum_{n=1}^{\infty} \frac{4(-1)^n}{n^2\pi^2} \cos n\pi x$ provides:

$a_1=-4/\pi^2$, $a_2=1/\pi^2$, $a_3=-4/9\pi^2$, $a_4=1/4\pi^2$, $a_5=-4/25\pi^2$, $a_6=1/9\pi^2$, $a_7=-4/49\pi^2$, $a_8=1/16\pi^2$, $a_9=-4/81\pi^2$, $a_{10}=1/25\pi_2$, $a_{11}=-4/121\pi^2$, $a_{12}=1/36\ \pi^2$, ...

and $\bar{g}(x) = -4\sum_{n=1}^{\infty} \frac{\cos(2n-1)\pi x}{(2n-1)^2 \pi^2}$ provides: (2d)

$c_1=-4/\pi^2$, $c_2=0$, $c_3=-4/9\pi^2$, $c_4=0$, $c_5=-4/25\pi^2$, $c_6=0$, $c_7=-4/49\pi^2$, $c_8=0$, $c_9=-4/81\pi^2$, $c_{10}=0$, $c_{11}=-4/121\pi^2$, $c_{12}=0$, ...

Using relations (2.2), we can develop the function $f_2(x)$ into an infinite series of triangle-functions: $\hat{f}_2(x) = \frac{1}{3} + \sum_{n=1}^{\infty} A_n \bar{g}_n(x)$, where:

$A_1=1$, $A_2=-1/4$, $A_3=0$, $A_4=-1/16$, $A_5=0$, $A_6=0$, $A_7=0$, $A_8=-1/64$, $A_9=0$, $A_{10}=0$, $A_{11}=0$, $A_{12}=0$, ...

$\hat{f}_2(x) = \frac{1}{3} + \bar{g}_1(x) - \sum_{n=1}^{\infty} 4^{-n} \bar{g}_{2^n}(x) = \frac{1}{3} + \left[-1 > -x - \frac{1}{2} < 0 > x - \frac{1}{2} < 1\right] -$

$- \sum_{n=1}^{\infty} \frac{1}{4^n} \left[-\frac{1}{2^n} > -2^n x - \frac{1}{2} < 0 > 2^n x - \frac{1}{2} < \frac{1}{2^n}\right]_{2^n}$

Figure 2 shows the first quasi-harmonics and the first partial sums.

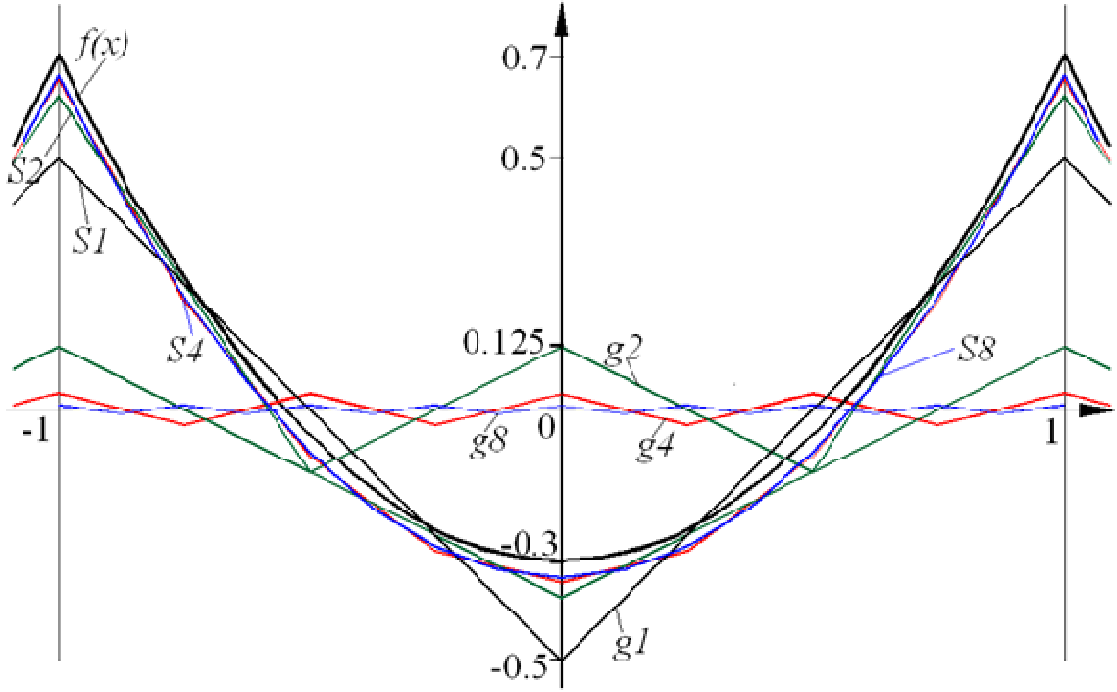

*Fig. 2. The approximation of the function $x^2 - 1/3$ by a sum of even triangle functions*
*gi: quasi-harmonics of order i; Si: partial sums*

The functions $f_e(x)$ and $g(x)$ being both pairs, we can also do the expansion in the opposite direction: $\hat{g}(x) = \sum_{n=1}^{\infty} A'_n \bar{f}_{en}(x)$, for $a_1 \neq 0$, where: (2.3)

$$A'_1 = \frac{c_1}{a_1},\ A'_2 = \frac{c_2 a_1 - c_1 a_2}{a_1^2},\ A'_3 = \frac{c_3 a_1 - c_1 a_3}{a_1^2},\ A'_4 = \frac{c_4 a_1 - c_2 a_2 - c_1 a_4}{a_1^2} + \frac{c_1 a_2^2}{a_1^3},\ A'_5 = \frac{c_5 a_1 - c_1 a_5}{a_1^2},\ \ldots$$

If the function $f_e(x)$ is even the function $f = \cos\omega_0 x$, we write:

$$\cos\omega_0 x = \sum_{n=1}^{\infty} A_n \bar{g}_n(x) = A_1 \bar{g}_1(x) + A_2 \bar{g}_2(x) + A_3 \bar{g}_3(x) + \ldots =$$
$$= A_1(c_1\cos\omega_0 x + c_2\cos 2\omega_0 x + c_3\cos 3\omega_0 x + \ldots) + A_2(c_1\cos 2\omega_0 x + c_2\cos 4\omega_0 x + c_3\cos 6\omega_0 x + \ldots) +$$
$$+ A_3(c_1\cos 3\omega_0 x + c_2\cos 6\omega_0 x + c_3\cos 9\omega_0 x + \ldots) + A_4(c_1\cos 4\omega_0 x + c_2\cos 8\omega_0 x + c_3\cos 12\omega_0 x + \ldots) + \ldots$$

This relation leads to the following system of coefficients: (2.4)

$$A_1 = \frac{1}{c_1},\ A_2 = -\frac{c_2}{c_1^2},\ A_3 = -\frac{c_3}{c_1^2},\ A_4 = \frac{-c_1 c_4 + c_2^2}{c_1^3},\ A_5 = -\frac{c_5}{c_1^2},\ A_6 = \frac{-c_1 c_6 + 2c_2 c_3}{c_1^3},\ A_7 = -\frac{c_7}{c_1^2},$$

...

If the pair function $g(x)$ that generates the basis of the expansion is a function with the mean value $g_0 \neq 0$ over the interval $[-L, L]$:

$$\cos\omega_n x = (\bar{g}_n - g_0 - c_2\cos 2\omega_n x - c_3\cos 3\omega_n x - c_4\cos 4\omega_n x - \ldots)/c_1,\ \ \text{for } n \in N$$

and if the function $f(x)$ also has an average value $f_0 \neq 0$ over this interval:

$$\hat{f}(x) = f_0 + \sum_{n=1}^{\infty} A_n[\bar{g}_n(x) - g_0] \tag{2.5}$$

## 2.2. Non sinusoidal periodic Fourier series of odd functions

The problem of the odd function $f_o(x)$, which belongs to the $F_O$ space of the odd functions of $L^2[-L, L]$ is treated in the same way.

**Theorem 2.** *The basis $B_h$ generated by any odd function $h(x)$ of $L^2[-L, L]$ is a complete basis for the $F_O$ system of all odd functions $f_o(x)$, real, of $L^2$-space, of $2L$-period.*

According to Fourier's thesis, the odd function $f_o(x)$ (whose average value over the interval $[-L, L]$ is always zero) can be unambiguously developed into an infinite sum of odd sinusoidal functions:

$$\bar{f}_o(x) = \sum_{n=1}^{\infty} b_n \sin(\omega_n x), \text{ where } \omega_n = n\omega_0 = n\frac{\pi}{L}. \tag{2.6}$$

Any other expansion of the function $f_o(x)$ must also be an infinite sum of odd functions:

$$\hat{f}_o(x) = \sum_{n=1}^{\infty} B_n \bar{h}_n(x) \text{ , where } \bar{h}_n(x) = G[-1/n < h(nx) > 1/n]_n \text{ , } n \in \mathbf{N}, \tag{2.7}$$

$\bar{h}_n(x)$ are F-functions $2L/n$-periodic. For each such function: $\bar{h}_n(x) = \sum_{i=1}^{\infty} d_i \sin(i\omega_n x)$

As in the previous demonstration, this system of equations allows us to determine the $B_n$ coefficients of expansion: (2.8)

$$B_1 = \frac{b_1}{d_1} = K_2, \ B_2 = K_2\left(\frac{b_2}{b_1} - \frac{d_2}{d_1}\right), \ B_3 = K_2\left(\frac{b_3}{b_1} - \frac{d_3}{d_1}\right), \ B_4 = K_2\left(\frac{b_4}{b_1} - \frac{b_2 d_2}{b_1 d_1} - \frac{d_4}{d_1} + \frac{d_2^2}{d_1^2}\right),$$

$$B_5 = K_2\left(\frac{b_5}{b_1} - \frac{d_5}{d_1}\right), \ B_6 = K_2\left(\frac{b_6}{b_1} - \frac{b_2 d_3}{b_1 d_1} - \frac{b_3 d_2}{b_1 d_1} - \frac{d_6}{d_1} + 2\frac{d_2 d_3}{d_1^2}\right), \ B_7 = K_2\left(\frac{b_7}{b_1} - \frac{d_7}{d_1}\right), \text{ etc.}$$

We can conclude that any odd function $f_o(x)$ of the subspace $L^2[-L, L]$, which can be developed in sinusoidal Fourier series (2.6), can also be developed in non-sinusoidal Fourier series (2.7). To calculate the coefficients of this expansion, (as well as those of the inverse expansion), it is necessary to know the coefficients $b_n$ of the Fourier expansion of the function $f_o(x)$, as well as those of the function $h(x)$, which involves the calculation of integrals $\int_{-L}^{L} f_o(x) \sin \omega_n x \, dx$, respectively $\int_{-L}^{L} h(x) \sin \omega_n x \, dx$.

To illustrate the calculation method, let be the odd functions
$f_o(x) = F_o[-1 > -1 < 0 > 1 < 1]$ (the odd rectangular pulses), for which:

$$\bar{f}_o(x) = \frac{4}{\pi} \sum_{n=1}^{\infty} \frac{\sin(2n-1)\pi x}{2n-1} \text{ , and } g_o(x) = G[-1 > x < 1] \text{ (the sawtooth wave), for which:}$$

$$\bar{g}_o(x) = \sum_{n=1}^{\infty} d_n \sin n\pi x = \frac{2}{\pi} \sum_{n=1}^{\infty} (-1)^{n+1} \frac{\sin n\pi x}{n} \text{ , from where:}$$

$b_1=4/\pi$, $b_2=0$, $b_3=4/3\pi$, $b_4=0$, $b_5=4/5\pi$, $b_6=0$, $b_7=4/7\pi$, $b_8=0$, $b_9=4/9\pi$, $b_{10}=0$, $b_{11}=4/10\pi$, ...,
$d_1=2/\pi$, $d_2=-2/2\pi$, $d_3=2/3\pi$, $d_4=-2/4\pi$, $d_5=2/5\pi$, $d_6=-2/6\pi$, $d_7=2/7\pi$, $d_8=-2/8\pi$, $d_9=2/9\pi$, $d_{10}=-2/10\pi$, $d_{11}=2/11\pi$, $d_{12}=-2/12\pi$, ...

Consequently, the expansion $\hat{f}_o(x) = \sum_{n=1}^{\infty} B_n \bar{g}_n(x)$ will have the following coefficients:

$B_1 = 2$, $B_2 = 1$, $B_3 = 0$, $B_4 = 1$, $B_5 = 0$, $B_6 = 0$, $B_7 = 0$, $B_8 = 1/2$, $B_9 = 0$, $B_{10} = 1/5$,
$B_{11} = 0$, $B_{12} = 0$, $B_{13}=0$, $B_{14}=0$, $B_{15}= -2/15$, ...,

Figure 3 presents the graphs of these two functions (**3a** and **3b**), the first quasi-harmonics of the expansion (**c**) and the first partial sums (**d**). We note that, as $N \rightarrow \infty$, the sum $S_N(x)$ tends very slowly towards the function $f_o(x)$.

The coefficients of the inverse expansions are (for $b_1 \neq 0$): (2.9)

$$B'_1 = \frac{d_1}{b_1},\ B'_2 = \frac{d_2 b_1 - b_1 d_2}{b_1^2},\ B'_3 = \frac{d_3 b_1 - d_1 b_3}{b_1^2},\ B'_4 = \frac{d_4 b_1 - d_2 b_2 - d_1 b_4}{b_1^2} + \frac{d_1 b_2^2}{b_1^3},\ B'_5 = \frac{d_5 b_1 - d_1 b_5}{b_1^2},$$

For the two previous functions $f_o(x)=F_o[-1>-1<0>1<1]$ et $g_o(x)=G[-1>x<1]$:
$B'_1 = 1/2$, $B'_2 = 1/4$, $B'_3 = 0$, $B'_4 = 1/8$, $B'_5 = 0$, $B'_6 = 0$, $B'_7 = 0$, $B'_8 = 1/16$, $B'_9 = 0$, $B'_{10} = 0$, $B'_{11} = 0$, $B'_{12} = 0$, $B'_{13} = 0$, $B'_{14} = 0$, $B'_{15} = 0$, $B'_{16} = 1/32$, ...,
and we can write:

$$\overline{g}_o(x) = \overline{G}[-1<x>1] = \frac{1}{2}\left(\overline{F}_{o1} - \sum_{n=1}^{\infty}\frac{1}{2^n}\overline{F}_{o2^n}\right) = \frac{L}{2}\overline{F}_o[-1>-1<0>1<1] - \frac{1}{2}\sum_{n=1}^{\infty}\frac{1}{2^n}\left[-\frac{1}{2^n}>-1<0>1<\frac{1}{2^n}\right]_{2^n}^{F}$$

If $f(x)=\sin(\omega_0 x)$, the coefficients of a non-sinusoidal expansion are: (2.10)

$$B_1 = \frac{1}{d_1},\ B_2 = -\frac{d_2}{d_1^2},\ B_3 = -\frac{d_3}{d_1^2},\ B_4 = \frac{-d_1 d_4 + d_2^2}{d_1^3},\ B_5 = -\frac{d_5}{d_1^2},\ B_6 = \frac{-d_1 d_6 + 2 d_2 d_3}{d_1^3},\ B_7 = -\frac{d_7}{d_1^2},$$

Where $d_n$ are the expansion coefficients $\overline{h}(x) = \sum_{n=1}^{\infty} d_n \sin(\omega_n x)$, where $\omega_n = n\frac{\pi}{L}$.

For the function $\Pi_1(\theta) = \Pi[-\pi>-1<0>1<\pi]$, for which $\Pi_n(\theta) = \Pi[-\pi/n>-1<0>1<\pi/n]_n$, and $\overline{\Pi}_1 = \frac{4}{\pi}\sum_{n=1}^{\infty}\frac{\sin(2n-1)\theta}{2n-1}$, to write the expansion $\sin\theta = \hat{f}(\theta) = \sum_{n=1}^{\infty} B_n \overline{\Pi}_n(\theta)$, we use the formulas (2.10) and we obtain:

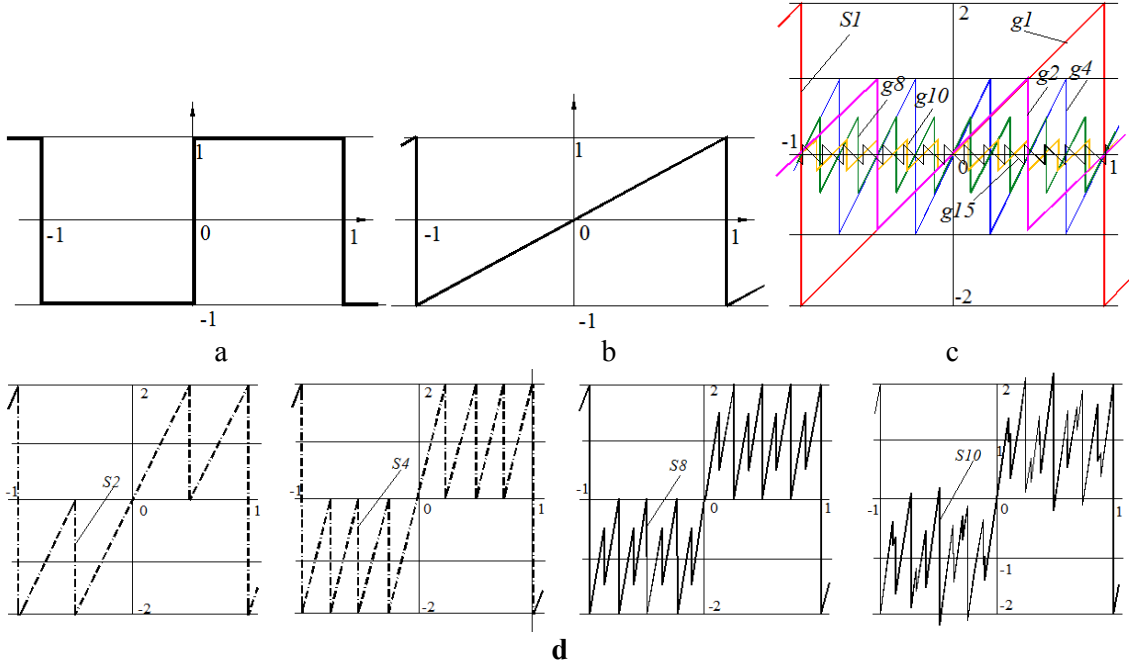

Fig. 3. The expansion of the function $f_o(x)$ into the base $g_o(x)$
*a*: $f_o(x)$ *b*: $g_o(x)$ *c*: quasi-harmonics $g_1(x)$, $g_2(x)$, $g_4(x)$, $g_8(x)$, $g_{10}(x)$, $g_{15}(x)$
*d*: the first partial sums

$$B_1 = \frac{\pi}{4},\ B_2 = 0,\ B_3 = -\frac{\pi}{12},\ B_4 = 0,\ B_5 = -\frac{\pi}{20},\ B_6 = 0,\ B_7 = -\frac{\pi}{28},\ B_8 = 0,\ B_9 = 0,\ B_{10} = 0,\ ...$$

$$B_{2n-1} = \frac{\pi}{4} \cdot \frac{-1}{2n-1}\ ,\ B_{2n} = 0,\ \text{but } B_{n^2} = 0,\ \text{pour } n = 2, 3, .., \infty$$

Tous les coefficients du développement, sauf le coefficient de la fondamentale sont négatifs.
Pour la fonction $g_1(\theta)=X_o[-\pi<-\theta-\pi>-\pi/2<\theta>\pi/2<-\theta+\pi>\pi]$, le développement:

$$\bar{g}_1(\theta) = \frac{4}{\pi}\sum_{n=1}^{\infty}\frac{(-1)^{n+1}\sin(2n-1)\theta}{(2n-1)^2}$$, fournit les coefficients:

$d_1=4/\pi$, $d_2=0$, $d_3=-4/9\pi$, $d_4=0$, $d_5=4/25\pi$, $d_6=0$, $d_7=-4/49\pi$, $d_8=0$, $d_9=4/81\pi$, $d_{10}=0$, $d_{11}=-4/121\pi$, $d_{12}=0$, ..., pour lesquelles, pour $f(\theta)=\sin\theta$, on obtient:

$$B_1 = \frac{\pi}{4},\ B_2 = 0,\ B_3 = \frac{\pi}{36},\ B_4 = 0,\ B_5 = -\frac{\pi}{100},\ B_6 = 0,\ B_7 = \frac{\pi}{196},\ B_8 = 0,\ B_9 = 0,\ ...$$

(2.11)

$$B_{2n-1} = \frac{\pi}{4}\cdot\frac{(-1)^n}{(2n-1)^2}\ ,\ B_{2n} = 0,\ \text{but}\ B_{n^2} = 0,\ \text{for}\ n=2,\ 3,\ ..,\ \infty$$

### 2.3. Non-sinusoidal periodic Fourier series of any function

In the general case, a certain function $f(x)$ of $L^2[-L, L]$-space, can be written as the sum of its mean value $f_0$ over this interval, of its even component $f_e(x)$ (by definition, of zero mean value over the interval $[-L, L]$) and of its odd component $f_o(x)$: $f(x)=f_0+f_e(x)+f_o(x)$. Following the two previous theorems, we can say:

**Theorem 3.** *Any function $f(x)$ of $L^2[-L, L]$, can be developed in non-sinusoidal Fourier series into a base composed of $f_0$ (its mean value over this interval), of a certain even base $B_{g-g0}$ and of a certain odd base $B_h$ of $L^2[-L, L]$, where the even function $g(x)$ and the odd function $h(x)$, are any functions of $L^2[-L, L]$:*

$$\hat{f}(x) = f_0 + \sum_{n=1}^{\infty}A_n[\bar{g}_n(x)-g_0] + \sum_{n=1}^{\infty}B_n\bar{h}_n(x),\ \text{where}\ g_0 = \int_{-L}^{L}g(x)dx$$

We can see that the sinusoidal Fourier expansion is a special case of the non-sinusoidal Fourier expansion. To illustrate, let be the function $f(x)=F[-1>0<-1/2>-2<0>0<1/2>2<1]$ which is the sum of $f_0$ (=0), of the even function $f_e=F_e[-1>1<-1/2>-1<1/2>1<1]$ and of the odd function $f_o=F_o[-1>-1<0>1<1]$, whose Fourier expansions are [5]:

$$\bar{f}_e(x) = \sum_{n=1}^{\infty}a_n\cos(n\pi x) = -\frac{4}{\pi}\sum_{n=1}^{\infty}\frac{(-1)^{n+1}\cos(2n-1)\pi x}{2n-1}\ ,\ \text{respectively}$$

$$\bar{f}_o(x) = \sum_{n=1}^{\infty}b_n\sin(n\pi x) = \frac{4}{\pi}\sum_{n=1}^{\infty}\frac{\sin(2n-1)\pi x}{2n-1}$$

The following coefficients are obtained:

$a_1=4/\pi$, $a_2=0$, $a_3=-4/3\pi$, $a_4=0$, $a_5=4/5\pi$, $a_6=0$, $a_7=-4/7\pi$, $a_8=0$, $a_9=4/9\pi$, $a_{10}=0$, $a_{11}=-4/11\pi$,...

$b_1=4/\pi$, $b_2=0$, $b_3=4/3\pi$, $b_4=0$, $b_5=4/5\pi$, $b_6=0$, $b_7=4/7\pi$, $b_8=0$, $b_9=4/9\pi$, $b_{10}=0$, $b_{11}=4/11\pi$, ...

For the expansion of the function $f(x)$ in an exponential basis $g(x)=e^x$, given that over the interval $[-1, 1]$, $g_0=\sinh 1$, we will choose the set consisting of the functions $g_e(x)=\cosh x - \sinh 1$ and $g_o(x)=\sinh x$, whose expansions in Fourier series are:

$$\bar{g}_e(x) = \sinh 1\sum_{n=1}^{\infty}\frac{2(-1)^n}{1+n^2\pi^2}\cos(n\pi x),\ \text{and}\ \bar{g}_o(x) = \pi\cdot\sinh 1\sum_{n=1}^{\infty}\frac{2n\cdot(-1)^{n+1}}{1+n^2\pi^2}\sin(n\pi x),\ \text{so:}$$

$c_1=-2\sinh 1/(1+\pi^2)$, $c_2=2\sinh 1/(1+4\pi^2)$, $c_3=-2\sinh 1/(1+9\pi^2)$,
$c_4=2\sinh 1/(1+16\pi^2)$, $c_5=-2\sinh 1/(1+25\pi^2)$
$c_6=2\sinh 1/(1+36\pi^2)$, $c_7=-2\sinh 1/(1+49\pi^2)$, $c_8=2\sinh 1/(1+64\pi^2)$, $c_9=-2\sinh 1/(1+81\pi^2)$,
$c_{10}=2\sinh 1/(1+100\pi^2)$, $c_{11}=-2\sinh 1/(1+121\pi^2)$, $c_{12}=2\sinh 1/(1+144\pi^2)$, ...

$d_1=2\pi\sinh1/(1+\pi^2)$, $d_2=-4\pi\sinh1/(1+4\pi^2)$, $d_3=6\pi\sinh1/(1+9\pi^2)$, $d_4=-8\pi\sinh1/(1+16\pi^2)$,
$d_5=10\pi\sinh1/(1+25\pi^2)$, $d_6=-12\pi\sinh1/(1+36\pi^2)$, $d_7=14\pi\sinh1/(1+49\pi^2)$,
$d_8=-16\pi\sinh1/(1+64\pi^2)$, $d_9=18\pi\sinh1/(1+81\pi^2)$, $d_{10}=-20\pi\sinh1/(1+100\pi^2)$,
$d_{11}=22\pi\sinh1/(1+121\pi^2)$, $d_{12}=-24\pi\cdot\sinh1/(1+144\pi^2)$

This results in a non-sinusoidal Fourier series expansion of the form:

$$\hat{f}(x)=\hat{f}_e(x)+\hat{f}_o(x)=\sum_{n=1}^{\infty}A_n[\cosh_{Fn}(x)-\sinh1]+\sum_{n=1}^{\infty}B_n\sinh_{Fn}(x)$$

where $\cosh_{Fn}$ and $\sinh_{Fn}$ are the extensions on the real axis of the F-functions $\cosh_F(nx)$, respectively $\sinh_F(nx)$, defined on the intervals $[-1/n, 1/n]$), and the coefficients $A_n$, $B_n$ are:

for $K_1 = \dfrac{a_1}{c_1} = -\dfrac{2(1+\pi^2)}{\pi\cdot\sinh1}$ and $K_2 = \dfrac{b_1}{d_1} = \dfrac{2(1+\pi^2)}{\pi^2\cdot\sinh1} = -\dfrac{K_1}{\pi}$ :

$A_1 = K_1$, $A_2 = K_1\dfrac{1+\pi^2}{1+4\pi^2}$, $A_3 = -K_1\dfrac{4}{3}\cdot\dfrac{1+3\pi^2}{1+9\pi^2}$, $A_4 = K_1\left[\dfrac{1+\pi^2}{1+16\pi^2}+\dfrac{(1+\pi^2)^2}{(1+4\pi^2)^2}\right]$, $A_5 = -K_1\dfrac{4}{5}\cdot\dfrac{1-5\pi^2}{1+25\pi^2}$

$A_6 \approx -K_1\dfrac{(1+\pi^2)\cdot(7+16\pi^2)}{36\pi^2(1+4\pi^2)}$, $A_7 = -K_1\dfrac{8}{7}\cdot\dfrac{1+7\pi^2}{1+49\pi^2}$, $A_8 \approx K_1\dfrac{(1+\pi^2)\cdot(7+17\pi^2+12\pi^4)}{4(1+4\pi^2)^3}$,...

$B_1 = K_2$, $B_2 = 2K_2\dfrac{1+\pi^2}{1+4\pi^2}$, $B_3 = -K_2\cdot\dfrac{8}{3}\cdot\dfrac{1}{1+9\pi^2}$, $B_4 = 4K_2\cdot\dfrac{(1+\pi^2)\cdot(2+25\pi^2+32\pi^4)}{(1+16\pi^2)\cdot(1+4\pi^2)^2}$

$B_5 = -K_2\cdot\dfrac{24}{5}\cdot\dfrac{1}{1+25\pi^2}$, $B_6 \approx -K_2\cdot\dfrac{3(1+\pi^2)\cdot(3+2\pi^2)}{(1+4\pi^2)\cdot(1+9\pi^2)}$, $B_7 = -K_2\cdot\dfrac{48}{7}\cdot\dfrac{1}{1+49\pi^2}$,

$B_8 \approx K_2\cdot\dfrac{2(1+\pi^2)\cdot(9+21\pi^2+16\pi^4)}{(1+4\pi^2)^3}$, $B_9 \approx -K_2\cdot\dfrac{1}{9\pi^4}$, $B_{10} \approx -K_2\cdot\dfrac{1+\pi^2}{1+4\pi^2}\cdot\dfrac{3-4\pi^2}{10\pi^2}$, ...

Over the interval $[-1/2, 1/2]$: $g_e(x)_2 = \cosh x - g_{02} = 1/2(e^x + e^{-x}) - 2\sinh(1/2)$ and $g_o(x)_2 = \sinh x$,

so $\bar{g}_e(x)_2 = \sinh\dfrac{1}{2}\cdot\sum_{n=1}^{\infty}\dfrac{4(-1)^n}{1+4n^2\pi^2}\cos(2n\pi x)$, and $\bar{g}_o(x)_2 = \pi\cdot\sinh\dfrac{1}{2}\cdot\sum_{n=1}^{\infty}\dfrac{4n\cdot(-1)^{n+1}}{1+4n^2\pi^2}\sin(2n\pi x)$

We can note that for values $L<1$, the function $g_{es}(x)=g_e(x)_L/\sinh L$ is approximated with acceptable deviations by the function $g_{ep}(x)=x^2$, and the function $g_{os}(x)=g_o(x)_L/\sinh L$ is approximated with acceptable deviations by the function $g_{op}(x)=x$, the deviations being so small that $L$ is smaller.

When we ask for the expansion of the function $f(x)=f_0+f_e+f_o$ into a base generated by any function $g(x)=g_0+g_e(x)+g_o(x)$ of $L^2[-L, L]$, we must find the expansion coefficients $C_n$ of:

$$\hat{f}(x) = f_0 + \sum_{n=0}^{\infty} C_n[g_{Fn}(x)-g_0], \quad \text{where} \quad g_{Fn}(x)=G_F[-L/n<g_F(nx)>L/n]_n \quad , \quad n\in\mathbf{N}.$$

(2.11)

To simplify, consider the particular case $f_0 = g_0 = 0$:

$$\hat{f}(x) = \hat{f}_e(x) + \hat{f}_o(x) = \sum_{n=1}^{\infty} A_n g_{Fen}(x) + \sum_{n=1}^{\infty} B_n g_{Fon}(x) =$$

$$= \sum_{n=1}^{\infty} [A_n g_{Fn}(x) + (B_n - A_n)g_{Fon}(x)] = \sum_{n=1}^{\infty} [(A_n - B_n)g_{Fen}(x) + B_n g_{Fn}(x)]$$

equality which coincides with (2.11), only if $A_n = B_n = C_n$. Therefore, no function $g_F(x)$ alone can generate a base for the entire space $L^2[-L, L]$, requiring the help of another base, generated by a function $h_F(x)$ with a different parity. If we consider the identities:

$g_{en}(x) = \dfrac{1}{2}[g_n(x)+g_n(-x)]$ and $g_{on}(x) = \dfrac{1}{2}[g_n(x)-g_n(-x)]$ , one obtains, in the general case:

$$\hat{f}(x) = f_0 + \sum_{n=1}^{\infty} \left[ \frac{A_n + B_n}{2} [g_{Fn}(x) - g_0] + \frac{A_n - B_n}{2} [g_{Fn}(-x) - g_0] \right] \quad , \quad \text{or}$$

(2.12)

$$\hat{f}(x) = f_0 + \sum_{n=-\infty}^{\infty} C_n [\bar{g}_{Fn}(x) - g_0]$$

In conclusion, we can formulate the following theorem:

**Theorem 4.** *Any function f(x) of $L^2[-L, L]$, can be developed into non-sinusoidal Fourier series, into a base composed of $f_0$ (its mean value over this interval) and the bases generated by the functions $[g(x)-g_0]$ and $[g(-x)-g_0]$. Here, g(x) is anything function of $L^2[-L, L]$ which has both two components (even and odd) non-zero, g(-x) is also of $L^2[-L, L]$, and $g_0$ is the average value of g(x) over $[-L, L]$.*

Consequently, the expansion of the function $f(x)$ analyzed earlier, on the interval $[-L, L]$, can be performed in a base generated by the functions $e^x$ et $e^{-x}$:

$$\hat{f}(x) = f_0 + \hat{f}_e(x) + \hat{f}_o(x) = f_0 + \sum_{n=1}^{\infty} A_n [\cosh_{Fn}(x) - \sinh 1] + \sum_{n=1}^{\infty} B_n \sinh_{Fn}(x) = f_0 + \sum_{n=1}^{\infty} \frac{A_n + B_n}{2} e_{Fn}^x +$$

$$+ \sum_{n=1}^{\infty} \left( \frac{A_n - B_n}{2} e_{Fn}^{-x} - A_n \sinh 1 \right) = f_0 + \sum_{n=1}^{\infty} \left[ \frac{A_n + B_n}{2} (e_{Fn}^x - \sinh 1) + \frac{A_n - B_n}{2} (e_{Fn}^{-x} - \sinh 1) \right]$$

where $e_{Fn}^x$ and $e_{Fn}^{-x}$ are the quasi−harmonics of the *n* order of the F-functions $(e^x)_n$, respectively $(e^{-x})_n$ (the extensions on the real axis of the F-functions $e^{nx}$, respectively $e^{-nx}$, defined on the intervals $[-1/n, 1/n]$).

The palette of functions that can serve as a basis for non-sinusoidal Fourier expansion is extremely wide:

- if g(x) is a polynomial in $[-L, L]$, its even component $g_e(x)$ contains the even powers of x, while its odd component $g_o(x)$ contains the odd powers
- if g(x) is an exponential function, $g_{Fe}(x)$ can be an even function $G_{Fe}(\cosh x)$ and $g_{Fo}(x)$ can be an odd function $G_{Fo}(\sinh x)$
- if g(x) is logarithmic: $\ln(A+x)$ (where $A>0$), f(x) can be developed only on a sub-interval $[a, b]$, included in the interval $(-A, A)$, with the bases rational functions of the form:

$$G_{Fe}(x) = \frac{1}{2} \ln(A^2 - x^2) - g_0 \quad \text{and} \quad G_{Fo}(x) = \frac{1}{2} \ln \frac{A+x}{A-x}$$

- if g(x) is a rational function of the form $1/(A+x)$, $A>0$, the function f(x) can be developed only on a sub-interval $[a, b]$, included in the interval $(-A, A)$, with the bases:

$$G_{Fe}(x) = \frac{1}{2} \left[ \frac{1}{A+x} + \frac{1}{A-x} \right] - g_0 = \frac{A}{A^2 - x^2} - g_0 \quad \text{and} \quad G_{Fo}(x) = \frac{1}{2} \left[ \frac{1}{A+x} - \frac{1}{A-x} \right] = -\frac{x}{A^2 - x^2}$$

- if g(x) is an irrational function of the form $\sqrt{A+x}$, $A>0$, the function f(x) can be developed only on a sub-interval $[a, b]$, included in the interval $(-A, A)$, with the bases:

$$G_{Fe}(x) = \frac{1}{2} \left[ \sqrt{A+x} + \sqrt{A-x} \right] - g_0 \quad \text{and} \quad G_{Fo}(x) = \frac{1}{2} \left[ \sqrt{A+x} - \sqrt{A-x} \right]$$

### 3. Quasi-sinusoidal periodic Fourier series

Another way to combine two functions g(x) and h(x), in order to form a basis for the $L^2[-L, L]$-space, is to choose the function h(x) as a translation of the function g(x): $h(x)=g(x+\alpha T)$, where $\alpha \in (0, 1)$. Among the bases obtained by this modality, there are a few that are complete and have a special property: the function g(x) has a single component (even or odd), and for $\alpha=1/4$, h(x) has the opposite parity, property that they have the functions $sin(\omega_0 x)$, respectively $cos(\omega_0 x)$, also. Because this property is extremely useful for solving some practical problems, in this section we will pay some attention to them.

If $g(x)$ is a function of $L^2$-space defined on the interval $[0, L/2]$, we can construct the functions $g_\&(x)$, composed of 4 segments, each segment explicitly defined with the help of the function $g(x)$, on a quarter of the interval $[-L, L]$. We impose that the functions $g_\&(x)$ obtained have the mean value zero over the interval $[-L, L]$, and that they have internal symmetries similar to those of the *sine*, respectively *cosine* functions: the two branches (this for $x<0$ and this for $x>0$) of the odd $g_\&(x)$ functions are symmetrical with respect to their **mid-axis**, and the two branches of the $g_\&(x)$ even functions are symmetric with respect to their **mid-point**. In addition, by translation, to the left or to the right, with $L/2$, an opposite parity function is obtained. We will call the functions $g(x)$, **the kernel** of the expansion, and the $g_\&(x)$ derived functions, the **quasi-sinusoids.** We will use the notations:

$g_s(x) = \mathbf{S}[g(x)]_L = G_s[-L < -g(x+L) > -L/2 < -g(-x) > 0 < g(x) > L/2 < g(L-x) > L]$, respectively
$g_c(x) = \mathbf{C}[g(x)]_L = G_c[-L < -g(x+L) > -L/2 < g(-x) > 0 < g(x) > L/2 < -g(L-x) > L]$

For example, if $g(x) = x$, for $x \in [0, \pi/2]$
$g_s(x) = \mathbf{S}[x]_\pi = G_s[-\pi < -x - \pi) > -\pi/2 < x > \pi/2 < \pi - x > \pi]$,
$g_c(x) = \mathbf{C}[x]_\pi = G_c[-\pi < -x - \pi) > -\pi/2 < -x > 0 < x > \pi/2 < -\pi + x > \pi]$

The functions $g_s(x)$ and $g_c(x)$ satisfy all the conditions required by Theorem 3, therefore we can conclude:

**Theorem 5**: *Let $g(x)$ be any function of $L^2[0, L/2]$-space. Any function $f(x)$ of $L^2[0, L/2]$-space can be developed in a quasi-sinusoidal Fourier series, into a base composed of $f_0$ (its mean value over this interval) and the bases generated by the functions $\mathbf{S}[g(x)]_L$ and $\mathbf{C}[g(x+L/2)]_L$, or $\mathbf{C}[g(x)]_L$ and $\mathbf{S}[g(x+L/2)]_L$.*

In the general case, the quasi-sinusoids and/or their first order derivatives have discontinuities which can be eliminated by adding quasi-sinusoids formed by rectangular pulses (one for each jump) and/or quasi-sinusoids formed by ramp-functions correctly chosen (for the odd quasi-sinusoids with a discontinuity at the origin, the ramp is $-2[g(L/2)-g(0)]/L$, while for even quasi-sinusoids, the ramp has the value of $-dg/dx$ at point $x=0$). For the functions obtained, we will use the name of **smooth quasi-sinusoids**, or **almost-sinusoids**, or **approximate sinusoids**, respectively **almost-cosinusoids**, or **approximate cosinusoids**. They are particularly useful in certain practical problems (resolution of certain differential equations with partial derivatives).

For example, to obtain the almost-sinusoids $\mathbf{S}[x^2-2x]_2$ of Fig.4a, we will choose two even functions $f_{2p}(x) = -x^2$ and $f_{2p}(x) = x^2$, for $-1 \leq x \leq 1$, and by vertical (by adding rectangular waves) and horizontal (by changes of variables) translations, we superimpose at the origin ($x=0$), the last point of the negative wave, with the first point of the positive wave. The almost-cosine $\mathbf{C}[1-x^2]_2$ is obtained by changing the variable $x$ with $x - 1$. The functions obtained are defined over the interval $[-2, 2]$, therefore having half of the angular frequency of the initial functions:

$g_s(x)_{2L} = \mathbf{S}[x^2-2x]_{2L} = G_s[-2L < -x^2 - 2x > 0 < x^2 - 2x > 2L]$, respectively
$g_c(x)_{2L} = \mathbf{C}[1-x^2]_{2L} = G_c[-2L < -1 + (x+2)^2 > -L < 1 - x^2 > L < -1 + (x-2)^2 > L]$

or, going back to the initial interval of definition:
$g_s(x)_L = \mathbf{S}[x^2-2x]_L = G_s[-L < -4x^2 - 4x > 0 < 4x^2 - 4x > L]$, respectively
$g_c(x)_L = \mathbf{C}[1-x^2]_L = G_c[-L < -1 + 4(x+1)^2 > -L/2 < 1 - 4x^2 > L/2 < -1 + 4(x-1)^2 > L]$

The function $g_c(x)_L$ is shown in Fig. 4b.

Similar to the Fourier series expansions, the presence of discontinuities inside or at the ends of the definition interval of the function $f(x)$ produces, for non-sinusoidal series expansions, additional terms (provided with coefficients of the $a_n/n$ form) and an effect similar to the Gibbs phenomenon. Likewise, the discontinuities of the first derivative generate other additional terms (provided with coefficients of the $a_n/n^2$ form) and additional oscillation phenomena with a significant amplitude.

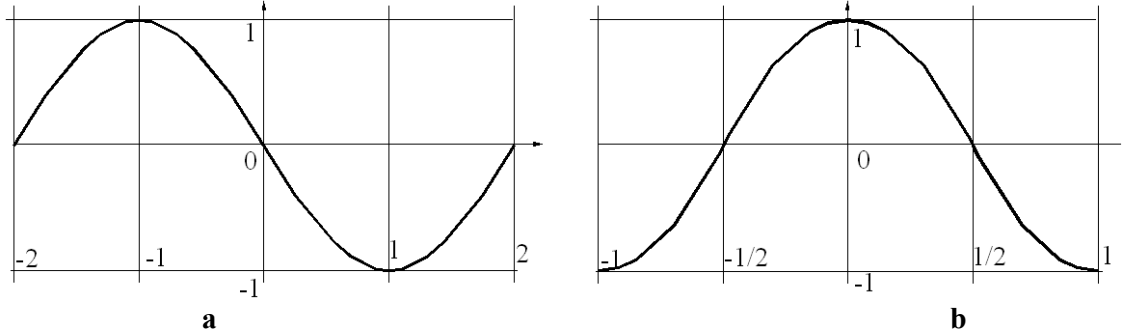

*Fig. 4: Almost-sinusoids of the 2nd degree*
**a**: *the odd function $S[x^2-2x]_2$;* **b**: *the even function $C[1-x^2]_1$*

The almost-sinusoids are part of the $C^1$ class of regularity (functions whose first derivative is continuous) and, because of their symmetry properties, similar to those of the sine and cosine functions, are best suited when series expansions of Fourier are required. Similarly, as in the practice of Fourier expansions, when the other data of the problem allow it, it is advantageous to construct for the developed function *f(x)* also, defined on an interval $[x_1, x_2]$, an extension quasi-sinusoidal smooth, defined on an interval $[x_{1e}, x_{2e}]$, which includes the interval of definition. By such an approach, the "residual terms" of expansion are removed.

For clarification, we will develop in non-sinusoidal series the almost-sinusoid:
$g_c(x)_L = C[1-x^2]_L = G_c[-L < -1 + 4(x+1)^2 > -L/2 < 1-4x^2 > L/2 < -1+4(x-1)^2 > L]$, what have the

expansion in Fourier series: $\bar{g}_c(x) = \sum_{n=1}^{\infty} \left\{ \frac{32(-1)^{n+1}}{(2n-1)^3 \pi^3} \cos(2n-1)\pi x \right\}$ for two bases:

1) the rectangular wave (2.c): $g_e = G_e[-1 > -1 < -1/2 > 1 < 1/2 > -1 < 1]$ and
2) the triangle wave (2.d): $f_{12}(x) = F_{12}[-1 > -x - 1/2 < 0 > x - 1/2 < 1]$ of the section 2.1.
Their expansions in Fourier sinusoidal series are:

$\bar{g}_e(x) = \sum_{n=1}^{\infty} c_n \cos(2n-1)\pi x = -\frac{2}{\pi} \sum_{n=1}^{\infty} \frac{(-1)^n \cos(2n-1)\pi x}{2n-1}$ and $\bar{f}_{12}(x) = -4 \sum_{n=1}^{\infty} \frac{\cos(2n-1)\pi x}{(2n-1)^2 \pi^2}$

In order to compare these expansions with those of the section 2.1, we will use the coefficients of the Fourier expansion of the function $-1/2 \cdot g_c(x)$, defined on the same interval. Its expansion in Fourier series generates the coefficients:
$a_1 = -16/\pi^3$, $a_2 = 0$, $a_3 = 16/27\pi^3$, $a_4 = 0$, $a_5 = -16/125\pi^3$, $a_6 = 0$, $a_7 = 16/343\pi^3$, $a_8 = 0$,
$a_9 = -16/729\pi^3$, $a_{10} = 0$, $a_{11} = 16/1331\pi^3$, $a_{12} = 0$, ...
For the rectangular wave, the relationships (2.2) lead us to:
$A_1 = -8/\pi^2$, $A_2 = 0$, $A_3 = -64/27\pi^2$, $A_4 = 0$, $A_5 = 192/125\pi^2$, $A_6 = 0$, $A_7 = -384/343\pi^2$, $A_8 = 0$,
$A_9 = 64/729\pi^2$, $A_{10} = 0$, $A_{11} = -960/1331\pi^2$, $A_{12} = 0$, ...
For the triangle wave, with the same relations, we calculate:
$A_1 = 4/\pi$, $A_2 = 0$, $A_3 = -16/27\pi$, $A_4 = 0$, $A_5 = -16/125\pi$, $A_6 = 0$, $A_7 = -32/343\pi$, $A_8 = 0$, $A_9 = 16/729\pi$,
$A_{10} = 0$, $A_{11} = -48/1331\pi$, $A_{12} = 0$,...

With the values thus obtained, we can construct the secondary quasi-harmonics and the first partial sums of the corresponding quasi-sinusoidal expansions: Figures 5 and 6.

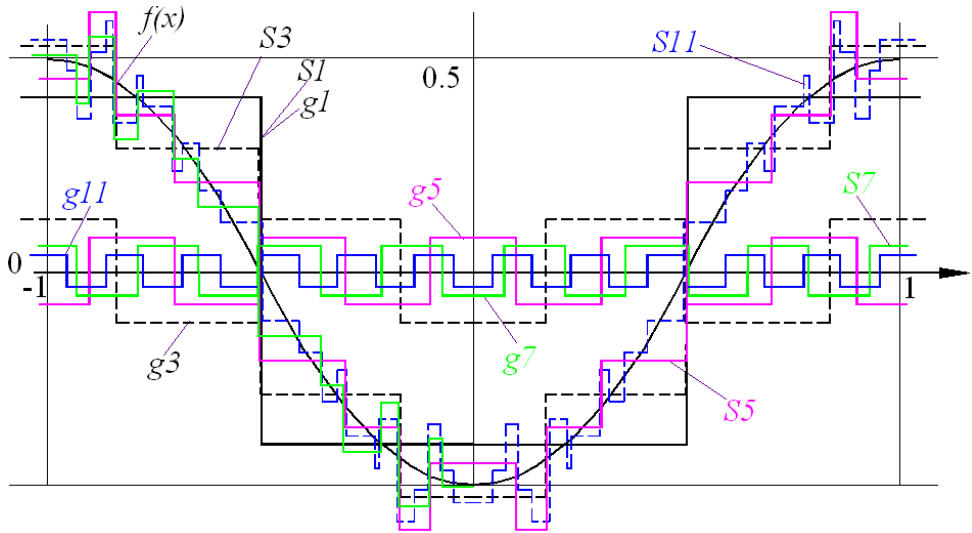

*Fig.5. The approximation of the quasi-sinusoid $-g_c/2$ by a sum of rectangular waves*

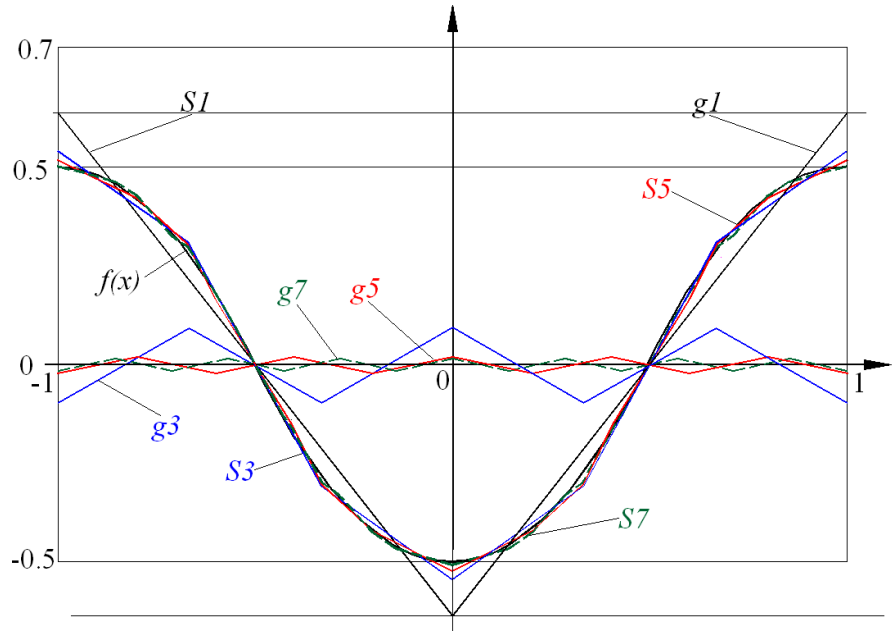

*Fig.6 The approximation of the quasi-sinusoid $-g_c/2$ by a sum of triangle waves*

The aspect of the expansion of quasi-sinusoids in almost-sinusoidal series is reflected in the inverse expansion of the functions analyzed above. Figures 7 and 8 show the first quasi-harmonics and the partial sums of order 12 for the inverse expansions $\hat{g}_e(x) = \sum_{n=1}^{\infty} A'_{dn}(g_c)_n$ and $\hat{f}_{12}(x) = \sum_{n=1}^{\infty} A'_m(g_c)_n$. According to (2.3):

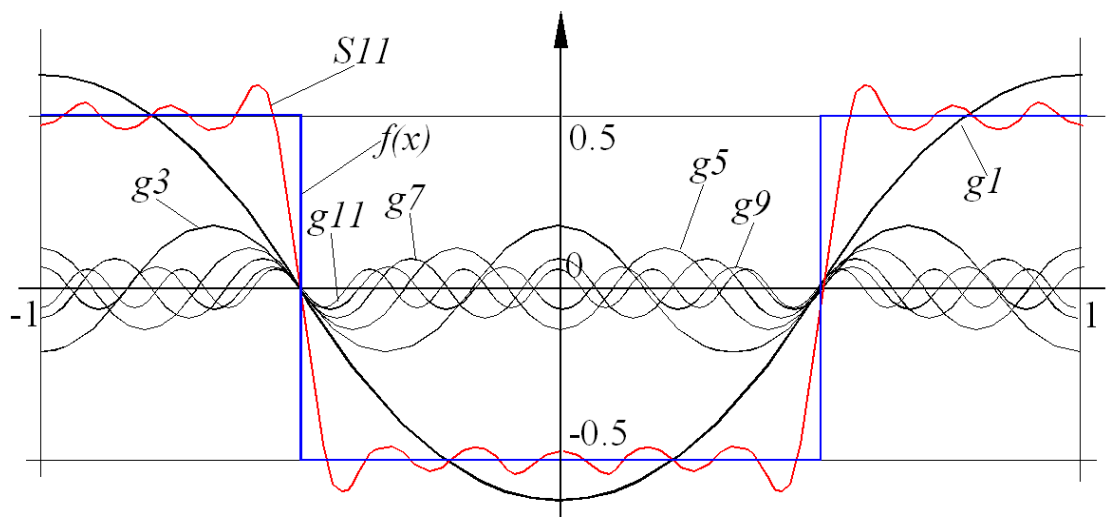

*Fig.7. The approximation of the rectangular waves by almost-sinusoids*

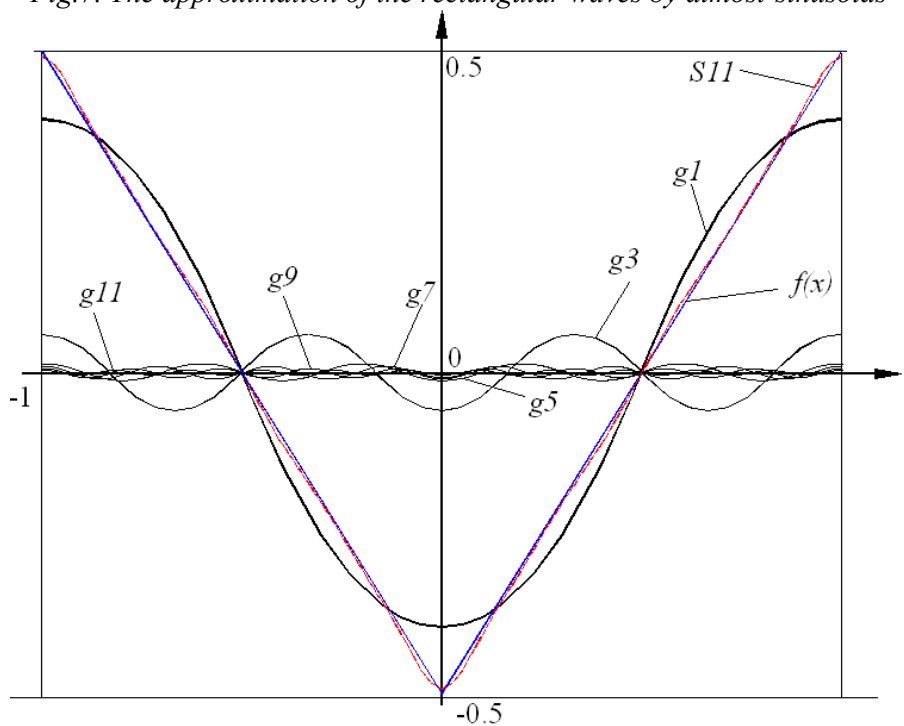

*Fig.8. The approximation of the triangle waves by almost-sinusoids*

$A'_1=-\pi^2/8$, $A'_2=0$, $A'_3=\pi^2/27$, $A'_4=0$, $A'_5=-3\pi^2/125$, $A'_6=0$, $A'_7=6\pi^2/343$, $A'_8=0$, $A'_9=-\pi^2/81$, $A'_{10}=0$, $A'_{11}=15\pi^2/1331$, $A'_{12}=0$, ... , respectivement:
$A'_1=-\pi/4$, $A'_2=0$, $A'_3=-\pi/27$, $A'_4=0$, $A'_5=-\pi/125$, $A'_6=0$, $A'_7=-2\pi/343$, $A'_8=0$, $A'_9=-\pi/243$, $A'_{10}=0$, $A'_{11}=-3\pi/1331$, $A'_{12}=0$, ...

    In both cases, the coefficients of the expansions have values close to those obtained by the expansions in sinusoidal series (Fourier). The differences get even smaller if the core of expansion is replaced by $g(x)=1-x^{1,75}$.

    In the following example, we will develop the rectangular wave (Fig.9.a) :
$f_o=F_o[-1>-1<0>1<1]$ into two non-sinusoidal series: one generated by the function $g_{o1}(x)=G_{o1}[-1>sinhx<1]$ (Fig.9.b), and the other by $g_{o2}(x)$ :

$$g_{o2}(x) = G_{o2}[-1 < \frac{[\cosh 0.5 - \cosh(x-0.5)]}{\cosh 0.5 - 1} > 0 < \frac{[-\cosh 0.5 + \cosh(x+0.5)]}{\cosh 0.5 - 1} > 1]$$

$g_{o2}(x)$ is generated by the odd almost-sinusoid of exponential type (Fig. 10.a) which has a kernel of the type $g(x)=K+\cosh(x+T/4)$, where $K=constant$:

In the interval $[-1, 1]$: $f_{o0}=g_{o10}=g_{o20}=0$, $\bar{f}_o(x) = \sum_{n=1}^{\infty} b_n \cos n\pi x = \frac{4}{\pi} \sum_{n=1}^{\infty} \frac{\sin(2n-1)\pi x}{2n-1}$ , so:

$b_1=4/\pi$, $b_2=0$, $b_3=4/3\pi$, $b_4=0$, $b_5=4/5\pi$, $b_6=0$, $b_7=4/7\pi$, $b_8=0$, ... ,

$$\bar{g}_{o1}(x)_1 = \sum_{n=1}^{\infty} [c_n \sin(n\pi x)] = \pi \cdot \sinh 1 \sum_{n=1}^{\infty} \frac{2n \cdot (-1)^{n+1}}{1+n^2\pi^2} \sin(n\pi x)$$

$$\bar{g}_{o2}(x) = \sum_{n=1}^{\infty} d_n \sin n\pi x = \sum_{n=1}^{\infty} \frac{4 \cdot \cosh 0.5}{\pi(\cosh 0.5 - 1) \cdot (2n-1) \cdot [1+(2n-1)^2 \pi^2]} \sin n\pi x$$ , so:

$d_1=8.8372 \cdot (4/\pi)$, $d_2=0$, $d_3=8.8372 \cdot (4/3\pi) \cdot 1/(1+4\pi^2)$, $d_4=0$, $d_5=8.8372 \cdot (4/5\pi) \cdot 1/(1+9\pi^2)$, $d_6=0$, $d_7=8.8372 \cdot (4/7\pi) \cdot 1/(1+16\pi^2)$, $d_8=0$, ... ,

For $K=8.8372$:
$B_1=K$, $B_2=0$, $B_3=K/(1+4\pi^2)$, $B_4=0$, $B_5=K/(1+9\pi^2)$, $B_6=0$, $B_7=K/(1+16\pi^2)$, $B_8=0$, $B_9=K/(1+25\pi^2)$, $B_{10}=0$, $B_{11}=K/(1+36\pi^2)$, $B_{12}=0$, ...

A graphic representation of the partial sum S12 is given in the figure 10.b.

For sub unitary values of $L$, the exponential quasi-sinusoid $g_{o2}(x)$ is satisfactorily approximated by the quadratic quasi-sinusoid $g_c(x)_1 = C[1-x^2]_1$, described in Figure 4.b of section 3. Between the coefficients of the expansion of the function $f_o(x)$ in the two bases generated by the quasi-sinusoids $g_c(x)_1$ and $g_{o2}(x)$, the differences are negligible.

### 4. Orthogonal bases composed of non-sinusoidal periodic functions

Neither the even quasi-harmonics $g_n(x)-g_0$ nor the odd $h_n(x)$, analyzed in the previous sections, are orthogonal to each other, which does not allow the calculation of the coefficients of these expansions from formulas similar to Euler's formulas. But, any even quasi-harmonic is orthogonal to all odd quasi-harmonics. This allows us, by the Gram-Schmidt orthogonalization process [4], to construct an orthogonal basis (which can be normalized by the same method) for each of the systems generated by the Fourier-functions $g_{Fn}(x)-g_0$ and $h_{Fn}(x)$. By combining them and adding the function $f_0$, a complete biorthogonal basis is obtained. The Gram-Schmidt orthogonalization process does not claim, for the g-harmonics of the non-orthogonal basis, the need to be continuous (to be Fourier-functions), but this functionality is imposed by our intention to create a complete base for the space $L^2[-L, L]$. Because of this, we consider the functions $g_{Fn}(x)$ and $h_{Fn}(x)$ to be Fourier-functions by definition (1.d).

For example, from some two functions $g(x)$ even and $h(x)$ odd, which have the mean zero value, defined over a certain interval $[a, b]$, we obtain a biorthogonal basis formed by the functions $1, \Phi_n(x)$ et $\Psi_n(x)$, $n=1, 2, 3,...$, where:

$$\Phi_n(x) = g_{Fn}(x) - \sum_{j=1}^{n-1} \frac{\int_a^b g_n(x)\Phi_j(x)dx}{\|\Phi_j\|^2} \Phi_j(x) = g_{Fn}(x) - \sum_{i=1}^{n-1} C_{in} g_{Fi}(x) \text{ and}$$

$$\Psi_n(x) = h_{Fn}(x) - \sum_{j=1}^{n-1} \frac{\int_a^b h_n(x)\Psi_j(x)dx}{\|\Psi_j\|^2} \Psi_j(x) = h_{Fn}(x) - \sum_{i=1}^{n-1} D_{in} h_{Fi}(x) \quad (4.1)$$

with $\|\Phi_j\|^2 = \int_a^b \Phi_j^2(x)dx$ and $\|\Psi_j\|^2 = \int_a^b \Psi_j^2(x)dx$

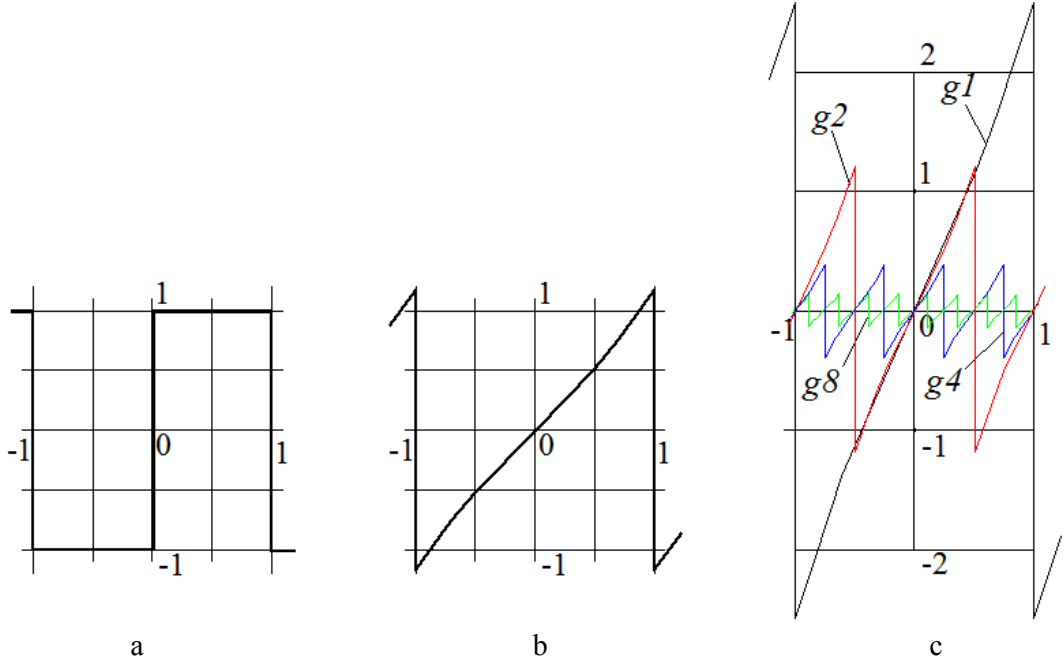

Fig. 9. *The expansion of the function $f_o(x)$ in the base $g_{o1}(x)$*
*a*: $f_{o1}(x)$  *b*: $g_{o1}(x)=sinh_1 x$  *c*: *the quasi-harmonics $g_1(x)$, $g_2(x)$, $g_4(x)$, $g_8(x)$*

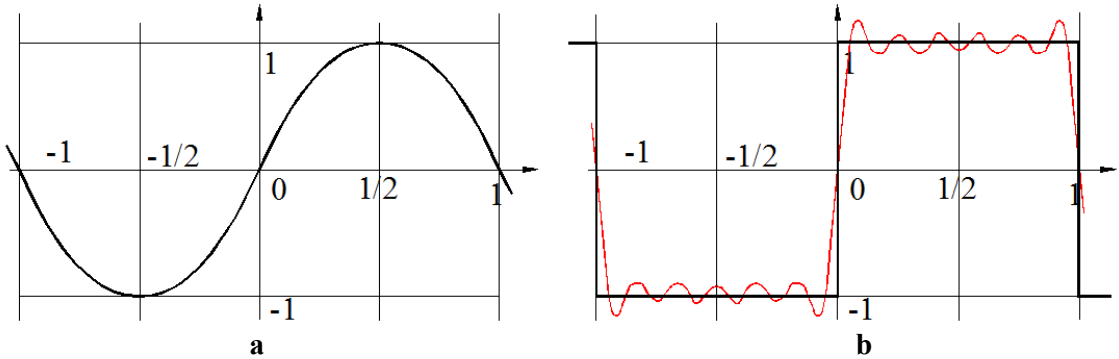

Fig. 10. *The expansion of the function $f_o(x)$ in the base $g_{o2}(x)$*
*a*: $g_{o2}(x)$  *b*: *the partial sum $S_{12}(x)$*

These considerations allow us to formulate the
**Theorem 6**: *Let be two some F-functions $g_F(x)$ even and $h_F(x)$ odd of $L^2[a, b]$-space. Any function $f(x)$ of $L^2[a, b]$-space can be developed into a complete series, based on the bi orthogonal system 1, $\Phi_n(x)$ and $\Psi_n(x)$, where $\Phi_n(x)$ and $\Psi_n(x)$ are generated by the functions $g_{Fn}(x)-g_0$ and $h_{Fn}(x)$ by an orthogonalization process:*

$$\tilde{f}(x) = A_0 + \sum_{n=1}^{\infty}\left[A_n^0 \Phi_n(x) + B_n^0 \Psi_n(x)\right] \tag{4.2}$$

Thanks to the orthogonality of the system, to calculate the coefficients of this expansion, are valid the formulas of Euler:

$$A_0 = \frac{1}{b-a}\int_a^b f(x)dx$$

$$A_n^0 = \frac{1}{\|\Phi_n\|^2} \int_a^b f(x)\Phi_n(x)dx \quad n=1, 2, 3,... \tag{4.3}$$

$$B_n^0 = \frac{1}{\|\Psi_n\|^2} \int_a^b f(x)\Psi_n(x)dx \quad n=1, 2, 3,...$$

Obtaining these expressions is based on relationships:
$$\langle f(x)|\Phi_n(x)\rangle = \langle [f_e(x)+f_o(x)]|\Phi_n(x)\rangle = \langle f_e(x)|\Phi_n(x)\rangle$$
$$\langle f(x)|\Psi_n(x)\rangle = \langle [f_e(x)+f_o(x)]|\Psi_n(x)\rangle = \langle f_o(x)|\Psi_n(x)\rangle$$
which are true for $\Phi_n(x)=cos(n\omega_0 x)$, respectively $\Psi_n(x)=sin(n\omega_0 x)$, too.

We can see that $\Phi_n(x)$, the components of order $n$ ($n=1, 2, 3,...$) of the orthogonal system, generated by the even functions $g_{Fn}(x)-g_0$, as well as $\Psi_n(x)$, generated by the odd functions $h_{Fn}(x)$, are linear combinations between the quasi-harmonics of order $n$ and the quasi-harmonics of lower order of the respective non-orthogonal expansions. Consequently, we can establish a correspondence between the coefficients $A_n$ and $B_n$ of the expansion in the non-orthogonal base generated by the functions $g_F(x)$ and $h_F(x)$ and those of the expansion in the orthogonal base $\Phi_n(x)$ and $\Psi_n(x)$:

$$\tilde{f}(x) = A_0 + \sum_{n=1}^{\infty} A_n^0 \left[g_{Fn}(x) - \sum_{i=1}^{n-1} C_{in}g_{Fi}(x)\right] + \sum_{n=1}^{\infty} B_n^0 \left[h_{Fn}(x) - \sum_{i=1}^{n-1} D_{in}h_{Fi}(x)\right] =$$

$$= A_0 + \sum_{n=1}^{\infty}\left(A_n^0 - \sum_{i=n+1}^{\infty} A_n^0 C_{ni}\right)g_{Fn}(x) + \sum_{n=1}^{\infty}\left(B_n^0 - \sum_{i=n+1}^{\infty} B_n^0 D_{ni}\right)h_{Fn}(x) = A_0 + \sum_{n=1}^{\infty} A_n g_{Fn}(x) + \sum_{n=1}^{\infty} B_n h_{Fn}(x)$$

We note that for the calculation of the coefficients $A_n = A_n^0 - \sum_{i=n+1}^{\infty} A_n^0 C_{ni}$ and $B_n = B_n^0 - \sum_{i=n+1}^{\infty} B_n^0 D_{ni}$, it is necessary to calculate certain integrals of the type:

$$\int_a^b f_e(x)g_n(x)dx, \int_a^b f_o(x)h_n(x)dx, \int_a^b g_i(x)g_j(x)dx, \int_a^b h_i(x)h_j(x)dx, \text{ for } i, j=1, 2, 3, ...$$

By this method of calculating the coefficients, it is no longer necessary to know the coefficients of the expansions in sinusoidal series neither for the function $f(x)$, nor for the functions $g(x)$ and $h(x)$.

We are going to exemplify by building an orthogonal base, starting from a base generated by the system of periodic square unitary pair functions:
$f_e=F_e[-1>1<-1/2>-1<1/2>1<1]$

In the case chosen here, the calculation will be simplified thanks to the symmetry properties of the chosen quasi-sinusoid. Thanks to relations (4.1), it follows:

$$\Phi_1(x) = g_1(x), \qquad \|\Phi_1\|^2 = \int_{-1}^{1} g_1^2(x)dx = 2$$

$$\Phi_2(x) = g_2(x) - \frac{\int_{-1}^{1} g_2(x)g_1(x)dx}{2}g_1(x) = g_2(x), \qquad \|\Phi_2\|^2 = \int_{-1}^{1} g_2^2(x)dx = 2$$

$$\Phi_3(x) = g_3(x) - \frac{\int_{-1}^{1} g_3(x)g_1(x)dx}{2}g_1(x) - \frac{\int_{-1}^{1} g_3(x)g_2(x)dx}{2}g_2(x) = g_{F3}(x) - \frac{1}{3}g_{F1}(x) = g_3(x) - C_{13}g_1(x)$$

$$\|\Phi_3\|^2 = \int_{-1}^{1}\left[g_3(x) - \frac{1}{3}g_1(x)\right]^2 dx = \frac{16}{9}$$

$$\Phi_4(x) = g_4(x) - \frac{\int_{-1}^{1} g_4(x)g_1(x)dx}{2} g_1(x) - \frac{\int_{-1}^{1} g_4(x)g_2(x)dx}{2} g_2(x) -$$

$$- \frac{9\int_{-1}^{1} g_4(x)\left[g_3(x) - \frac{1}{3}g_1(x)\right]dx}{16} \left[g_3(x) - \frac{1}{3}g_1(x)\right] = g_4(x) - C_{14}g_1(x) - C_{24}g_2(x) - C_{34}g_3(x) = g_4(x)$$

$$\|\Phi_4\|^2 = \int_{-1}^{1} g_4^2(x)dx = 2$$

$$\Phi_5(x) = g_5(x) - \frac{\int_{-1}^{1} g_5(x)g_1(x)dx}{2} g_1(x) - \frac{9\int_{-1}^{1} g_5(x)\left[g_3(x) - \frac{1}{3}g_1(x)\right]dx}{16} \left[g_3(x) - \frac{1}{3}g_1(x)\right] =$$

$$= g_5(x) - \frac{9}{40} g_1(x) - \frac{3}{40} g_3(x) = g_5(x) - C_{15}g_1(x) - C_{35}g_3(x)$$

$$\|\Phi_5\|^2 = \int_{-1}^{1} \left[g_5(x) - \frac{17}{40} g_1(x) - \frac{3}{40} g_3(x)\right]^2 dx = \frac{1607}{800} \quad \text{and so on.}$$

Through similar relationships: $\Psi_n(x) = h_{Fn}(x) - \sum_{i=1}^{n-1} D_{in} h_{Fi}(x)$ is obtained the orthogonal system $\Psi_n(x)$, from the function $f_o(x)$, unitary odd periodic rectangular waves: $f_o(x) = F_o[-1 > -1 < 0 > 1 < 1]$

## 5. Properties of non-sinusoidal Fourier series

Work on expansions in trigonometric (sinusoidal) Fourier series has shown that a real function $f(x)$ from the space $L^2[-L, L]$, can be expressed as a sum of its projections on the components of an orthogonal base of a space of functions, if it fulfills several conditions. The series resulting from these expansions have several properties: convergence, summability, differentiability, integrability.

The expansions into non-sinusoidal Fourier series in $g(x)$-basis of the function $f(x)$ analyzed in the previous sections have resulted from the sinusoidal series of this function, by a redistribution of its coefficients. This redistribution reconstruct the coefficients of expansion into sinusoidal series of the components $g_n(x)$ of a complete base of non-orthogonal functions. Therefore, new expansions in non-sinusoidal Fourier series transfer from the expansions in sinusoidal Fourier series a palette of conditionalities and properties. Undoubtedly, this subject merits further study, but for the moment, we boil down to a few obvious conclusions:

- if the function $f(x)$ is integrable, the sequence of the coefficients of its expansion in non-sinusoidal series converges towards $0$ (the Riemann-Lebesgue theorem)

- if $\tilde{f}(x) = A_0 + \sum_{n=1}^{\infty} \left[A_n^0 \Phi_n(x) + B_n^0 \Psi_n(x)\right]$ is the expansion of the function $f(x)$ into a base of periodic non-sinusoidal orthogonal functions, then

$$A_0^2 + \frac{1}{2} \sum_{n=1}^{\infty} \left[(A_n^0)^2 + (B_n^0)^2\right] = \frac{1}{2L} \int_{-L}^{L} |f(x)|^2 dx \quad \text{(the Parceval theorem)}$$

- if the $2L$-periodic functions $f(x)$, $g(x)$, and $h(x)$, in the interval $[-L, L]$ are continuous by pieces and differentiable on the left and on the right in all the point of the interval, the series $\hat{f}(x) = f_0 + \sum_{n=1}^{\infty} A_n g_n(x) + \sum_{n=1}^{\infty} B_n h_n(x)$ converges to

$$\hat{f}(x_i) = \frac{1}{2}\left[\lim_{x \to x_i^+} f(x) + \lim_{x \to x_i^-} f(x)\right] \text{ in all points } x_i.$$

- the Fourier series resulting from the expansion in any non-sinusoidal base, of an *f(x)* function, *2L*-periodic, continuous and differentiable by pieces, converges uniformly on **R** towards this function.
- the non-sinusoidal Fourier series of a *2L*-periodic, square-integrable function that can be integrated over a period, converges in norm of $L^2$ to the considered function (the Riesz–Fischer theorem)
- the non-sinusoidal Fourier series of a square-integrable functions converges almost everywhere to this function (Carleson's theorem)
- two *2L*-periodic functions, having the same coefficients of their expansion in the same non-sinusoidal Fourier base, are equal almost everywhere. In particular, in the case of continuity by pieces, they coincide in all the points of the interval [−L, L], except a finite number
- let be *f(x)* a function *2L*-periodic, continues in the interval [−L, L]. Its Fourier expansion $f(x) = f_0 + \sum_{n=1}^{\infty} A_n g_n(x) + \sum_{n=1}^{\infty} B_n h_n(x)$, sinusoidal or not, convergent or not, can be integrated term by term, between all integration limits:

$$\int f(x)dx = d_0 + f_0 x + \sum_{n=1}^{\infty} A_n \int g_n(x)dx + \sum_{n=1}^{\infty} B_n \int h_n(x)dx, \text{ where } d_0 \text{ is an arbitrary constant.}$$

- let be *f(x)* a *2L*-periodic function, continuous in the interval [−L, L], with *f(−L)=f(L)* and with the derivative *f'(x)* smooth by pieces in this interval. The Fourier expansion, sinusoidal or not, of the function *f'(x)*, can be obtained by deriving term by term the Fourier expansion of the function *f(x)*. The series obtained converges punctually towards *f'(x)* in all the points of continuity and towards *[f'(x)+ f'(−x)]/2* in those of discontinuity. If $f(x) = f_0 + \sum_{n=1}^{\infty} A_n g_n(x) + \sum_{n=1}^{\infty} B_n h_n(x) \to f'(x) = \sum_{n=1}^{\infty} A_n g_n'(x) + \sum_{n=1}^{\infty} B_n h_n'(x)$

The condition *f(−L)=f(L)* imposed in this statement is quite restrictive and reduces the usefulness of the theorem. We can get around this condition if we take into account that the even type component $f_e$ of the function *f(x)* always satisfies the differentiability condition, and that the odd component $f_o$ can be written as a sum of the differentiable function $f_{os}$ and of the ramp-function: $f_r = x \cdot f_o(L)/L$. So:

$$\frac{d}{dx} f_o(x) = \frac{d}{dx}\left[f_{os}(x) + \frac{f_o(L)}{L}x\right] = \frac{d}{dx} f_{os}(x) + \frac{f_o(L)}{L}$$

For example, in the case of sinusoidal Fourier series expansion:

$$f'(x) = \frac{f_o(L)}{L} + \sum_{n=1}^{\infty}\left[b_n \omega_n \cos(\omega_n x) - \left(a_n \omega_n + 2(-1)^n \frac{f_o(L)}{L}\right)\sin(\omega_n x)\right]$$

which requires the knowledge of the boundary conditions $f_o(-L)$ and $f_o(L)$. This relation makes it possible to solve certain differential equations by determining the coefficients of the sinusoidal Fourier series expansion of the unknown function (similar to the expansion in series of Taylor).

### 6. Conclusions

We noted in the previous sections that any function *f(x)*:[−L, L], *2L*-periodic, which belongs to the $L^2$-subspace, can be developed, in a similar way to that indicated by Fourier, there is over 200 years ago, in a multitude of variants, in bases formed by $f_0$, the average value of the function *f(x)* in the interval [−L, L] and two sets of quasi-harmonics: a set of even functions $g_n(x)$ and a set of odd functions $h_n(x)$, ($n=1, 2, 3, ..., \infty$), periodic functions, with the period *2L/n*. In the most general case, the fundamental quasi-harmonics (for *n=1*) are all functions which satisfy the Dirichlet conditions. They can therefore be non-trigonometric functions and the base can be non-orthogonal.

The expansion in Fourier sinusoidal series is only one particular case of this expansion, namely the case where the fundamental quasi-harmonics are sinusoidal:
$g_{o1}(x)=sin(\omega_0 x)$ and $g_{e1}(x)=cos(\omega_0 x)$.

These results generate a wide range of theoretical results. First, a new, extremely broad perspective opens up in the analysis of function spaces, in their spectral analysis, in the development of new types of integral transforms, in the construction of wavelet function systems, etc.

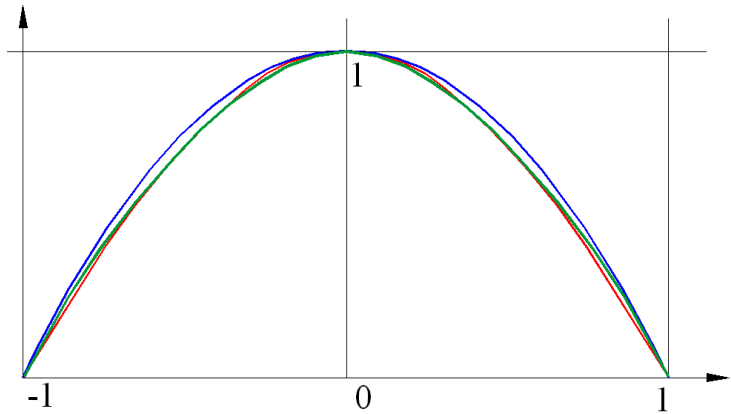

*Fig.11. Comparison between the curves $cos(x\pi/2)$ (red), $1-|x|^{1.75}$ (green) and $1-x^2$ (blue) on the interval [−1, 1]*

The comparison of the expansions of a function *f(x)* in different complete bases (the number of accessible bases has now become very large), makes it possible to solve new problems of convergence of the numerical series and series of functions, to find new correlations between different types of functions, etc. Second, the practical consequences are also extremely important. This new type of expansion leads to the development of new methods of approximation of functions, in which the precision of the approximation can be increased by the possibility of choosing from a wider range of possibilities. For example, in Figure 11, two possibilities of approximating a cosine are illustrated.

A perspective of extremely promising use of these types of series expansion is offered by the field of numerical and analytical resolution of broad categories of ordinary differential equations and with partial derivatives, linear and nonlinear.

# SÉRIES DE FOURIER NON SINUSOÏDALES

Arpad Török[1], Stoian Petrescu[2], Michel Feidt[3]

[1] PhD student, The Polytechnic University of Bucharest, Department of Engineering Thermodynamics, 313, Splaiul Independentei, 060042 Bucharest, Romania, e-mail: arpi_torok@yahoo.com
[2] Prof. Dr. Eng., Polytechnic University of Bucharest, Department of Engineering Thermodynamics, România
[3] Prof. Dr. Eng., L.E.M.T.A., U.R.A. C.N.R.S. 7563, Université de Loraine Nancy 12, avenue de la Foret de Haye, 54516 Vandeuvre-lès-Nancy, France

**Résumé**: Selon l'analyse harmonique, toute fonction $f(x)$, périodique sur l'intervalle $[-L, L]$, qui satisfait les conditions de Dirichlet, peut être développée en une somme infinie (connue dans la littérature sous le nom de **série trigonométrique** et pour laquelle, pour des raisons qui deviendront évidentes au cours de ce travail, nous utiliserons le nom de **série sinusoïdale**), constituée des composantes pondérées, d'une base biortogonale complète, formée de la fonction unitaire $1$, des harmoniques fondamentales, $2L$-périodiques, $sin(\pi x/L)$ paire et $cos(\pi x/L)$ impaire et des harmoniques secondaires, $2L/n$-périodiques, $sin(n\pi x/L)$ et $cos(n\pi x/L)$, où $n \in \mathbf{Z}^+$. Les coefficients de ces développements (coefficients de Fourier) peuvent être calculés à l'aide des formules d'Euler. Cette affirmation peut être généralisée, si nous prouvons que la fonction $f(x)$ peut également être développée en des **séries périodiques non sinusoïdales**, formées de la somme infinie des composantes pondérées d'une base complète, non orthogonale: la fonction unitaire $1$, les quasi-harmoniques fondamentales $g(x)$, paire et $h(x)$, impaire (des fonctions $2L$-périodiques, avec la valeur moyenne nulle sur l'intervalle de définition) et des quasi-harmoniques secondaires $g_n(x)$ et $h_n(x)$ (fonctions $2L/n$-périodiques), où $n \in \mathbf{Z}^+$. Les quasi-harmoniques fondamentales $g(x)$ et $h(x)$ sont n'importe quoi fonctions qui admettent des développements en série sinusoïdale (elles satisfont les conditions de Dirichlet, ou elles appartiennent à l'espace $L^2$). Les coefficients de ces développements sont obtenus avec l'aide de certains relations algébriques entre les coefficients de Fourier des développements des fonctions $f(x)$, $g(x)$ et $h(x)$. En plus de leur importance théorique évidente, ces types de développements peuvent avoir une importance pratique dans l'approximation des fonctions et dans la résolution numériques et analytiques de certains classes d'équations différentielles ordinaires et avec des dérivées partielles.

**Mots clefs**: séries de Fourier sinusoïdales, séries de Fourier non sinusoïdales, bases indépendants, bases orthogonales, approximation des fonctions, équations différentielles

### 1. Introduction

Pour définir simultanément à la fois la fonction $g(x)$ et son domaine de définition, pour une quelconque fonction $g(x)$ de $\mathbf{R}$, définie sur un intervalle réel $I$ (ouvert, fermé ou semi-ouvert), délimité par $x_1$ et $x_{m+1}$, introduite par une seule expression $g(x)=g_1(x)$ sur toute $I$ ($m=1$), ou par $m$ différentes expressions: $g(x)=g_i(x)$, $i=1, 2, ..., m,$ valables chacune sur l'un des sous-intervalles disjoints adjacents bornés par $x_i$ et $x_{i+1}$, $1 \leq i \leq m$, nous utiliserons un système pour marquer simultanément les limites $x_i$ des (sous)intervalles de définition et aussi les expressions de la fonction valables sur ces (sous)intervalles:

$g(x) = G[x_1^+ > g_1(x) < x_2^-]$, ou $g(x) = G[x_1 \geq g_1(x) \leq x_2]$, ou $g(x) = G[x_1^+ * g_1(x) * x_2^-]$, ou

$g(x) = G[x_1^+ (g_1(x)) x_2^-]$, au lieu de $g(x) = g_1(x)$ si $x_1 \leq x \leq x_2$. Ici $I = [x_1, x_2]$

$g(x) = G[x_1^+ > g_1 < x_2 \cup x_2^+ > ... < x_i \cup x_i^+ > g_i < x_{i+1} \cup x_{i+1}^+ > ... < x_m \cup x_m^+ > g_m < x_{m+1}^-]$, ou

$g(x) = G[x_1 \geq g_1 < x_2 \cup x_2 \geq ... < x_i \cup x_i \geq g_i < x_{i+1} \cup x_{i+1} \geq ... < x_m \cup x_m \geq g_m \leq x_{m+1}]$, ou

$$g(x) = G\left[x_1^+(g_1)x_2^- \cup x_2^+(g_2)...(g_{i-1})x_i \cup x_i^+(g_i)x_{i+1} \cup x_{i+1}^+(g_{i+1})...(g_{m-1})x_m \cup x_m^+(g_m)x_{m+1}^-\right], \text{etc}$$

au lieu de $g(x) = \{g_1(x)$ si $x_1 \leq x < x_2, ..., g_i(x)$ si $x_i \leq x < x_{i+1}, ..., g_m(x)$ si $x_m \leq x \leq x_{m+1}\}$

Ici $I = [x_1, x_2) \cup [x_2, x_3) \cup ... \cup [x_i, x_{i+1}) \cup ... \cup [x_m, x_{m+1}] = [x_1, x_2]$.

Les indices supérieurs attachés aux bornes $x_i$ de ces intervalles signifient:

- $+$ : $g(x_i) = \lim_{x \to x_i^+} g(x)$
- $-$ : $g(x_i) = \lim_{x \to x_i^-} g(x)$
- sans indice : la fonction $g(x)$ est indéfinie en $x_i$, ou elle a une valeur fixé $g(x_i) = a \in \mathbf{R}$.

Par conséquent:

- pour $g(x):(x_1, x_2)$ nous allons utiliser l'une des notations:

$g(x) = G[x_1 > g < x_2]$, ou $g(x) = G[x_1*(g)*x_2]$, ou $g(x) = G[x_1(g)x_2]$

- pour $g(x):[x_1, x_2]$, avec $g(x) = g_1(x)$ si $x \in (x_1, x_2)$ ; $g(x_1) = a$, $g(x_2) = b$, nous allons utiliser :

$g(x) = G[a < x_1 > g_1 < x_2 > b]$, $g(x) = G[(a)*x_1*(g_1)*x_2*(b)]$ , ou $g(x) = G[(a)x_1(g_1)x_2(b)]$

- pour $g(x):(x_1, x_2) \cup (x_2, x_3)$, avec $g(x) = g_1(x)$ si $x \in (x_1, x_2)$ et $g(x) = g_2(x)$ si $x \in (x_2, x_3)$ :

$g(x) = G[x_1 > g_1 < x_2 > g_2 < x_3]$, etc.

- pour $g(x):[x_1, x_2) \cup (x_3, x_4]$, avec $g(x) = g_1(x)$ si $x \in (x_1, x_2)$ et $g(x) = g_2(x)$ si $x \in (x_3, x_4)$:

$g(x) = G[x_1^+ > g_1 < x_2 > \cup < x_3 > g_2 < x_4^-]$, ou $g(x) = G[x_1 \geq g_1 < x_2 > \cup < x_3 > g_2 \leq x_4]$

- pour $g_a(x):[x_1, x_{m+1}]$, avec $g_a(x) = g_i(x)$ si $x \in (x_i, x_{i+1})$, et $g_a(x_i) = a_i \neq \infty$, $i = 1, 2, ..., m+1$:

$$g_a(x) = G[a_1 < x_1 > g_1 < x_2 > a_2 < x_2 > g_2 < x_3 > a_3 ... a_m < x_m > g_m < x_{m+1} > a_{m+1}] \quad (a)$$

- pour $g_b(x): \bigcup_{i=1}^{m}(x_i, x_{i+1})$, avec $g_b(x) = g_i(x)$ si $x \in (x_i, x_{i+1})$, $i = 1, 2, ..., m+1$, et $g_b(x_i)$ indéfinies:

$$g_b(x) = G[x_1 > g_1 < x_2 > g_2 < x_3 > ... < x_m > g_m < x_{m+1}]. \quad (b)$$

Pour ces deux derniers exemples, si $g_i(x)$ sont les mêmes, pour $i = 1, 2, ..., m+1$, les deux fonctions (ainsi que d'autres fonctions $g_r(x)$ qui accomplies cette condition et, encore la condition $g_r(x_i) = r_i \neq \infty$), sont *égaux presque partout*. S'ils sont de carré intégrable, leurs développements de Fourier $\overline{g}_a(x)$, $\overline{g}_b(x)$ et $\overline{g}_r(x)$ ont la même expression $\overline{g}(x)$ quelle, en les points de discontinuité, converge vers:

$$\overline{g}(x_1) = \overline{g}(x_{m+1}) = \frac{1}{2}\left[\lim_{x \to x_1^+} f(x) + \lim_{x \to x_{m+1}^-} f(x)\right] \text{ et } \overline{g}(x_i) = \frac{1}{2}\left[\lim_{x \to x_i^+} f(x) + \lim_{x \to x_i^-} f(x)\right], \text{ pour } i \neq 1, m+1$$

Nous noterons $g_F(x)$ une fonction de type $g_r(x)$, pour laquelle $g_r(x_i) = \overline{g}(x_i)$, $i = 1, 2, ..., m+1$:

$$g_F(x) = G[\overline{g}(x_1) < x_1 > g_1 < x_2 > \overline{g}(x_2) < x_2 > g_2 < x_3 > \overline{g}(x_3) ... \overline{g}(x_m) < x_m > g_m < x_{m+1} > \overline{g}(x_1)] \quad (c)$$

On peut noter $g_a(x) \stackrel{F}{=} g_b(x) \stackrel{F}{=} g_r(x) \stackrel{F}{=} g_F(x) \stackrel{F}{=} \lim_{N \to \infty} \overline{g}(x)$, où $\stackrel{F}{=}$ veux dire que l'égalitée est *presque partout*. Le développement de Fourier de la fonction discontinue $g(x)$ est une fonction continue $\overline{g}(x)$ qui s'approche autant que possible de $g_F(x)$. On peut assigner à la fonction continue $\overline{g}_F(x)$ (qui est une approximation autant que nous voulons de proche de la fonction $g_F(x)$), l'appellation de **Fourier-fonction,** ou **F-fonction.** Dans des nombreuses situations, y compris la majorité de cet article, les valeurs $r_i$ de la fonction $g_r(x)$ en les points de discontinuité ne sont pas pertinentes. Par conséquent, lorsque nous analyserons des fonctions de ce type, sans perdre le caractère de généralité, mais dans un souci de simplification de l'exposition, nous allons considérer toujours (sauf les cas spécifiés expressément) qu'il s'agit d'une fonction du type $g_F(x)$, et nous allons utiliser la plus simple notation équivalente, celle de $g_b(x)$, la relation (b):

$g_F(x) = G[x_1 > g_1 < x_2 > g_2 < x_3 > ... < x_m > g_m < x_{m+1}]$

Autres exemples:

pour la fonction Heaviside: $H(x) = G[-\infty > 0 < 0 > 1/2 < 0 > 1 < \infty]$

pour la fonction Dirac: $\delta(x) = \Delta[-\infty > 0 < 0 > \infty < 0 > 0 < \infty]$

pour la fonction Haar: $\psi_H(x) = \Psi[0^+ > 1 < 1/2^+ > -1 < 1]$

Pour toute fonction réelle $g(x) : I \to \mathbf{Im}(g)$, où $I = [x_1, x_2]$, pour laquelle $g(x_1) = g(x_2) = a$, on peut construire, une **extension périodique** sur $\mathbf{R}$: $g_p(x_R) = \sum_k g_{pk}(x_R, k)$, tel que :

$g_{pk}(x_R, k) = g_{pk}(x + kT) = g(x_R - kT) = g(x)$ si $x_R \in [x_1 + kT, x_2 + kT]$ et (1)
$g_{pk}(x_R) = 0$, si $x_R \notin [x_1 + kT, x_2 + kT]$, pour $\forall k \in \mathbf{Z}$ et $\forall x \in [x_1, x_2]$. Ici $T = x_2 - x_1$
Si $k = 0$, on obtient pour $x \in [x_1, x_2]$: $g_{p0}(x, 0) = g(x)$.

Pour une définition dans laquelle la valeur de $k$ (dépendant de $x_R$), apparaît implicitement, on peut appeler à la **fonction partie entière (par défaut)** $E(x) = \lfloor x \rfloor =$ le plus grand entier inférieur ou égal à $x$: $E(x) \leq x < E(x) + 1$ (*floor function*, dans la littérature anglo-saxonne) [1]. Pour tout $x_R \in \mathbf{R}$, on définit la fonction $K(x_R) = E((x_R - x_1)/T)$ (1a)
Alors, pour $\forall x_R \in R$, $\exists x \in [x_1, x_2]$, $x_R = x + KT$. Par définition: $g(x_R) = g(x + KT) = g(x)$.

Pour la fonction $g_P(x) = \sin(x)$, définie sur l'axe $\mathbf{R}$, la relation (1) est vrai, sous la forme $\sin(x_R) = \sin(x)$, pour $\forall x_R = x + 2\pi k$, **implicitement**, **simultanément**, pour toute $k \in \mathbf{Z}$ et pour tous les intervalles $[(2k-1)\pi, (2k+1)\pi]$ que leurs correspondent. Pour une certain fonction $g(x):[x_1, x_2]$, en dehors de cet intervalle $g_p(x)$ doit être spécifiée **explicitement**, par des translations **successives**, pour tous les intervalles de validité $x_1 + kT < x < x_1 + (k+1)T$: $g_P(x_R) = g(x_R - kT) = g(x)$, pour toute $k \in \mathbf{Z}$, ou **implicitement**, **simultanément**:
$g_P(x_R) = g[x_R - T \cdot E((x_R - x_1)/T)] = g(x)$. (1b)

Selon la théorie développée par Fourier, la fonction $2L$-périodique $f_p(x):(-\infty, \infty)$, l'extension sur l'axe réel de la fonction de carré intégrable $f(x):[-L, L]$, peut être décomposée en une somme: $\bar{f}(x) = f_0 + \sum_{n=1}^{\infty} [a_n \cos(\omega_n x) + b_n \sin(\omega_n x)]$ où, pour toute $n \in \mathbf{N}$, $\omega_n = n\frac{\pi}{L}$. Ici, $f_0$ est la valeur moyenne de la fonction $f(x)$, sur l'intervalle $[-L, L]$, $\cos(\omega_n x)$ et $\sin(\omega_n x)$ sont des fonctions continues (appelées harmoniques secondaires paires unitaires, respectivement harmoniques secondaires impaires unitaires), qui proviennent des fonctions continues $\cos(\omega_0 x)$ et $\sin(\omega_0 x)$, pour $\omega_0 = \pi/L$ (appelées harmonique fondamentale paire unitaire, respectivement harmonique fondamentale impaire unitaire) par la multiplication de leurs argument avec un entier naturel positifs $n \in \mathbf{N}^+$. Les fonctions de valeur moyenne nulle, $\sin(n\omega_0 x)$ et $\cos(n\omega_0 x)$ prennent sur l'intervalle $[-L/n, L/n]$, les mêmes valeurs que les harmoniques fondamentales les prennent sur l'intervalle $[-L, L]$ et ils satisfont implicitement:
$\sin(n\omega_0 x) = \sin[n(\omega_0 x + 2\pi k)]$ et $\cos(n\omega_0 x) = \cos[n(\omega_0 x + 2\pi k)]$, pour tous les intervalles $[(2k-1)\omega_0/n, (2k+1)\omega_0/n]$, où $k$ est un entier relatif.

Soit la fonction $g(x)$ de variable réelle $x \in [-L, L]$, qui a des valeurs réelles finies dans l'ensemble image $\mathbf{Im}(g)$ et vérifie l'égalité $g(-L) = g(L)$. Semblable aux sinusoïdes, à partir de la fonction $g_P(x) \to (-\infty, \infty)$, qui est l'extension $2L$-périodique sur l'axe réel de la fonction $g(x)$, on peut obtenir pour chaque $n \in \mathbf{N}^+$, par dilatation, une fonction $2L/n$-périodique:

$g_n(x) = g_p(nx) = \sum_{k=-\infty}^{\infty} g_{nk}(nx, k): (-\infty, \infty)$, où pour tout entier relatif $k$, $g_{nk}(nx, k)$ est une fonction définie sur l'intervalle $[(2k-1)L/n, (2k+1)L/n]$. Dans cet intervalle, $g_{nk}(nx, k)$ prend les mêmes valeurs que celles que les prend $g_{n0}(nx, 0) = g(nx)$ sur l'intervalle $[-L/n, L/n]$ et $g(x)$ sur l'intervalle $[-L, L]$. Explicitement et successivement: $g_{nk}(nx, k) = g(nx + 2kL) = g_{n0}(nx, 0) = g(nx)$ pour $x \in [(2k-1)L/n, (2k+1)L/n]$ et $g_{nk}(nx, k) = 0$, pour $x \notin [(2k-1)L/n, (2k+1)L/n]$. La relation de périodicité devient: $g_n(x) = g_n(x + 2kL/n)$, pour tout $k \in Z$, ou implicitement:
$g_n(x) = g_n(x + 2L \cdot E(n(x+L)/2L))$. La fonction $g(x)$ et les fonctions $g_n(x)$ ont, sur l'intervalle $[-L, L]$, la même valeur moyenne $g_0$.

Nous appellerons la fonction $g_n(x)$, restreinte à l'intervalle $[-L, L]$, **la g-harmonique de l'ordre $n$** de la fonction $g(x)$ et la fonction $g_1(x)=g(x)$, **la g-harmonique fondamental.** Nous allons également, introduire pour la g-harmonique de l'ordre $n$ une notation réduite:
$$g_n(x)=G[-L/n<g(nx)>L/n]_n, \ n \in \mathbf{N}^+. \tag{1c}$$

Ces opérations de translation et de dilatation sont similaires à celles utilisées pour créer les fonctions ondelettes $\psi_{nk}(x)$, à partir d'une fonction *mère $\Psi(x)$* [3]:
$\psi_{nk}(x)=k\Psi[(x-b)/a]$, pour $b=2k/n$ et $a=1/n$.

Si la fonction $g(x)$ a, dans l'intervalle $[x_1, x_{m+1}]$, un nombre fini $m$ de discontinuités, la fonction $g_n(x)$ (la g-harmonique de l'ordre $n$) aura un nombre $m \cdot n$ de telles discontinuités, qui tend vers l'infini si $n \to \infty$. Pour cette raison, la fonction $g(x)$ est impropre pour générer une base pour un sous-espace de fonctions. Mais, si la fonction $g(x)$ est de carré intégrable (appartiennent à l'espace $L^2[x_1, x_{m+1}]$), ou si elle satisfait les conditions aux limites de Dirichlet, elle peut être développé en une série [1, 9]:

$$\bar{g}(x) = g_0 + \sum_{n=1}^{\infty}[a_n \cos(\omega_n x) + b_n \sin(\omega_n x)], \text{ où } \omega_n=n\pi/L, \ \forall n \in \mathbf{N}.$$

Ici, $a_n = \frac{1}{L}\int_{x_1}^{x_{m+1}} g(x)\cos(\omega_n x)dx$ et $b_n = \frac{1}{L}\int_{x_1}^{x_{m+1}} g(x)\sin(\omega_n x)dx$; $g_0 = \frac{1}{L}\int_{x_1}^{x_{m+1}} g(x)\cos dx$

Parce que $\bar{g}(x)$ est une série convergente des fonctions continues, elle est une fonction continue (une Fourier-fonction $\bar{g}_F(x)$) et peut être pris en compte pour générer une base pour les fonctions de l'espace $L^2[x_1, x_{m+1}]$. Évidemment, toutes les fonctions continues sont des F-fonctions. En tous les points de continuités, $\bar{g}(x) \to g_F(x)$, alors que dans les environs d'un point de discontinuité, $\bar{g}(x-h) \to \lim_{x \to x_i^-} g_F(x)$, et $\bar{g}(x+h) \to \lim_{x \to x_i^+} g_F(x)$, si $h \to 0$. Dans l'intervalle $[x_i-h, x_i+h]$, pour $h \to 0$, la fonction $\bar{g}(x)$ approche la droit $g_F(x)=x[g_F(x_i+h)+g_F(x_i-h)]/2h$ et $\bar{g}(x_i)$ approche la valeur $[g_F(x_i+h)+g_F(x_i-h)]/2$. En consequence, toutes les g-harmoniques $g_n(x)$, fondamental ou secondaires sont des fonctions continues dans tout l' intervalle *I.*

Une F-fonction peut être construite par définition: soit, une fonction $g(x)$, définie dans l'intervalle $[x_1, x_2]$, avec une discontinuité de saut au point $x_i$. La F-fonction correspondante est: $g_F(x) = \lim_{h \to 0} G\left[x_1 * (g) * x_d - h * \left(x \frac{g(x_d - h) + g(x_d + h)}{2h}\right) * x_d + h * (g) * x_2\right]$, $h$ réel. (1d)

Comme nous l'avons déjà mentionné, dans cet article, lorsque nous analyserons des g-harmoniques $g_n(x)$, continues par morceaux, nous allons considérer toujours (sauf les cas spécifiés expressément) qu'il s'agit de fonctions $\bar{g}_F(x)$ de Fourier.

Pour les phénomènes de la Nature, décrits par l'évolution de certaines fonctions, du moins pour des considérations énergétiques, les fonctions discontinues cèdent leur place aux fonctions qui approchent les fonctions de Fourier.

C'est évident que les fonctions $\bar{g}_n(x)$, $n=1, 2, ..., \infty$ sont, deux à deux, **indépendantes**. Par conséquent, ils forment une base génératrice d'un sous-espace de $L^2$. Nous appellerons cette base: **la base générée par $g(x)$** ou, plus simplement, **la base $g(x)$**, notée ***$B_g$***.

### 2. Séries de Fourier périodiques non sinusoïdales

Dans la section précédente, nous avons constaté l'existence d'une quelconque analogie formelle entre les fonctions réelles et finies $cos(\omega_0 x), sin(\omega_0 x)$, définies sur l'intervalle $[-L,L]$, et les autres fonctions réelles et finies $g(x)$, définies sur le même intervalle. Dans cette section, nous allons essayer de découvrir ceux catégories des fonctions $g(x)$ qui accentuent cette analogie, de sorte qu'elle devienne une analogie fonctionnelle, utile pour créer des bases de fonctions complètes et indépendantes.

Nous utiliserons les notations $\bar{f}$, $\hat{f}$ et $\widetilde{f}$, pour les développements en série de Fourier sinusoïdale, en série de Fourier non sinusoïdale, respectivement en série de Fourier non sinusoïdale orthogonale. Pour les formules de développement des fonctions en série de Fourier et pour leures propriétés, nous avons consulté des travaux réputés [4-13].

## 2.1. Séries de Fourier périodiques non sinusoïdales des fonctions paires

**Théorème 1**. *La base $B_g$ d'une fonction paire $g(x)$ définie sur l'intervalle $[-L, L]$ d'une espace $L^2$ (notée $L^2[-L, L]$), ayant la valeur moyenne $g_0$ nulle sur cet intervalle, constitue une base complète pour le système $F_E$ des toutes les fonctions paires $f_e(x)$, réelles, de $L^2$, périodique de période 2L, ayant la valeur moyenne nulle sur cet intervalle.*

La démonstration de ce théorème inclut également, la manière de calculer les coefficients $A_n$ du développement de Fourier non sinusoïdale de la fonction paire $f_e(x)$ de $L^2$:

$$\hat{f}_e(x) = \sum_{n=1}^{\infty} A_n \bar{g}_n(x), \quad \text{où} \quad \bar{g}_n(x) \text{ sont des séries de Fourier} \tag{2}$$

La fonction $f_e(x)$ qui est, par définition, de valeur moyenne nulle sur l'intervalle $[-L, L]$, peut être développée, selon la thèse de Fourier, de manière univoque, en une somme infinie de fonctions cosinus paires:

$$\bar{f}_e(x) = \sum_{n=1}^{\infty} a_n \cos(\omega_n x), \quad \text{où} \quad \omega_n = n\frac{\pi}{L} = n\omega_0. \tag{2.1}$$

Dans le même temps, toutes les quasi-harmoniques $g_n(x)$ peuvent être écrites comme une combinaison linéaire de la fonction $\cos\omega_n x$ et des autres cosinusoïdes de rang supérieur:

$$\bar{g}_1(x) = c_1 \cos\omega_0 x + c_2 \cos 2\omega_0 x + c_3 \cos 3\omega_0 x + c_4 \cos 4\omega_0 x + ...$$
$$\bar{g}_2(x) = c_1 \cos 2\omega_0 x + c_2 \cos 4\omega_0 x + c_3 \cos 6\omega_0 x + c_4 \cos 8\omega_0 x + ...$$
$$\bar{g}_3(x) = c_1 \cos 3\omega_0 x + c_2 \cos 6\omega_0 x + c_3 \cos 9\omega_0 x + c_4 \cos 12\omega_0 x + ...$$
........................
$$\bar{g}_n(x) = c_1 \cos n\omega_0 x + c_2 \cos 2n\omega_0 x + c_3 \cos 3n\omega_0 x + c_4 \cos nN\omega_0 x + ...$$
........................

A partir de ces relations, pour le cas général, on obtient, pour $c_1 \neq 0$:

$$\cos\omega_n x = (\bar{g}_n - c_2 \cos 2\omega_n x - c_3 \cos 3\omega_n x - c_4 \cos 4\omega_n x - ....)/c_1, \text{ pour } n=1, 2, ..., \infty$$

Ici, toutes les fonctions $\bar{g}_n(x)$ sont des F-fonctions (donc, continues)

$$\bar{f}_e(x) = a_1 \cos\omega_0 x + a_2 \cos 2\omega_0 x + a_3 \cos 3\omega_0 x + a_4 \cos 4\omega_0 x + ... + a_i \cos i\omega_0 x + ... + a_n \cos n\omega_0 x + ... =$$

$$= \frac{a_1}{c_1}(\bar{g}_1 - c_2 \cos 2\omega_0 x - c_3 \cos 3\omega_0 x - ...) + \frac{a_2}{c_1}(\bar{g}_2 - c_2 \cos 4\omega_0 x - c_3 \cos 6\omega_0 x - ...) +$$

$$+ \frac{a_3}{c_1}(\bar{g}_3 - c_2 \cos 6\omega_0 x - c_3 \cos 9\omega_0 x - ...) + \frac{a_4}{c_1}(\bar{g}_4 - c_2 \cos 8\omega_0 x - c_3 \cos 12\omega_0 x - ...) + ... +$$

$$+ \frac{a_n}{c_1}(\bar{g}_n - c_2 \cos 2n\omega_0 x - c_3 \cos 3n\omega_0 x - ...) =$$

$$= \hat{f}_e(x) = A_1 \bar{g}_1(x) + A_2 \bar{g}_2(x) + A_3 \bar{g}_3(x) + A_4 \bar{g}_4(x) + A_5 \bar{g}_5(x) + ... = \sum_{n=1}^{\infty} A_n \bar{g}_n(x)$$

L'égalité de $\bar{f}_e(x)$ avec $\hat{f}_e(x)$ est univoque, ce qui devait être prouvé. Donc: (2.2)

$A_1 = \dfrac{a_1}{c_1} = K_1$, $A_2 = K_1\left(\dfrac{a_2}{a_1} - \dfrac{c_2}{c_1}\right)$, $A_3 = K_1\left(\dfrac{a_3}{a_1} - \dfrac{c_3}{c_1}\right)$, $A_4 = K_1\left(\dfrac{a_4}{a_1} - \dfrac{a_2}{a_1}\dfrac{c_2}{c_1} - \dfrac{c_4}{c_1} + \dfrac{c_2^2}{c_1^2}\right)$, $A_5 = K_1\left(\dfrac{a_5}{a_1} - \dfrac{c_5}{c_1}\right)$

$A_6 = K_1\left(\dfrac{a_6}{a_1} - \dfrac{a_2}{a_1}\dfrac{c_3}{c_1} - \dfrac{a_3}{a_1}\dfrac{c_2}{c_1} - \dfrac{c_6}{c_1} + 2\dfrac{c_2 c_3}{c_1^2}\right)$, $A_7 = K_1\left(\dfrac{a_7}{a_1} - \dfrac{c_7}{c_1}\right)$,

$A_8 = K_1\left(\dfrac{a_8}{a_1} - \dfrac{a_2}{a_1}\dfrac{c_4}{c_1} - \dfrac{a_4}{a_1}\dfrac{c_2}{c_1} - \dfrac{c_8}{c_1} + 2\dfrac{c_2 c_4}{c_1^2} + \dfrac{a_2}{a_1}\dfrac{c_2^2}{c_1^2} - \dfrac{c_2^3}{c_1^3}\right)$, $A_9 = K_1\left(\dfrac{a_9}{a_1} - \dfrac{a_3}{a_1}\dfrac{c_3}{c_1} - \dfrac{c_9}{c_1} + \dfrac{c_3^2}{c_1^2}\right)$,

$A_{10} = K_1\left(\dfrac{a_{10}}{a_1} - \dfrac{a_2}{a_1}\dfrac{c_5}{c_1} - \dfrac{a_5}{a_1}\dfrac{c_2}{c_1} - \dfrac{c_{10}}{c_1} + \dfrac{2c_2 c_5}{c_1^2}\right)$, $A_{11} = \dfrac{a_{11}c_1 - a_1 c_{11}}{c_1^2} = K_1\left(\dfrac{a_{11}}{a_1} - \dfrac{c_{11}}{c_1}\right)$,

$A_{12} = \dfrac{a_{12}c_1 - a_2 c_6 - a_3 c_4 - a_4 c_3 - a_6 c_2 - a_1 c_{12}}{c_1^2} + \dfrac{2a_1 c_3 c_4 + 2a_1 c_2 c_6 + 2a_2 c_2 c_3 + a_3 c_2^2}{c_1^3} - \dfrac{3a_1 c_2^2 c_3}{c_1^4} =$

$= K_1\left(\dfrac{a_{12}}{a_1} - \dfrac{a_2}{a_1}\dfrac{c_6}{c_1} - \dfrac{a_3}{a_1}\dfrac{c_4}{c_1} - \dfrac{a_4}{a_1}\dfrac{c_3}{c_1} - \dfrac{a_6}{a_1}\dfrac{c_2}{c_1} - \dfrac{c_{12}}{c_1} + 2\dfrac{c_3 c_4 + c_2 c_6}{c_1^2} + 2\dfrac{a_2}{a_1}\dfrac{c_2 c_3}{c_1^2} + \dfrac{a_3}{a_1}\dfrac{c_2^2}{c_1^2} - 3\dfrac{c_2^2 c_3}{c_1^3}\right)$, ...

En conclusion, parce que toute fonction paire $f_e(x)$ du sous-espace $L^2[-L, L]$, peut être développée en série de Fourier sinusoïdale (2.1), elle peut également être développée en série de Fourier non sinusoïdale (2). Pour calculer les coefficients de ce développement, il est nécessaire de connaître les coefficients $a_n$ du développement de Fourier de la fonction $f_e(x)$ ainsi que $c_n$, les coefficients de la fonction $g(x)$, ce qui implique le calcul des intégrales $\int_{-L}^{L} f_e(x)\cos\omega_n x\, dx$, respectivement $\int_{-L}^{L} g(x)\cos\omega_n x\, dx$. Bien sûr, pour une autre fonction paire $f(x)=f_0(x)+f_e(x)$, $f_0(x)\neq 0$, $\widehat{f}(x)=f_0(x)+\sum_{n=1}^{\infty} A_n \bar{g}_n(x)$. L'aproximation d'ordre $N$ s'écrit :

$\widehat{f}_N(x) = f_0(x) + \sum_{n=1}^{N} A_n \bar{g}_n(x)$, où $\bar{g}_n(x) = \sum_{m=1}^{N} c_m \cos(m\omega_n x)$. Si $N\to\infty$, $\bar{g}_n(x)\to g(x)$

Pour illustrer la méthode de calcul, soit la fonction $f_2(x)=G[-1^+>x^2<1^-]$, qui est un polynôme du second degré, sans des discontinuités, et qui a le développement de Fourier:

$$\bar{f}_2(x) = f_0 + \sum_{n=1}^{\infty} a_n \cos n\pi x = \frac{1}{3} + \sum_{n=1}^{\infty} \frac{4(-1)^n}{n^2 \pi^2} \cos n\pi x \tag{2a}$$

Nous voulons de la développer en une base générée par la fonction paire (impulsions rectangulaires de valeur moyenne nulle) $g_e = g^{dr} = G_e[-1>-1<-1/2>1<1/2>-1<1]$ :

$$\widehat{f}_2(x) = \frac{1}{3} + \sum_{n=1}^{\infty} A_n \bar{g}_n^{dr}(x), \tag{2b}$$

Le développement en série trigonométrique de la fonction $g_e(x)$ est:

$$\bar{g}_e(x) = \sum_{n=1}^{\infty} c_n \cos(2n-1)\pi x = -\frac{2}{\pi}\sum_{n=1}^{\infty} \frac{(-1)^n \cos(2n-1)\pi x}{2n-1} \tag{2c}$$

Les relations (2a) et (2c), fournissent les coefficients suivants:
$a_1 = -4/\pi^2$, $a_2 = 1/\pi^2$, $a_3 = -4/9\pi^2$, $a_4 = 1/4\pi^2$, $a_5 = -4/25\pi^2$, $a_6 = 1/9\pi^2$, $a_7 = -4/49\pi^2$, $a_8 = 1/16\pi^2$, $a_9 = -4/81\pi^2$, $a_{10} = 1/25\pi_2$, $a_{11} = -4/121\pi^2$, $a_{12} = 1/36\,\pi^2$, ... et
$c_1 = 2/\pi$, $c_2 = 0$, $c_3 = -2/3\pi$, $c_4 = 0$, $c_5 = 2/5\pi$, $c_6 = 0$, $c_7 = -2/7\pi$, $c_8 = 0$, $c_9 = 2/9\pi$, $c_{10} = 0$, $c_{11} = -2/11\pi$, $c_{12} = 0$, ...
Selon (2.2), les coefficients du développement (2b) sont:
$A_1 = -2/\pi$, $A_2 = 1/2\pi$, $A_3 = -8/9\pi$, $A_4 = 1/8\pi$, $A_5 = 8/25\pi$, $A_6 = 2/9\pi$, $A_7 = -16/49\pi$, $A_8 = 1/32\pi$, $A_9 = -8/81\pi$, $A_{10} = -2/25\pi$, $A_{11} = -24/121\pi$, $A_{12} = 1/18\,\pi$, ...
La représentation des quasi-harmoniques correspondantes et des sommes partielles résultantes est donnée à la figure 1. Ici, nous avons représenté les fonction $\bar{g}_{en}(x)$ au lieu des fonctions $\bar{g}_{en}(x)$. Puisque la fonction $g_e(x)$ a deux points de discontinuité, les sommes partielles $S_N$ du développement non sinusoïdal, ont des points de saut en nombre croissant, à mesure que le rang $N$ augmente.

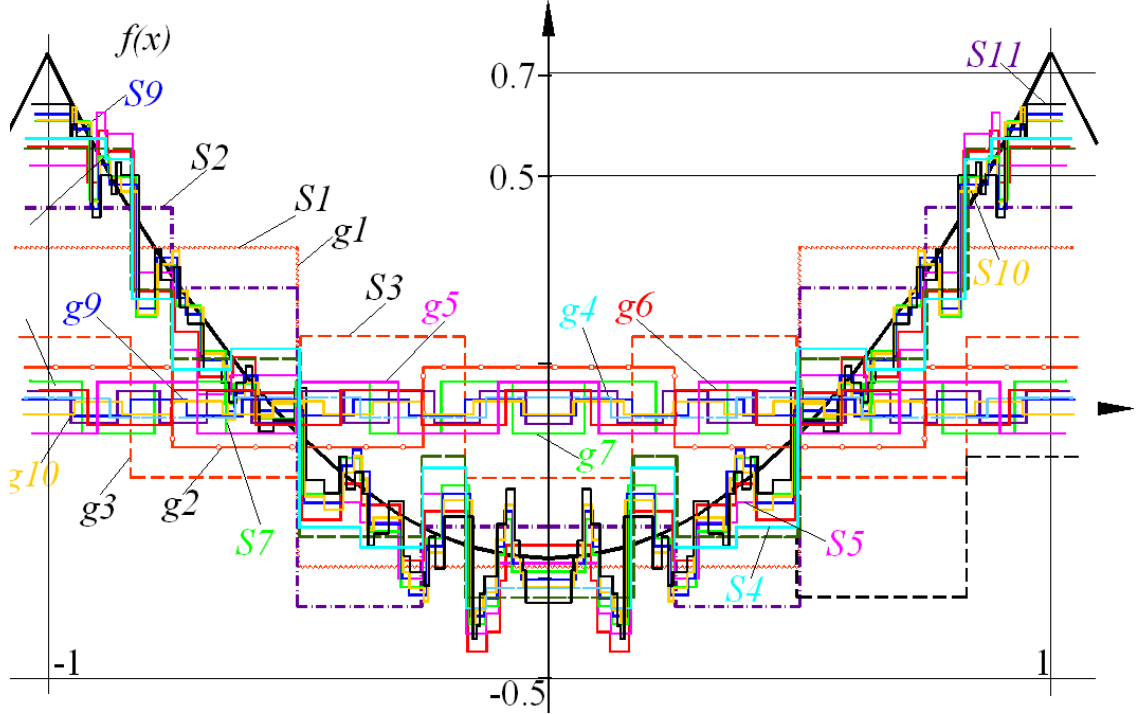

*Fig.1. L'approximation de la fonction $x^2-1/3$ par une somme d'impulsions rectangulaires gi: composante i des sommes partielles, Si: somme partielle de l'ordre i*

On peut remarquer une basse vitesse de convergence, par rapport à la méthode classique d'approximation de la même courbe, par le biais de segments des droits horizontaux.

Examinons, pour une comparaison, le développement de la même fonction par une somme des fonctions-rampe continues. Soit la même fonction $f_2(x)=F_2[-1^+>x^2<1^-]$ et la fonction $g(x)=G[-1^+>x-1/2<0>x-1/2<1^-]$. $\bar{f}_2(x)=\frac{1}{3}+\sum_{n=1}^{\infty}\frac{4(-1)^n}{n^2\pi^2}\cos n\pi x$ fournit:

$a_1=-4/\pi^2$, $a_2=1/\pi^2$, $a_3=-4/9\pi^2$, $a_4=1/4\pi^2$, $a_5=-4/25\pi^2$, $a_6=1/9\pi^2$, $a_7=-4/49\pi^2$, $a_8=1/16\pi^2$, $a_9=-4/81\pi^2$, $a_{10}=1/25\pi_2$, $a_{11}=-4/121\pi^2$, $a_{12}=1/36\,\pi^2$, ...

et $\bar{g}(x)=-4\sum_{n=1}^{\infty}\frac{\cos(2n-1)\pi x}{(2n-1)^2\pi^2}$ fournit: (2d)

$c_1=-4/\pi^2$, $c_2=0$, $c_3=-4/9\pi^2$, $c_4=0$, $c_5=-4/25\pi^2$, $c_6=0$, $c_7=-4/49\pi^2$, $c_8=0$, $c_9=-4/81\pi^2$, $c_{10}=0$, $c_{11}=-4/121\pi^2$, $c_{12}=0$, ...

A l'aide des relations (2.2), nous pouvons développer la fonction $f_2(x)$ en une série infinie de signal triangle: $\widehat{f}_2(x)=\frac{1}{3}+\sum_{n=1}^{\infty}A_n\bar{g}_n(x)$, où:

$A_1=1$, $A_2=-1/4$, $A_3=0$, $A_4=-1/16$, $A_5=0$, $A_6=0$, $A_7=0$, $A_8=-1/64$, $A_9=0$, $A_{10}=0$, $A_{11}=0$, $A_{12}=0$, ...

$$\widehat{f}_2(x)=\frac{1}{3}+\bar{g}_1(x)-\sum_{n=1}^{\infty}4^{-n}\bar{g}_{2^n}(x)=\frac{1}{3}+\left[-1>-x-\frac{1}{2}<0>x-\frac{1}{2}<1\right]-$$

$$-\sum_{n=1}^{\infty}\frac{1}{4^n}\left[-\frac{1}{2^n}>-2^n x-\frac{1}{2}<0>2^n x-\frac{1}{2}<\frac{1}{2^n}\right]_{2^n}$$

La figure 2 montre les premiers quasi-harmoniques et les premières sommes partielles

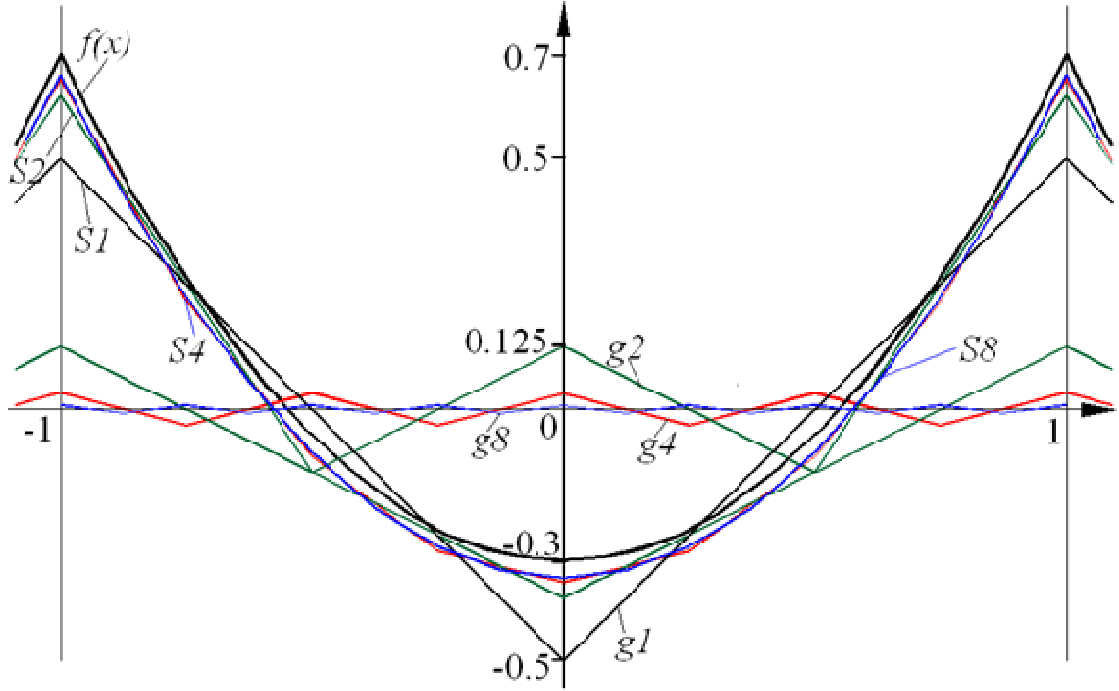

*Fig.2. L'approximation de la fonction $x^2 - 1/3$ par une somme de signal triangle paires
gi: quasi-harmoniques d'ordre i; Si: sommes partielles*

Les fonctions $f_e(x)$ et $g(x)$ étant tous les deux paires, on peut faire aussi le développement en le sens inverse: $\hat{g}(x) = \sum_{n=1}^{\infty} A'_n \bar{f}_{en}(x)$, pour $a_1 \neq 0$, où :

(2.3)
$$A'_1 = \frac{c_1}{a_1}, \ A'_2 = \frac{c_2 a_1 - c_1 a_2}{a_1^2}, \ A'_3 = \frac{c_3 a_1 - c_1 a_3}{a_1^2}, \ A'_4 = \frac{c_4 a_1 - c_2 a_2 - c_1 a_4}{a_1^2} + \frac{c_1 a_2^2}{a_1^3}, \ A'_5 = \frac{c_5 a_1 - c_1 a_5}{a_1^2}, \ ...$$

Si la fonction $f_e(x)$ est même la fonction $f = \cos\omega_0 x$, on écrit:

$$\cos\omega_0 x = \sum_{n=1}^{\infty} A_n \bar{g}_n(x) = A_1 \bar{g}_1(x) + A_2 \bar{g}_2(x) + A_3 \bar{g}_3(x) + ... =$$
$$= A_1(c_1 \cos\omega_0 x + c_2 \cos 2\omega_0 x + c_3 \cos 3\omega_0 x + ...) + A_2(c_1 \cos 2\omega_0 x + c_2 \cos 4\omega_0 x + c_3 \cos 6\omega_0 x + ...) +$$
$$+ A_3(c_1 \cos 3\omega_0 x + c_2 \cos 6\omega_0 x + c_3 \cos 9\omega_0 x + ...) + A_4(c_1 \cos 4\omega_0 x + c_2 \cos 8\omega_0 x + c_3 \cos 12\omega_0 x + ...) + ...$$

Cette relation conduit au suivant système de coefficients: (2.4)

$$A_1 = \frac{1}{c_1}, A_2 = -\frac{c_2}{c_1^2}, A_3 = -\frac{c_3}{c_1^2}, \ A_4 = \frac{-c_1 c_4 + c_2^2}{c_1^3}, \ A_5 = -\frac{c_5}{c_1^2}, \ A_6 = \frac{-c_1 c_6 + 2 c_2 c_3}{c_1^3}, \ A_7 = -\frac{c_7}{c_1^2}, \ ...$$

Si la fonction paire $g(x)$ qui génère la base du développement est une fonction avec valeur moyenne $g_0 \neq 0$ sur l'intervalle $[-L, L]$:
$$\cos\omega_n x = (\bar{g}_n - g_0 - c_2 \cos 2\omega_n x - c_3 \cos 3\omega_n x - c_4 \cos 4\omega_n x - ....)/c_1, \quad \text{pour } n \in N$$
et si la fonction $f(x)$ a aussi une valeur moyenne $f_0 \neq 0$ sur cet intervalle:
$$\hat{f}(x) = f_0 + \sum_{n=1}^{\infty} A_n [\bar{g}_n(x) - g_0] \tag{2.5}$$

**2.2. Séries de Fourier périodiques non sinusoïdales des fonctions impaires**

De la même manière est traité le problème de la fonction impaire $f_o(x)$, qui appartienne à l'espace $F_O$ des fonctions impaires de $L^2[-L, L]$.

**Théorème 2**. *La base $B_h$ d'une quiconque fonction impaire $h(x)$ de $L^2[-L, L]$ constitue une base complète pour le système $F_O$ des toutes les fonctions impaires $f_o(x)$, réelles, de $L^2$, périodique de période 2L.*

Selon la thèse de Fourier, la fonction impaire $f_o(x)$ (dont la valeur moyenne sur l'intervalle $[-L, L]$ est toujours zéro) peut être développée d'une manière univoque en une somme infinie de fonctions sinusoïdales impaires:

$$\bar{f}_o(x) = \sum_{n=1}^{\infty} b_n \sin(\omega_n x), \text{ où } \omega_n = n\omega_0 = n\frac{\pi}{L}. \tag{2.6}$$

Tout autre développement de la fonction $f_o(x)$ doit être également une somme infinie de fonctions impairs: $\hat{f}_o(x) = \sum_{n=1}^{\infty} B_n \bar{h}_n(x)$, où $\bar{h}_n(x)$ sont des séries de Fourier (2.7)

Ici, $\bar{h}_n(x) \to G[-1/n < h(nx) > 1/n]_n$, où $n \in \mathbf{N}$, sont des F-fonctions $2L/n$-périodiques. Pour chacun telle fonction: $\bar{h}_n(x) = \sum_{i=1}^{\infty} d_i \sin(i\omega_n x)$, où $\omega_n = n\frac{\pi}{L}$.

Comme dans la démonstration précédente, ce système d'équations nous permet de déterminer les coefficients $B_n$ du développement: (2.8)

$$B_1 = \frac{b_1}{d_1} = K_2, \; B_2 = K_2\left(\frac{b_2}{b_1} - \frac{d_2}{d_1}\right), \; B_3 = K_2\left(\frac{b_3}{b_1} - \frac{d_3}{d_1}\right), \; B_4 = K_2\left(\frac{b_4}{b_1} - \frac{b_2}{b_1}\frac{d_2}{d_1} - \frac{d_4}{d_1} + \frac{d_2^2}{d_1^2}\right),$$

$$B_5 = K_2\left(\frac{b_5}{b_1} - \frac{d_5}{d_1}\right), \; B_6 = K_2\left(\frac{b_6}{b_1} - \frac{b_2}{b_1}\frac{d_3}{d_1} - \frac{b_3}{b_1}\frac{d_2}{d_1} - \frac{d_6}{d_1} + 2\frac{d_2 d_3}{d_1^2}\right), \; B_7 = K_2\left(\frac{b_7}{b_1} - \frac{d_7}{d_1}\right), \text{ etc.}$$

Nous pouvons conclure, que toute fonction impaire $f_o(x)$ du sous-espace $L^2[-L, L]$, pouvant être développée en série de Fourier sinusoïdale (2.6), elle peut également être développée en série de Fourier non sinusoïdale (2.7). Pour calculer les coefficients de ce développement, (ainsi que ceux du développement inverse), il est nécessaire de connaître les coefficients $b_n$ du développement de Fourier de la fonction $f_o(x)$, ainsi que celles de la fonction $h(x)$, ce qui implique le calcul des intégrales $\int_{-L}^{L} f_o(x) \sin \omega_n x \, dx$, respectivement $\int_{-L}^{L} h(x) \sin \omega_n x \, dx$.

Cette fois, pour illustrer la méthode de calcul, soit les fonctions impaires $f_o(x) = F_o[-1 > -1 < 0 > 1 < 1]$ (la fonction créneau), pour laquelle

$$\bar{f}_o(x) = \frac{4}{\pi} \sum_{n=1}^{\infty} \frac{\sin(2n-1)\pi x}{2n-1}, \text{ et } g_o(x) = G[-1 > x < 1] \text{ (la fonction dents de scie), pour laquelle}$$

$$\bar{g}_o(x) = \sum_{n=1}^{\infty} d_n \sin n\pi x = \frac{2}{\pi} \sum_{n=1}^{\infty} (-1)^{n+1} \frac{\sin n\pi x}{n}, \text{ d'où:}$$

$b_1=4/\pi, b_2=0, b_3=4/3\pi, b_4=0, b_5=4/5\pi, b_6=0, b_7=4/7\pi, b_8=0, b_9=4/9\pi, b_{10}=0, b_{11}=4/10\pi, ...,$
$d_1=2/\pi, d_2=-2/2\pi, d_3=2/3\pi, d_4=-2/4\pi, d_5=2/5\pi, d_6=-2/6\pi, d_7=2/7\pi, d_8=-2/8\pi, d_9=2/9\pi,$
$d_{10}=-2/10\pi, d_{11}=2/11\pi, d_{12}=-2/12\pi, ...$

Par conséquence, le développement $\hat{f}_o(x) = \sum_{n=1}^{\infty} B_n \bar{g}_n(x)$ aura les coefficients suivants:

$B_1 = 2, \; B_2 = 1, \; B_3 = 0, \; B_4 = 1, \; B_5 = 0, \; B_6 = 0, \; B_7 = 0, \; B_8 = 1/2, \; B_9 = 0, \; B_{10} = 1/5,$
$B_{11} = 0, \; B_{12} = 0, \; B_{13} = 0, \; B_{14} = 0, \; B_{15} = -2/15, \; ...,$

La figure 3 présente les graphiques de ces deux fonctions (**3a** et **3b**), les premières quasi-harmoniques du développement (**c**) et les premières sommes partielles (**d**). On remarque que, au fur et à mesure que $N\rightarrow\infty$, la somme $S_N(x)$ tend très lentement vers la fonction $f_o(x)$.

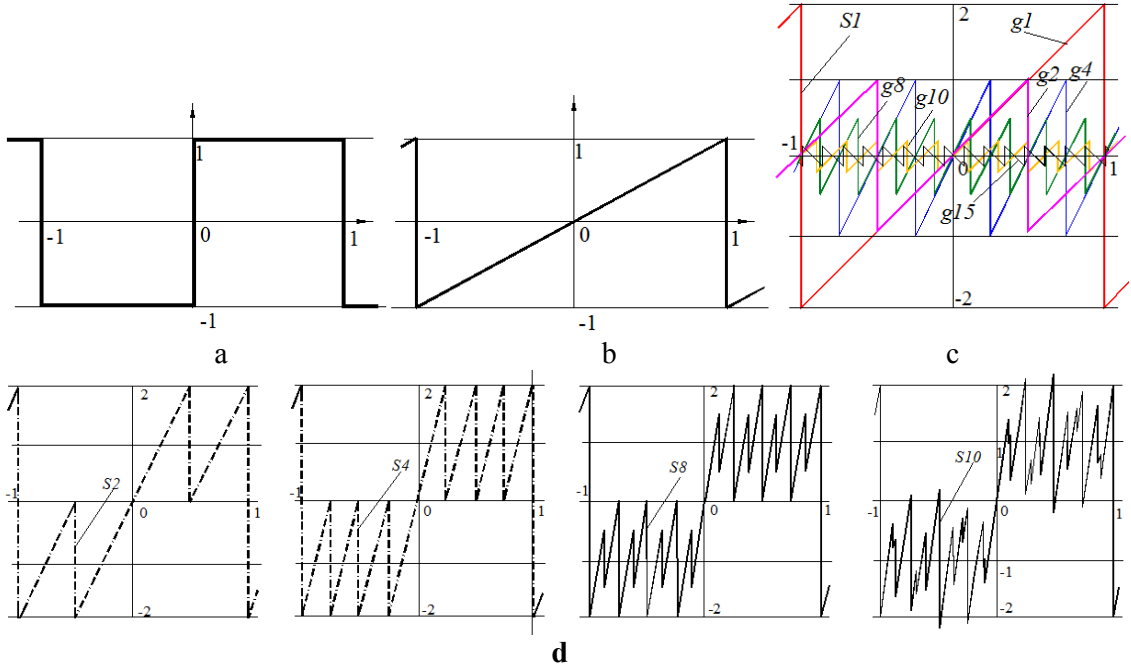

*Fig. 3. Le développement de la fonction $f_o(x)$ en la base $g_o(x)$*
**a**: $f_o(x)$ **b**: $g_o(x)$ **c**: les quasi-harmoniques $g_1(x)$, $g_2(x)$, $g_4(x)$, $g_8(x)$, $g_{10}(x)$, $g_{15}(x)$
**d**: les premières sommes partielles

Les coefficients des développements inverses sont (pour $b_1 \neq 0$): (2.9)

$$B'_1 = \frac{d_1}{b_1},\ B'_2 = \frac{d_2 b_1 - b_1 d_2}{b_1^2},\ B'_3 = \frac{d_3 b_1 - d_1 b_3}{b_1^2},\ B'_4 = \frac{d_4 b_1 - d_2 b_2 - d_1 b_4}{b_1^2} + \frac{d_1 b_2^2}{b_1^3},\ B'_5 = \frac{d_5 b_1 - d_1 b_5}{b_1^2},$$

Pour les deux fonctions précédentes $f_o(x)=F_o[-1>-1<0>1<1]$ et $g_o(x)=G[-1>x<1]$:
$B'_1 = 1/2$, $B'_2 = 1/4$, $B'_3 = 0$, $B'_4 = 1/8$, $B'_5 = 0$, $B'_6 = 0$, $B'_7 = 0$, $B'_8 = 1/16$, $B'_9 = 0$, $B'_{10} = 0$, $B'_{11} = 0$, $B'_{12} = 0$, $B'_{13} = 0$, $B'_{14} = 0$, $B'_{15} = 0$, $B'_{16} = 1/32$, ...,
et nous pouvons écrire:

$$\overline{g}_o(x) = \overline{G}[-1 < x > 1] = \frac{1}{2}\left(\overline{F}_{o1} - \sum_{n=1}^{\infty}\frac{1}{2^n}\overline{F}_{o2^n}\right) =$$

$$= \frac{L}{2}\overline{F}_o[-1>-1<0>1<1] - \frac{1}{2}\sum_{n=1}^{\infty}\frac{1}{2^n}\left[-\frac{1}{2^n} > -1 < 0 > 1 < \frac{1}{2^n}\right]_{2^n}^F$$

Si $f(x)=\sin(\omega_0 x)$, les coefficients d'un développement non sinusoïdal sont: (2.10)

$$B_1 = \frac{1}{d_1},\ B_2 = -\frac{d_2}{d_1^2},\ B_3 = -\frac{d_3}{d_1^2},\ B_4 = \frac{-d_1 d_4 + d_2^2}{d_1^3},\ B_5 = -\frac{d_5}{d_1^2},\ B_6 = \frac{-d_1 d_6 + 2 d_2 d_3}{d_1^2},\ B_7 = -\frac{d_7}{d_1^2},$$

où $d_n$ sont les coefficients du développement $\overline{h}(x) = \sum_{n=1}^{\infty} d_n \sin(\omega_n x)$, où $\omega_n = n\frac{\pi}{L}$.

Pour la fonction $\Pi_1(\theta) = \Pi[-\pi > -1 < 0 > 1 < \pi]$, pour laquelle $\Pi_n(\theta) = \Pi[-\pi/n > -1 < 0 > 1 < \pi/n]_n$, et $\overline{\Pi}_1 = \dfrac{4}{\pi} \sum_{n=1}^{\infty} \dfrac{\sin(2n-1)\theta}{2n-1}$, pour écrire le développement $\sin\theta = \hat{f}(\theta) = \sum_{n=1}^{\infty} B_n \overline{\Pi}_n(\theta)$, on utilise les formules (2.10) et on obtient:

$$B_1 = \frac{\pi}{4}, B_2 = 0, B_3 = -\frac{\pi}{12}, B_4 = 0, B_5 = -\frac{\pi}{20}, B_6 = 0, B_7 = -\frac{\pi}{28}, B_8 = 0, B_9 = 0, B_{10} = 0, \ldots$$

$$B_{2n-1} = \frac{\pi}{4} \cdot \frac{-1}{2n-1}, \; B_{2n} = 0, \text{ mais } B_{n^2} = 0, \text{ pour } n = 2, 3, \ldots, \infty$$

Tous les coefficients du développement, sauf le coefficient de la fondamentale sont négatifs.

Pour la fonction $g_1(\theta) = X_o[-\pi < -\theta - \pi > -\pi/2 < \theta > \pi/2 < -\theta + \pi > \pi]$, le développement:

$$\overline{g}_1(\theta) = \frac{4}{\pi} \sum_{n=1}^{\infty} \frac{(-1)^{n+1} \sin(2n-1)\theta}{(2n-1)^2}, \text{ fournit les coefficients:}$$

$d_1 = 4/\pi$, $d_2 = 0$, $d_3 = -4/9\pi$, $d_4 = 0$, $d_5 = 4/25\pi$, $d_6 = 0$, $d_7 = -4/49\pi$, $d_8 = 0$, $d_9 = 4/81\pi$, $d_{10} = 0$, $d_{11} = -4/121\pi$, $d_{12} = 0$, ..., pour lesquelles, pour $f(\theta) = \sin\theta$, on obtient:

$$B_1 = \frac{\pi}{4}, B_2 = 0, B_3 = \frac{\pi}{36}, B_4 = 0, B_5 = -\frac{\pi}{100}, B_6 = 0, B_7 = \frac{\pi}{196}, B_8 = 0, B_9 = 0, \ldots \quad (2.11)$$

$$B_{2n-1} = \frac{\pi}{4} \cdot \frac{(-1)^n}{(2n-1)^2}, \; B_{2n} = 0, \text{ mais } B_{n^2} = 0, \text{ pour } n = 2, 3, \ldots, \infty$$

### 2.3. Séries de Fourier périodiques non sinusoïdales de quelconques fonctions

Dans le cas général, une certaine fonction $f(x)$ de $L^2[-L, L]$, peut être écrite comme la somme de sa valeur moyenne $f_0$ sur cet intervalle, de son composant pair $f_e(x)$ (par definition, de valeur moyenne nulle sur l'intervalle $[-L, L]$) et de son composant impair $f_o(x)$: $f(x) = f_0 + f_e(x) + f_o(x)$. À la suite des deux théorèmes précédents, on peut affirmer:

**Théorème 3**. *Tout fonction $f(x)$ de $L^2[-L, L]$, peut être développé en série de Fourier non sinusoïdale en une base composée de $f_0$ (sa valeur moyenne sur cet intervalle), d'une certain base paire $B_{g-g0}$ et d'une certain base impaire $B_h$ de $L^2[-L, L]$, où la fonction paire $g(x)$ et la fonction impaire $h(x)$, sont des quelconques fonctions de $L^2[-L, L]$:*

$$\hat{f}(x) = f_0 + \sum_{n=1}^{\infty} A_n [\overline{g}_n(x) - g_0] + \sum_{n=1}^{\infty} B_n \overline{h}_n(x), \text{ où } g_0 = \int_{-L}^{L} g(x) dx$$

On peut voir que le développement de Fourier sinusoïdale est un cas particulier du développement de Fourier non sinusoïdale.

Pour illustrer, soit la fonction $f(x) = F[-1 > 0 < -1/2 > -2 < 0 > 0 < 1/2 > 2 < 1]$ qui est la somme de $f_0$ (=0), de la fonction paire $f_e = F_e[-1 > 1 < -1/2 > -1 < 1/2 > 1 < 1]$ et de la fonction impaire $f_o = F_o[-1 > -1 < 0 > 1 < 1]$, dont les développements de Fourier sont [5]:

$$\overline{f}_e(x) = \sum_{n=1}^{\infty} a_n \cos(n\pi x) = -\frac{4}{\pi} \sum_{n=1}^{\infty} \frac{(-1)^{n+1} \cos(2n-1)\pi x}{2n-1}, \text{ respectivement}$$

$$\overline{f}_o(x) = \sum_{n=1}^{\infty} b_n \sin(n\pi x) = \frac{4}{\pi} \sum_{n=1}^{\infty} \frac{\sin(2n-1)\pi x}{2n-1}$$

Les coefficients suivants sont obtenus:

$a_1 = 4/\pi$, $a_2 = 0$, $a_3 = -4/3\pi$, $a_4 = 0$, $a_5 = 4/5\pi$, $a_6 = 0$, $a_7 = -4/7\pi$, $a_8 = 0$, $a_9 = 4/9\pi$, $a_{10} = 0$, $a_{11} = -4/11\pi$,...

$b_1 = 4/\pi$, $b_2 = 0$, $b_3 = 4/3\pi$, $b_4 = 0$, $b_5 = 4/5\pi$, $b_6 = 0$, $b_7 = 4/7\pi$, $b_8 = 0$, $b_9 = 4/9\pi$, $b_{10} = 0$, $b_{11} = 4/11\pi$, ...

Pour un développement de la fonction $f(x)$ en une base exponentielle $g(x)=e^x$, étant donné que sur l'intervalle $[-1, 1]$, $g_0=\sinh 1$, nous choisirons l'ensemble constitué des fonctions $g_e(x)=\cosh x - \sinh 1$ et $g_o(x)=\sinh x$, dont les développements en série de Fourier sont:

$$\overline{g}_e(x) = \sinh 1 \sum_{n=1}^{\infty} \frac{2(-1)^n}{1+n^2\pi^2} \cos(n\pi x), \text{ et } \overline{g}_o(x) = \pi \cdot \sinh 1 \sum_{n=1}^{\infty} \frac{2n \cdot (-1)^{n+1}}{1+n^2\pi^2} \sin(n\pi x), \text{ donc:}$$

$c_1=-2\sinh 1/(1+\pi^2)$, $c_2=2\sinh 1/(1+4\pi^2)$, $c_3=-2\sinh 1/(1+9\pi^2)$,
$c_4=2\sinh 1/(1+16\pi^2)$, $c_5=-2\sinh 1/(1+25\pi^2)$
$c_6=2\sinh 1/(1+36\pi^2)$, $c_7=-2\sinh 1/(1+49\pi^2)$, $c_8=2\sinh 1/(1+64\pi^2)$, $c_9=-2\sinh 1/(1+81\pi^2)$,
$c_{10}=2\sinh 1/(1+100\pi^2)$, $c_{11}=-2\sinh 1/(1+121\pi^2)$, $c_{12}=2\sinh 1/(1+144\pi^2)$, ...
$d_1=2\pi\sinh 1/(1+\pi^2)$, $d_2=-4\pi\sinh 1/(1+4\pi^2)$, $d_3=6\pi\sinh 1/(1+9\pi^2)$, $d_4=-8\pi\sinh 1/(1+16\pi^2)$,
$d_5=10\pi\sinh 1/(1+25\pi^2)$, $d_6=-12\pi\sinh 1/(1+36\pi^2)$, $d_7=14\pi\sinh 1/(1+49\pi^2)$,
$d_8=-16\pi\sinh 1/(1+64\pi^2)$, $d_9=18\pi\sinh 1/(1+81\pi^2)$, $d_{10}=-20\pi\sinh 1/(1+100\pi^2)$,
$d_{11}=22\pi\sinh 1/(1+121\pi^2)$, $d_{12}=-24\pi\cdot\sinh 1/(1+144\pi^2)$

On résulte un développement en série de Fourier non sinusoïdale de la forme

$$\hat{f}(x) = \hat{f}_e(x) + \hat{f}_o(x) = \sum_{n=1}^{\infty} A_n [\cosh_{Fn}(x) - \sinh 1] + \sum_{n=1}^{\infty} B_n \sinh_{Fn}(x)$$

où $\cosh_{Fn}$ et $\sinh_{Fn}$ sont les extensions sur l'axe réel des F-fonctions $\cosh_F(nx)$, respectivement $\sinh_F(nx)$, définies sur les intervalles $[-1/n, 1/n])$, et les coefficients $A_n$, $B_n$ sont:

pour $K_1 = \dfrac{a_1}{c_1} = -\dfrac{2(1+\pi^2)}{\pi \cdot \sinh 1}$ et $K_2 = \dfrac{b_1}{d_1} = \dfrac{2(1+\pi^2)}{\pi^2 \cdot \sinh 1} = -\dfrac{K_1}{\pi}$ :

$A_1 = K_1$, $A_2 = K_1 \dfrac{1+\pi^2}{1+4\pi^2}$, $A_3 = -K_1 \dfrac{4}{3} \cdot \dfrac{1+3\pi^2}{1+9\pi^2}$, $A_4 = K_1 \left[ \dfrac{1+\pi^2}{1+16\pi^2} + \dfrac{(1+\pi^2)^2}{(1+4\pi^2)^2} \right]$, $A_5 = -K_1 \dfrac{4}{5} \cdot \dfrac{1-5\pi^2}{1+25\pi^2}$

$A_6 \approx -K_1 \dfrac{(1+\pi^2) \cdot (7+16\pi^2)}{36\pi^2(1+4\pi^2)}$, $A_7 = -K_1 \dfrac{8}{7} \dfrac{1+7\pi^2}{1+49\pi^2}$, $A_8 \approx K_1 \dfrac{(1+\pi^2) \cdot (7+17\pi^2+12\pi^4)}{4(1+4\pi^2)^3}$, ...

$B_1 = K_2$, $B_2 = 2K_2 \dfrac{1+\pi^2}{1+4\pi^2}$, $B_3 = -K_2 \dfrac{8}{3} \cdot \dfrac{1}{1+9\pi^2}$, $B_4 = 4K_2 \cdot \dfrac{(1+\pi^2) \cdot (2+25\pi^2+32\pi^4)}{(1+16\pi^2) \cdot (1+4\pi^2)^2}$

$B_5 = -K_2 \cdot \dfrac{24}{5} \cdot \dfrac{1}{1+25\pi^2}$, $B_6 \approx -K_2 \cdot \dfrac{3(1+\pi^2) \cdot (3+2\pi^2)}{(1+4\pi^2) \cdot (1+9\pi^2)}$, $B_7 = -K_2 \cdot \dfrac{48}{7} \cdot \dfrac{1}{1+49\pi^2}$,

$B_8 \approx K_2 \cdot \dfrac{2(1+\pi^2) \cdot (9+21\pi^2+16\pi^4)}{(1+4\pi^2)^3}$, $B_9 \approx -K_2 \cdot \dfrac{1}{9\pi^4}$, $B_{10} \approx -K_2 \cdot \dfrac{1+\pi^2}{1+4\pi^2} \cdot \dfrac{3-4\pi^2}{10\pi^2}$, ...

Dans l'intervalle $[-1/2, 1/2]$: $g_e(x)_2 = \cosh x - g_{02} = 1/2(e^x + e^{-x}) - 2\sinh(1/2)$ et $g_o(x)_2 = \sinh x$, donc $\overline{g}_e(x)_2 = \sinh\dfrac{1}{2} \cdot \sum_{n=1}^{\infty} \dfrac{4(-1)^n}{1+4n^2\pi^2} \cos(2n\pi x)$, et $\overline{g}_o(x)_2 = \pi \cdot \sinh\dfrac{1}{2} \cdot \sum_{n=1}^{\infty} \dfrac{4n \cdot (-1)^{n+1}}{1+4n^2\pi^2} \sin(2n\pi x)$

Nous pouvons notez que pour des valeurs $L<1$, la fonction $g_{es}(x)=g_e(x)_L/\sinh L$ est approximée avec des déviations acceptables par la fonction $g_{ep}(x)=x^2$, et la fonction $g_{os}(x)=g_o(x)_L/\sinh L$ est approximé avec des déviations acceptables par la fonction $g_{op}(x)=x$, les déviations étant tant petites que $L$ est plus petites.

Lorsque on demande le développement de la fonction $f(x)=f_0+f_e+f_o$ en une base générée par une fonction quelconque $g(x)=g_0+g_e(x)+g_o(x)$ de $L^2[-L, L]$, il faut trouver les coefficients $C_n$ du développement :

$$\hat{f}(x) = f_0 + \sum_{n=0}^{\infty} C_n [g_{Fn}(x) - g_0], \text{ où } g_{Fn}(x) = G_F[-L/n<g_F(nx)>L/n]_n, \; n \in \mathbf{N}. \qquad (2.11)$$

Pour simplifier, considérons le cas particulier $f_0 = g_0 = 0$:

$$\hat{f}(x) = \hat{f}_e(x) + \hat{f}_o(x) = \sum_{n=1}^{\infty} A_n g_{Fen}(x) + \sum_{n=1}^{\infty} B_n g_{Fon}(x) =$$

$$= \sum_{n=1}^{\infty} [A_n g_{Fn}(x) + (B_n - A_n) g_{Fon}(x)] = \sum_{n=1}^{\infty} [(A_n - B_n) g_{Fen}(x) + B_n g_{Fn}(x)]$$

égalité qui coïncide avec (2.11), seulement si $A_n=B_n=C_n$. Donc, aucune fonction $g_F(x)$ ne peut pas seul générer une base pour l'espace $L^2[-L, L]$ entière, exigeant l'aide d'une autre base, générée par une fonction $h_F(x)$ avec un indice de parité différent. Si on considère les identités:

$$g_{en}(x) = \frac{1}{2}[g_n(x) + g_n(-x)] \text{ et } g_{on}(x) = \frac{1}{2}[g_n(x) - g_n(-x)], \text{ on obtient, dans le cas général:}$$

$$\hat{f}(x) = f_0 + \sum_{n=1}^{\infty} \left[ \frac{A_n + B_n}{2} [g_{Fn}(x) - g_0] + \frac{A_n - B_n}{2} [g_{Fn}(-x) - g_0] \right], \text{ ou} \qquad (2.12)$$

$$\hat{f}(x) = f_0 + \sum_{n=-\infty}^{\infty} C_n [\bar{g}_{Fn}(x) - g_0]$$

En conclusion, nous pouvons formuler le théorème suivant:

**Théorème 4.** *Tout fonction $f(x)$ de $L^2[-L, L]$, peut être développé en série de Fourier non sinusoïdale, en une base composée de $f_0$ (sa valeur moyenne sur cet intervalle) et les bases générées par les fonctions $[g(x)-g_0]$ et $[g(-x)-g_0]$. Ici, $g(x)$ est n'importe quoi fonction de $L^2[-L, L]$ qui a tous ces deux composants (paire et impaire) non nulles, $g(-x)$ est aussi de $L^2[-L, L]$, et $g_0$ est la valeur moyenne de $g(x)$.*

En conséquence, le développement de la fonction $f(x)$ analysé plus tôt, sur l'intervalle $[-L, L]$, peut être effectué en une base générée par les fonctions $e^x$ et $e^{-x}$:

$$\hat{f}(x) = f_0 + \hat{f}_e(x) + \hat{f}_o(x) = f_0 + \sum_{n=1}^{\infty} A_n [\cosh_{Fn}(x) - \sinh 1] + \sum_{n=1}^{\infty} B_n \sinh_{Fn}(x) = f_0 + \sum_{n=1}^{\infty} \frac{A_n + B_n}{2} e_{Fn}^x +$$

$$+ \sum_{n=1}^{\infty} \left( \frac{A_n - B_n}{2} e_{Fn}^{-x} - A_n \sinh 1 \right) = f_0 + \sum_{n=1}^{\infty} \left[ \frac{A_n + B_n}{2} (e_{Fn}^x - \sinh 1) + \frac{A_n - B_n}{2} (e_{Fn}^{-x} - \sinh 1) \right]$$

où $e_{Fn}^x$ et $e_{Fn}^{-x}$ sont les quasi−harmoniques de l'ordre $n$ des F-fonctions $(e^x)_n$, respectivement $(e^{-x})_n$ (les extensions sur l'axe réel des F-fonctions $e^{nx}$, respectivement $e^{-nx}$, définies sur les intervalles $[-1/n, 1/n]$).

L'éventail des fonctions pouvant servir comme base pour le développement de Fourier non sinusoïdal est extrêmement large:

- si $g(x)$ est un polynôme dans $[-L, L]$, son composant paire $g_e(x)$ contient les puissances paires de $x$, pendant que son composant impaire $g_o(x)$ contient les puissances impaires
- si $g(x)$ est une fonction exponentielle, $g_{Fe}(x)$ peut être une fonction paire $G_{Fe}(\cosh x)$ et $g_{Fo}(x)$ peut être une fonction impaire $G_{Fo}(\sinh x)$
- si $g(x)$ est logarithmique: $\ln(A+x)$ (où $A>0$), $f(x)$ peut être développée uniquement sur un sous−intervalle $[a, b]$, compris dans l'intervalle $(-A, A)$, avec les bases:

$$G_{Fe}(x) = \frac{1}{2} \ln(A^2 - x^2) - g_0 \text{ et } G_{Fo}(x) = \frac{1}{2} \ln \frac{A+x}{A-x}$$

- si $g(x)$ est une fonction rationnelle de la forme $1/(A+x)$, $A>0$, la fonction $f(x)$ peut être développée sur un sous-intervalle $[a, b]$ de l'intervalle $(-A, A)$, avec les bases:

$$G_{Fe}(x) = \frac{1}{2}\left[\frac{1}{A+x} + \frac{1}{A-x}\right] - g_0 = \frac{A}{A^2 - x^2} - g_0 \text{ et } G_{Fo}(x) = \frac{1}{2}\left[\frac{1}{A+x} - \frac{1}{A-x}\right] = -\frac{x}{A^2 - x^2}$$

- si $g(x)$ est une fonction irrationnelle de la forme $\sqrt{A+x}$, $A>0$, la fonction $f(x)$ peut être développée en un sous-intervalle $[a, b]$ de l'intervalle $(-A, A)$, avec les bases:

$$G_{Fe}(x) = \frac{1}{2}[\sqrt{A+x} + \sqrt{A-x}] - g_0 \text{ et } G_{Fo}(x) = \frac{1}{2}[\sqrt{A+x} - \sqrt{A-x}]$$

## 3. Séries de Fourier périodiques quasi−sinusoïdales

Une autre façon de combiner deux fonctions $g(x)$ et $h(x)$, afin de constituer une base pour l'espace $L^2[-L, L]$, consiste à choisir la fonction $h(x)$ comme une translation de la fonction $g(x)$: $h(x)=g(x+\alpha T)$, où $\alpha \in (0, 1)$. Parmi les bases obtenues par cette modalité, il y a quelques-unes qui sont complètes et qu'ils ont une propriété spéciale: la fonction $g(x)$ a une seule composante (la paire ou la impaire), et pour $\alpha=1/4$, $h(x)$ a la parité opposée, propriété qu'ils l'ont les fonctions $sin(\omega_0 x)$, respectivement $cos(\omega_0 x)$, aussi. Parce que cette propriété est extrêmement utile pour résoudre certains problèmes pratiques, dans cette section, nous allons leur prêter certain attention.

Si $g(x)$ est une fonction de $L^2$ définie sur l'intervalle $[0, L/2]$, nous pouvons construire les fonctions $g_\&(x)$, composées de 4 segments, **chaque** segment explicitement définis avec l'aide de la fonction $g(x)$, sur un quart de l'intervalle $[-L, L]$. On impose que les fonctions $g_\&(x)$ obtenues ont la valeur moyenne nulle sur l'intervalle $[-L, L]$, et qu'ils ont des symétries internes similaires à celles des fonctions *sinus*, respectivement *cosinus*: les deux branches (ce pour $x<0$ et ce pour $x>0$) des fonctions $g_\&(x)$ impaires sont symétriques par rapport à leur **axe** médian, et les deux branches des fonctions $g_\&(x)$ paires sont symétriques par rapport au leur **point** médiane. De plus, par translation, à gauche ou à droite, avec $L/2$, une fonction de parité opposée est obtenue. Nous appellerons la fonction $g(x)$, **le noyau du développement,** et les fonctions $g_\&(x)$ dérivés, des **quasi−sinusoïdes**. Nous utiliserons les notations:
$g_s(x)=\mathbf{S}[g(x)]_L=G_s[-L<-g(x+L)>-L/2<-g(-x)>0<g(x)>L/2<g(L-x)>L]$, respectivement
$g_c(x)=\mathbf{C}[g(x)]_L=G_c[-L<-g(x+L)>-L/2<g(-x)>0<g(x)>L/2<-g(L-x)>L]$

Par exemple, si $g(x)=x$, pour $x \in [0, \pi/2]$
$g_s(x)=\mathbf{S}[x]_\pi=G_s[-\pi<-x-\pi)>-\pi/2<x>\pi/2<\pi-x>\pi]$,
$g_c(x)=\mathbf{C}[x]_\pi=G_c[-\pi<-x-\pi)>-\pi/2<-x>0<x>\pi/2<-\pi+x)>\pi]$

Les fonctions $g_s(x)$ et $g_c(x)$ satisfont toutes les conditions requises par le théorème 3, par conséquence on peut conclure:

**Théorème 5:** *Soit $g(x)$ une quelconque fonction de $L^2[0, L/2]$. Tout fonction $f(x)$ de $L^2[-L, L]$, peut être développé en série de Fourier quasi-sinusoïdale, en une base composée de $f_0$ (sa valeur moyenne sur cet intervalle) et les bases générées par les fonctions $\mathbf{S}[g(x)]_L$ et $\mathbf{C}[g(x+L/2)]_L$, ou $\mathbf{C}[g(x)]_L$ et $\mathbf{S}[g(x+L/2)]_L$.*

Dans le cas général, les quasi-sinusoïdes et/ou leurs dérivées de premier ordre présentent des discontinuités qui peuvent être supprimées en ajoutant des quasi-sinusoïdes formées des impulsions rectangulaires (un pour chaque saut) et/ou des quasi-sinusoïdes formées des fonctions-rampe correctement choisie (pour les quasi-sinusoïdes impaires avec une discontinuité en l'origine, la rampe est $-2[g(L/2)-g(0)]/L$, pendant que pour les quasi-sinusoïdes paires, la rampe a la valeur de $-dg/dx$ en le point $x=0$). Pour les fonctions obtenues, nous utiliserons le nom de **quasi-sinusoïdes lisses**, ou **presque-sinusoïdes**, ou **sinusoïdes approximatives**, respectivement **presque-cosinusoïdes**, ou **cosinusoïdes approximatives.** Ils sont particulièrement utiles dans certains problèmes pratiques (résolution des certaines équations différentielles avec des dérivées partielles).

Par exemple, pour obtenir la presque-sinusoïdes $\mathbf{S}[x^2-2x]_2$ de la Fig.4a, nous allons choisir deux fonctions paires $f_{2p}(x)=-x^2$ et $f_{2p}(x)=x^2$, pour $-1 \leq x \leq 1$, et par des translations verticales (en ajoutant des ondes rectangulaires) et horizontales (par changements de variables), nous superposons en l'origine ($x=0$), le dernier point de l'onde négative, avec le premier point de l'onde positive. La presque-cosinusoïdes $\mathbf{C}[1-x^2]_2$ s'obtient en changeant la variable $x$ avec $x-1$. Les fonctions obtenues sont définies sur l'intervalle $[-2, 2]$, ayant donc la moitié de la fréquence angulaire des fonctions initiales:

$g_s(x)_{2L} = \mathbf{S}[x^2-2x]_{2L} = G_s[-2L<-x^2-2x>0< x^2-2x>2L]$, respectivement
$g_c(x)_{2L} = \mathbf{C}[1-x^2]_{2L} = G_c[-2L<-1+(x+2)^2>-L<1-x^2>L<-1+(x-2)^2>L]$
ou, en revenant à l'intervalle de définition initial:
$g_s(x)_L = \mathbf{S}[x^2-2x]_L = G_s[-L<-4x^2-4x>0< 4x^2-4x>L]$, respectivement
$g_c(x)_L = \mathbf{C}[1-x^2]_L = G_c[-L<-1+4(x+1)^2>-L/2<1-4x^2>L/2<-1+4(x-1)^2>L]$
La fonction $g_c(x)_L$ est représentée dans la Fig.4b.

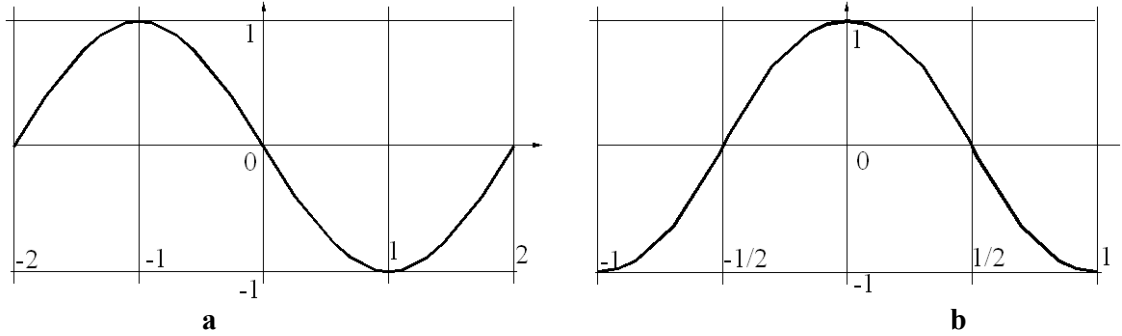

*Fig.4. Des presque−sinusoïdes de 2−ème degré*
**a**: *la fonction impaire $S[x^2-2x]_2$;* **b**: *la fonction paire $C[1-x^2]_1$*

Semblable aux développements en série de Fourier, la présence des discontinuités à l'intérieur ou aux extrémités de l'intervalle de définition de la fonction développées $f(x)$ produit, pour les développements en série non sinusoïdale, des termes supplémentaires (munies avec des coefficients de forme $a_n/n$) et un effet similaire au phénomène de Gibbs. De même, les discontinuités de la première dérivée génèrent des autres termes supplémentaires (munies avec des coefficients de forme $a_n/n^2$) et des phénomènes supplémentaires d'oscillation avec une amplitude significative. Les presque-sinusoïdes font partie de la classe $C^1$ de régularité (fonctions dont la première dérivée est continue) et, en raison des leurs propriétés de symétrie, semblable à ceux des fonctions sinus et cosinus, sont les mieux adaptés lorsque des développements en série de Fourier sont nécessaires. De même, comme dans la pratique des développements de Fourier, lorsque les autres données du problème le permettent, il est avantageux de construire pour la fonction développée $f(x)$ aussi, définie sur un intervalle $[x_1, x_2]$, une extension quasi-sinusoïdale lisse, définie sur un intervalle $[x_{1e}, x_{2e}]$, qui inclut l'intervalle de définition. Par une telle approche, les "termes résiduels" du développement sont supprimés.

Pour plus de clarté, nous allons développer en série non sinusoïdale la quasi-sinusoïde: $g_c(x)_L = C[1-x^2]_L = G_c[-L<-1+4(x+1)^2>-L/2<1-4x^2>L/2<-1+4(x-1)^2>L]$, lequel a le développement en série de Fourier: $\bar{g}_c(x) = \sum_{n=1}^{\infty}\left\{\frac{32(-1)^{n+1}}{(2n-1)^3 \pi^3}\cos(2n-1)\pi x\right\}$, pour deux bases :

1) l'onde rectangulaire (2.c): $g_e = G_e[-1>-1<-1/2>1<1/2>-1<1]$ et
2) l'onde triangulaire (2.d): $f_{12}(x) = F_{12}[-1>-x-1/2<0>x-1/2<1]$ de la section 2.1.
Leurs développements en des séries trigonométriques de Fourier sont:

$$\bar{g}_e(x) = \sum_{n=1}^{\infty} c_n \cos(2n-1)\pi x = -\frac{2}{\pi}\sum_{n=1}^{\infty}\frac{(-1)^n \cos(2n-1)\pi x}{2n-1} \quad \text{et} \quad \bar{f}_{12}(x) = -4\sum_{n=1}^{\infty}\frac{\cos(2n-1)\pi x}{(2n-1)^2 \pi^2}$$

Afin de pouvoir comparer ces développements avec ceux de la section 2.1, nous allons utiliser les coefficients du développement de Fourier de la fonction $-1/2 \cdot g_c(x)$, défini sur le même intervalle. Son développement en série de Fourier génère les coefficients:
$a_1=-16/\pi^3$, $a_2=0$, $a_3=16/27\pi^3$, $a_4=0$, $a_5=-16/125\pi^3$, $a_6=0$, $a_7=16/343\pi^3$, $a_8=0$, $a_9=-16/729\pi^3$, $a_{10}=0$, $a_{11}=16/1331\pi^3$, $a_{12}=0$, ...
Pour la fonction créneau, les relations (2.2) nous conduire à:

$A_1=-8/\pi^2$, $A_2=0$, $A_3=-64/27\pi^2$, $A_4=0$, $A_5=192/125\pi^2$, $A_6=0$, $A_7=-384/343\pi^2$, $A_8=0$, $A_9=64/729\pi^2$, $A_{10}=0$, $A_{11}=-960/1331\pi^2$, $A_{12}=0$, ...

Pour la fonction dent de scie, avec les mêmes relations, on calcule:

$A_1=4/\pi$, $A_2=0$, $A_3=-16/27\pi$, $A_4=0$, $A_5=-16/125\pi$, $A_6=0$, $A_7=-32/343\pi$, $A_8=0$, $A_9=16/729\pi$, $A_{10}=0$, $A_{11}=-48/1331\pi$, $A_{12}=0$,...

Avec les valeurs ainsi obtenues, nous pouvons construire les quasi-harmoniques secondaires et les premières sommes partielles des développements quasi-sinusoïdales correspondants: les figures 5 et 6.

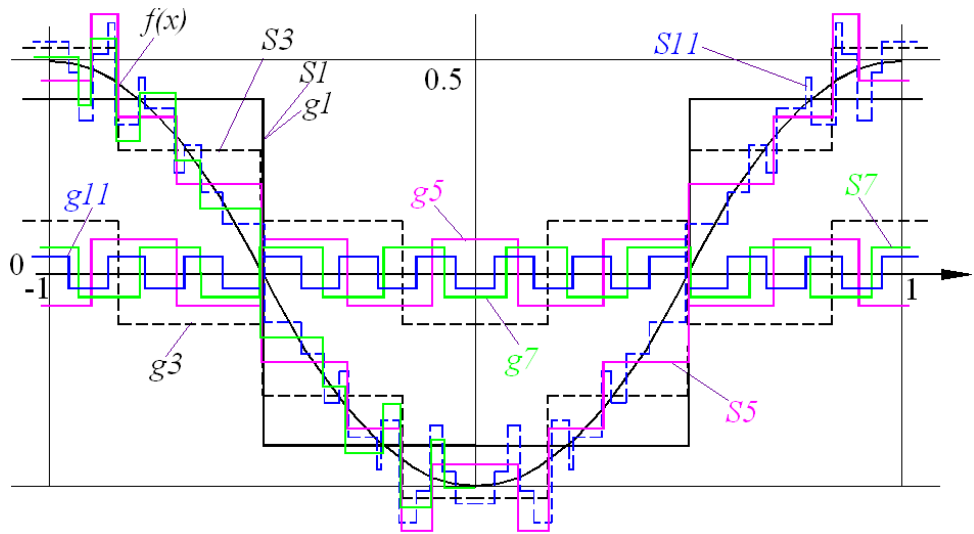

*Fig.5. L'approximation de la quasi−sinusoïde −$g_c$/2 par une somme des fonctions créneau*

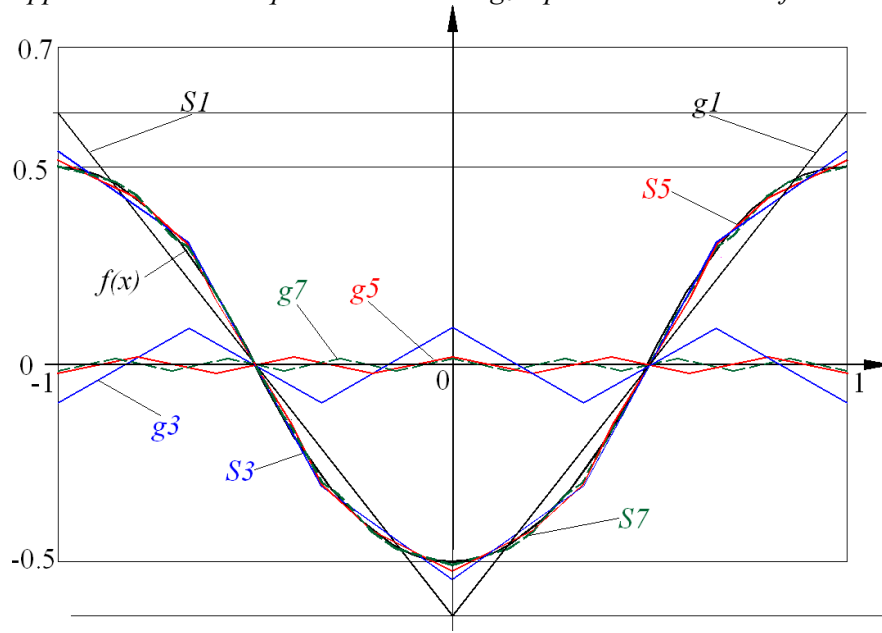

*Fig.6 L'approximation de la quasi-sinusoïde −$g_c$ /2 par une somme des fonctions dent-de-scie*

L'aspect des développements des quasi-sinusoïdes en série presque-sinusoïdale se reflète en les développements inverses des fonctions analysées précédemment. Aux figures 7 et 8, sont présentées les premières quasi-harmoniques et les sommes partielles d'ordre 12 pour les développements inverses $\hat{g}_e(x) = \sum_{n=1}^{\infty} A'_{dn}(g_c)_n$ et $\hat{f}_{12}(x) = \sum_{n=1}^{\infty} A'_{tn}(g_c)_n$. Selon (2.3):

$A'_1 = -\pi^2/8$, $A'_2 = 0$, $A'_3 = \pi^2/27$, $A'_4 = 0$, $A'_5 = -3\pi^2/125$, $A'_6 = 0$, $A'_7 = 6\pi^2/343$, $A'_8 = 0$,
$A'_9 = -\pi^2/81$, $A'_{10} = 0$, $A'_{11} = 15\pi^2/1331$, $A'_{12} = 0$, ..., respectivement:
$A'_1 = -\pi/4$, $A'_2 = 0$, $A'_3 = -\pi/27$, $A'_4 = 0$, $A'_5 = -\pi/125$, $A'_6 = 0$, $A'_7 = -2\pi/343$, $A'_8 = 0$, $A'_9 = -\pi/243$,
$A'_{10} = 0$, $A'_{11} = -3\pi/1331$, $A'_{12} = 0$, ...

    Dans les deux cas, les coefficients des développements ont des valeurs proches de celles obtenues par les développements en série sinusoïdale (Fourier). Les différences deviennent encore plus petites si le noyau du développement est remplacé par $g(x) = 1 - x^{1,75}$.

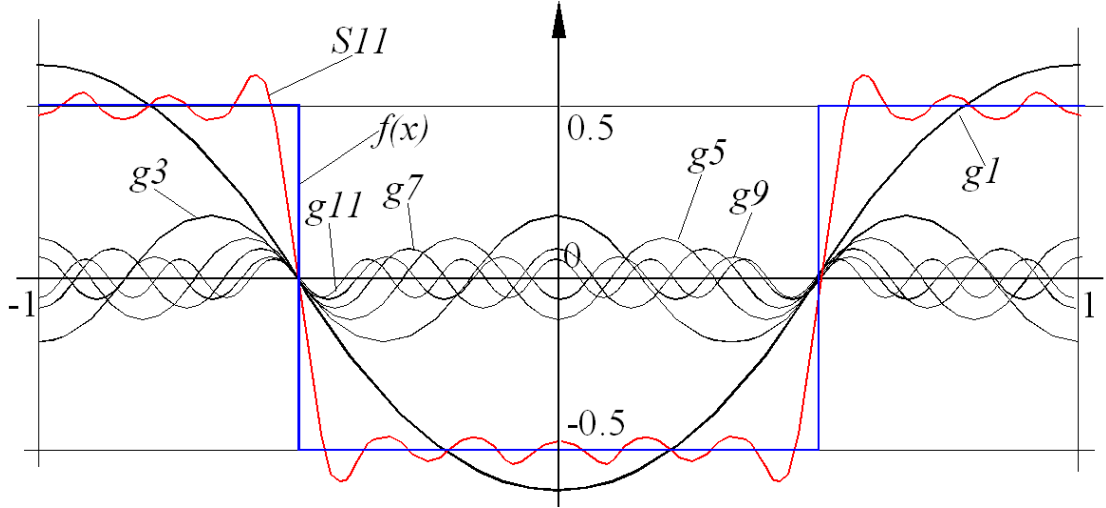

*Fig.7. L'approximation de la fonction créneau par des presque-sinusoïdes*

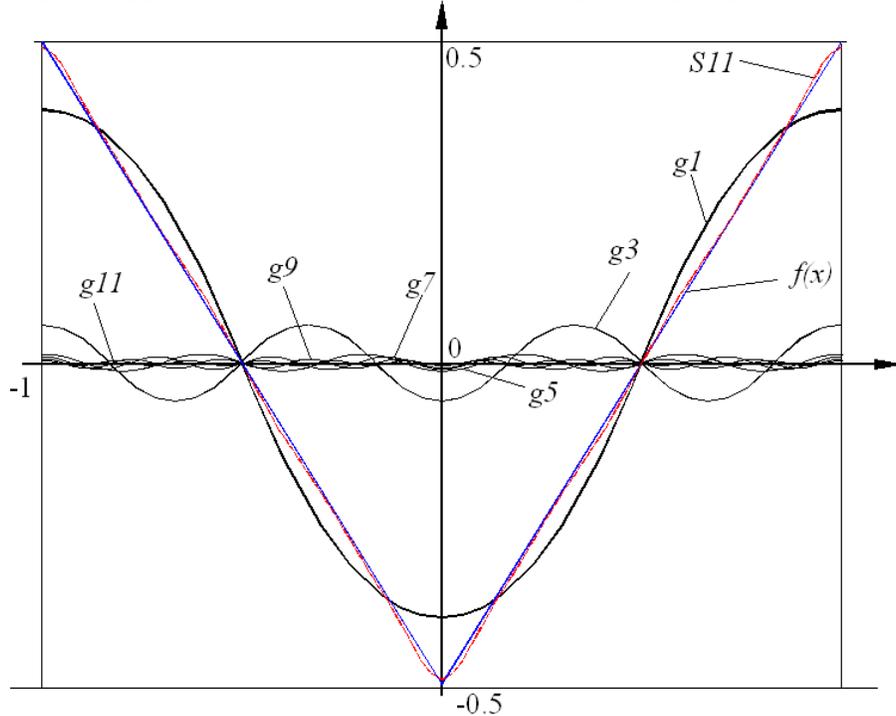

*Fig.8. L'approximation de la fonction dent de scie par des presque−sinusoïdes*

Dans l'exemple suivant, nous allons développer la fonction créneau (Fig.9.a) : $f_o=F_o[-1>-1<0>1<1]$ en des séries non sinusoïdales: l'un générée par la fonction (Fig.9.b) $g_{o1}(x)=G_{o1}[-1>sinhx<1]$, et l'autre par $g_{o2}(x)$, générée par la presque-sinusoïde impaire de type exponentiel (Fig. 10.a) qui a un noyau de type $g(x)=K+cosh(x+T/4)$, où $K=constat$:

$$g_{o2}(x) = G_{o2}[-1<\frac{[\cosh 0.5 - \cosh(x-0.5)]}{\cosh 0.5 - 1}>0<\frac{[-\cosh 0.5 + \cosh(x+0.5)]}{\cosh 0.5 - 1}>1]$$

Dans l'intervalle $[-1, 1]$: $f_{o0}=g_{o10}=g_{o20}=0$, $\bar{f}_o(x) = \sum_{n=1}^{\infty} b_n \cos n\pi x = \frac{4}{\pi}\sum_{n=1}^{\infty}\frac{\sin(2n-1)\pi x}{2n-1}$ , donc:

$b_1=4/\pi$, $b_2=0$, $b_3=4/3\pi$, $b_4=0$, $b_5=4/5\pi$, $b_6=0$, $b_7=4/7\pi$, $b_8=0$, ... ,

$$\bar{g}_{o1}(x)_1 = \sum_{n=1}^{\infty}[c_n \sin(n\pi x)] = \pi \cdot \sinh 1 \sum_{n=1}^{\infty}\frac{2n\cdot(-1)^{n+1}}{1+n^2\pi^2}\sin(n\pi x)$$

$$\bar{g}_{o2}(x) = \sum_{n=1}^{\infty}d_n \sin n\pi x = \sum_{n=1}^{\infty}\frac{4\cdot\cosh 0.5}{\pi(\cosh 0.5 - 1)\cdot(2n-1)\cdot[1+(2n-1)^2\pi^2]}\sin n\pi x$$ , donc:

$d_1=8.8372\cdot(4/\pi)$, $d_2=0$, $d_3=8.8372\cdot(4/3\pi)\cdot 1/(1+4\pi^2)$, $d_4=0$, $d_5=8.8372\cdot(4/5\pi)\cdot 1/(1+9\pi^2)$, $d_6=0$, $d_7=8.8372\cdot(4/7\pi)\cdot 1/(1+16\pi^2)$, $d_8=0$, ... ,

Pour $K=8.8372$:
$B_1=K$, $B_2=0$, $B_3=K/(1+4\pi^2)$, $B_4=0$, $B_5=K/(1+9\pi^2)$, $B_6=0$, $B_7=K/(1+16\pi^2)$, $B_8=0$, $B_9=K/(1+25\pi^2)$, $B_{10}=0$, $B_{11}=K/(1+36\pi^2)$, $B_{12}=0$, ...

Une représentation graphique de la somme partielle $S_{12}$ est donnée à la figure 10.b.

Pour les valeurs de $L$ sous-unitaire, la quasi-sinusoïde exponentiel $g_{o2}(x)$ est approximée de manière satisfaisante par la quasi-sinusoïde quadratique $g_c(x)_1=C[1-x^2]_1$, décrite dans la figure 4.b de la section 3. Entre les coefficients des développements de la fonction $f_o(x)$ en les deux bases générées par les quasi-sinusoïdes $g_c(x)_1$ et $g_{o2}(x)$, les différences sont négligeables.

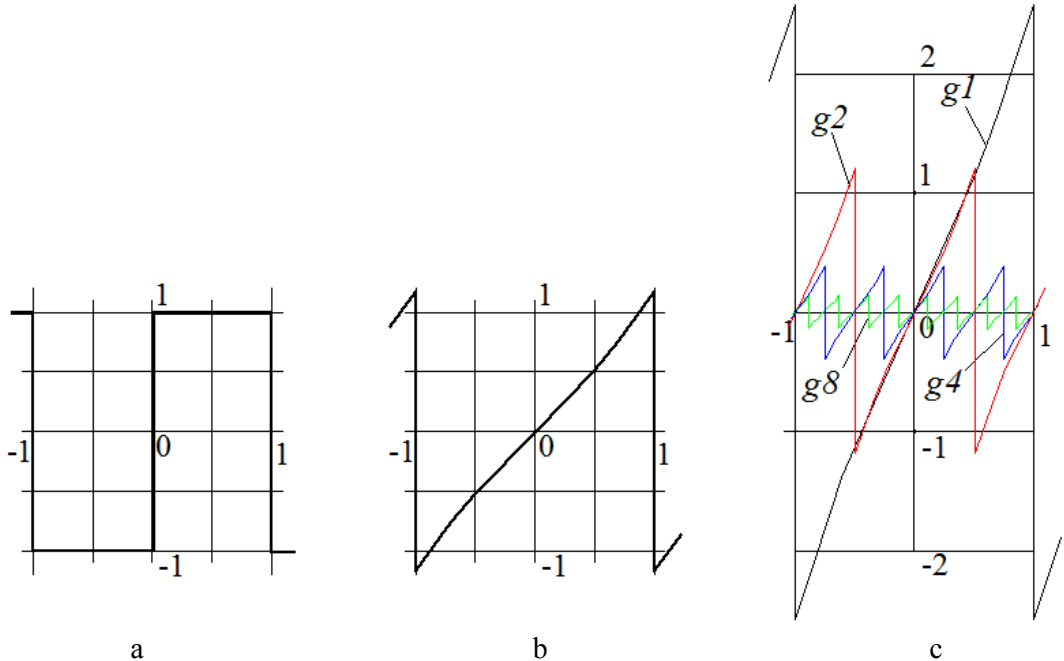

a         b         c

*Fig. 9. Le développement de la fonction $f_o(x)$ en la base $g_{o1}(x)$*
***a**: $f_{o1}(x)$ **b**: $g_{o1}(x)=sinh_1x$ **c**: les quasi-harmoniques $g_1(x)$, $g_2(x)$, $g_4(x)$, $g_8(x)$*

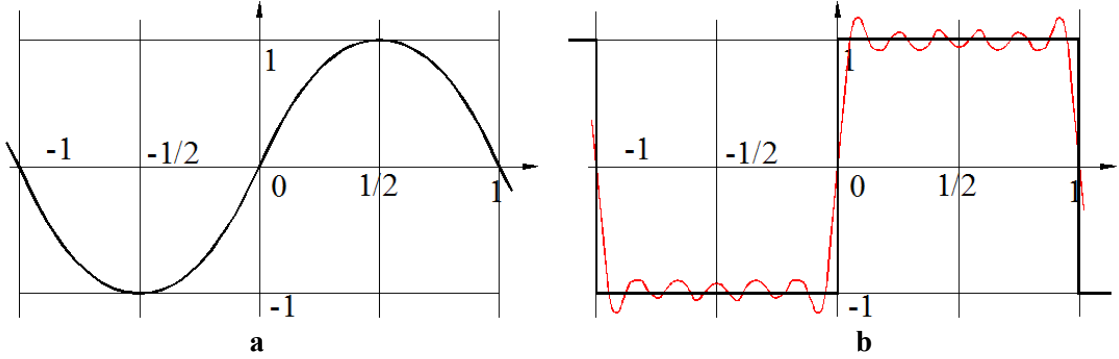

*Fig. 10. Le développement de la fonction $f_o(x)$ en la base $g_{o2}(x)$*
***a**: $g_{o2}(x)$  **b**: la somme partiels $S_{12}(x)$*

### 4. Bases orthogonales composées des fonctions périodiques non sinusoïdales

Ni les quasi-harmoniques pairs $g_n(x)−g_0$ ni les impairs $h_n(x)$, analysés dans les sections précédentes, ne sont pas orthogonaux les uns aux autres, ce qui ne permet pas le calcul des coefficients de ces développements à partir de formules similaires aux formules de Euler. Mais, tout quasi-harmonique pair est orthogonal par rapport aux tous quasi-harmoniques impairs. Cela nous permet, par le procédé d'orthogonalisation de Gram-Schmidt [4], de construire une base orthogonale (qu'on peut le normaliser par le même procédé) pour chacun des systèmes générés par les Fourier-fonctions $g_{Fn}(x)−g_0$ et $h_{Fn}(x)$. En les combinant et en ajoutant la fonction $f_0$, on obtient une base biorthogonale complète. Le procédé d'orthogonalisation de Gram-Schmidt ne prétend pas, pour les g-harmoniques de la base non orthogonale, la nécessité d'être continues (d'être Fourier-fonctions), mais cette fonctionnalité est imposée par notre intention de créer une base complète pour l'espace $L^2[−L, L]$. À cause de cela, nous considérons que les fonctions $g_{Fn}(x)$ et $h_{Fn}(x)$ sont des Fourier-fonctions par définition (1.d).

Par exemple, à partir des certaines deux fonctions $g(x)$ pair et $h(x)$ impair, qui ont la valeur moyenne nulle, définie sur un certain intervalle $[a, b]$, on obtient une base bi orthogonale formée par les fonctions $1$, $\Phi_n(x)$ et $\Psi_n(x)$, $n=1, 2, 3,...$, où:

$$\Phi_n(x) = g_{Fn}(x) - \sum_{j=1}^{n-1} \frac{\int_a^b g_n(x)\Phi_j(x)dx}{\|\Phi_j\|^2}\Phi_j(x) = g_{Fn}(x) - \sum_{i=1}^{n-1} C_{in}g_{Fi}(x) \text{ et}$$

$$\Psi_n(x) = h_{Fn}(x) - \sum_{j=1}^{n-1} \frac{\int_a^b h_n(x)\Psi_j(x)dx}{\|\Psi_j\|^2}\Psi_j(x) = h_{Fn}(x) - \sum_{i=1}^{n-1} D_{in}h_{Fi}(x) \quad (4.1)$$

avec $\|\Phi_j\|^2 = \int_a^b \Phi_j^2(x)dx$ et $\|\Psi_j\|^2 = \int_a^b \Psi_j^2(x)dx$

Ces considérations nous permettent de formuler le

**Théorème 6:** *Soit deux certaines F-fonctions $g_F(x)$-paire et $h_F(x)$-impaire de $L^2[a, b]$. Toute fonction $f(x)$ de $L^2[a, b]$ peut être développée en une série complète, basée sur le système bi orthogonal $1$, $\Phi_n(x)$ et $\Psi_n(x)$, où $\Phi_n(x)$ et $\Psi_n(x)$ sont généré par les fonctions $g_{Fn}(x)−g_0$ et $h_{Fn}(x)$ par une procédé d'orthogonalisation:*

$$\tilde{f}(x) = A_0 + \sum_{n=1}^{\infty}\left[A_n^0\Phi_n(x) + B_n^0\Psi_n(x)\right] \quad (4.2)$$

Grâce à l'orthogonalité du système, pour calculer les coefficients de ce développement, sont valables les formules de Euler:

$$A_0 = \frac{1}{b-a}\int_a^b f(x)dx$$

$$A_n^0 = \frac{1}{\|\Phi_n\|^2}\int_a^b f(x)\Phi_n(x)dx \quad n=1, 2, 3,... \tag{4.3}$$

$$B_n^0 = \frac{1}{\|\Psi_n\|^2}\int_a^b f(x)\Psi_n(x)dx \quad n=1, 2, 3,...$$

L'obtention de ces expressions est basée sur les relations:

$$\langle f(x)|\Phi_n(x)\rangle = \langle[f_e(x)+f_o(x)]|\Phi_n(x)\rangle = \langle f_e(x)|\Phi_n(x)\rangle$$

$$\langle f(x)|\Psi_n(x)\rangle = \langle[f_e(x)+f_o(x)]|\Psi_n(x)\rangle = \langle f_o(x)|\Psi_n(x)\rangle$$

qui sont vrais pour $\Phi_n(x)=cos(n\omega_0 x)$, respectivement $\Psi_n(x)=sin(n\omega_0 x)$, aussi.

On peut voir que les composantes $\Phi_n(x)$, d'ordre $n$ ($n=1, 2, 3,...$), du système orthogonal générés par les fonctions paires $g_{Fn}(x)-g_0$, ainsi que $\Psi_n(x)$, du système orthogonal généré par les fonctions impaires $h_{Fn}(x)$, sont des combinaisons linéaires entre la quasi-harmoniques d'ordre $n$ et les quasi-harmoniques d'ordre inférieur des respectifs développements non orthogonaux. Par conséquence, on peut établir une correspondance entre les coefficients $A_n$ et $B_n$ du développement en la base non orthogonale générés par les fonctions $g_F(x)$ et $h_F(x)$ et ceux du développement en la base orthogonale $\Phi_n(x)$ et $\Psi_n(x)$:

$$\tilde{f}(x) = A_0 + \sum_{n=1}^{\infty} A_n^0\left[g_{Fn}(x) - \sum_{i=1}^{n-1} C_{in}g_{Fi}(x)\right] + \sum_{n=1}^{\infty} B_n^0\left[h_{Fn}(x) - \sum_{i=1}^{n-1} D_{in}h_{Fi}(x)\right] =$$

$$= A_0 + \sum_{n=1}^{\infty}\left(A_n^0 - \sum_{i=n+1}^{\infty} A_n^0 C_{ni}\right)g_{Fn}(x) + \sum_{n=1}^{\infty}\left(B_n^0 - \sum_{i=n+1}^{\infty} B_n^0 D_{ni}\right)h_{Fn}(x) = A_0 + \sum_{n=1}^{\infty} A_n g_{Fn}(x) + \sum_{n=1}^{\infty} B_n h_{Fn}(x)$$

On constate que pour le calcul des coefficients $A_n = A_n^0 - \sum_{i=n+1}^{\infty} A_n^0 C_{ni}$ et $B_n = B_n^0 - \sum_{i=n+1}^{\infty} B_n^0 D_{ni}$, il est nécessaire de calculer certaines intégrales du type:

$$\int_a^b f_e(x)g_n(x)dx, \quad \int_a^b f_o(x)h_n(x)dx, \quad \int_a^b g_i(x)g_j(x)dx, \quad \int_a^b h_i(x)h_j(x)dx, \text{ pour } i,j=1, 2, 3, ...$$

Par cette méthode de calculer les coefficients, il n'est plus nécessaire de connaître les coefficients des développements en série sinusoïdale ni pour la fonction $f(x)$, ni pour les fonctions $g(x)$ et $h(x)$.

Nous allons exemplifier en construisant une base orthogonale, à partir d'une base générée par le système des fonctions créneau périodiques paires unitaires:
$f_e = F_e[-1>1<-1/2>-1<1/2>1<1]$

Dans le cas choisi ici, le calcul sera simplifié grâce aux propriétés de symétrie de la quasi-sinusoïde choisie. Grâce aux relations (4.1), il en résulte:

$$\Phi_1(x) = g_1(x), \qquad\qquad \|\Phi_1\|^2 = \int_{-1}^{1} g_1^2(x)dx = 2$$

$$\Phi_2(x) = g_2(x) - \frac{\int_{-1}^{1} g_2(x)g_1(x)dx}{2}g_1(x) = g_2(x), \qquad \|\Phi_2\|^2 = \int_{-1}^{1} g_2^2(x)dx = 2$$

$$\Phi_3(x) = g_3(x) - \frac{\int_{-1}^{1} g_3(x)g_1(x)dx}{2}g_1(x) - \frac{\int_{-1}^{1} g_3(x)g_2(x)dx}{2}g_2(x) = g_{F3}(x) - \frac{1}{3}g_{F1}(x) = g_3(x) - C_{13}g_1(x)$$

$$\|\Phi_3\|^2 = \int_{-1}^{1}\left[g_3(x) - \frac{1}{3}g_1(x)\right]^2 dx = \frac{16}{9}$$

$$\Phi_4(x) = g_4(x) - \frac{\int_{-1}^{1} g_4(x)g_1(x)dx}{2} g_1(x) - \frac{\int_{-1}^{1} g_4(x)g_2(x)dx}{2} g_2(x) -$$

$$-\frac{9\int_{-1}^{1} g_4(x)\left[g_3(x) - \frac{1}{3}g_1(x)\right]dx}{16}\left[g_3(x) - \frac{1}{3}g_1(x)\right] = g_4(x) - C_{14}g_1(x) - C_{24}g_2(x) - C_{34}g_3(x) = g_4(x)$$

$$\|\Phi_4\|^2 = \int_{-1}^{1} g_4^2(x)dx = 2$$

$$\Phi_5(x) = g_5(x) - \frac{\int_{-1}^{1} g_5(x)g_1(x)dx}{2} g_1(x) - \frac{9\int_{-1}^{1} g_5(x)\left[g_3(x) - \frac{1}{3}g_1(x)\right]dx}{16}\left[g_3(x) - \frac{1}{3}g_1(x)\right] =$$

$$= g_5(x) - \frac{9}{40}g_1(x) - \frac{3}{40}g_3(x) = g_5(x) - C_{15}g_1(x) - C_{35}g_3(x)$$

$$\|\Phi_5\|^2 = \int_{-1}^{1}\left[g_5(x) - \frac{17}{40}g_1(x) - \frac{3}{40}g_3(x)\right]^2 dx = \frac{1607}{800} \quad \text{et ainsi de suite.}$$

Par des relations similaires: $\Psi_n(x) = h_{F_n}(x) - \sum_{i=1}^{n-1} D_{in} h_{F_i}(x)$ est obtenu le système orthogonal $\Psi_n(x)$, à partir de $f_o(x)$, la fonction créneau périodique impaire unitaire:
$f_o(x) = F_o[-1 > -1 < 0 > 1 < 1]$.

### 5. Propriétés des séries de Fourier non sinusoïdales

Les travaux sur les développements en série de Fourier trigonométrique (sinusoïdales) ont montré qu'une fonction réelle $f(x):[-L, L]$, $2L$-périodique, peut être exprimée sous la forme d'une somme de ses projections sur les composantes d'une base orthogonal d'un espace de fonctions, s'il remplit plusieurs conditions. Les séries résultant à la suite de ces développements ont plusieurs propriétés: convergence, sommabilité, dérivabilité, intégrabilité.

Les séries de Fourier non sinusoïdales de la fonction $f(x)$ analysée dans les sections précédentes ont résulté à partir des série sinusoïdale de cette fonction, par une redistribution de ses coefficients. Cette redistribution reconstitue les coefficients des développements en des séries sinusoïdales des composantes d'une base complète des fonctions non orthogonaux. Par conséquent, les nouvelles méthodes de développement en série transfèrent des développements de Fourier sinusoïdales une série de conditionnalités et de propriétés. Sans aucun doute, ce sujet mérite une étude plus approfondie, mais pour le moment, nous nous résumons à quelques conclusions évidentes:

- toutes les fonctions $f(x)$ développées en séries de Fourier non sinusoïdales, ainsi que toutes les quasi-harmoniques ($g_n(x)$, $g_n(-x)$, $h_n(x)$, $\Phi_n(x)$, $\Psi_n(x)$, etc.) de ces développements sont des fonctions de Fourier.
- si la fonction $f(x)$ est intégrable, la suite des coefficients de son développement en série non sinusoïdale converge vers 0 (le théorème de Riemann-Lebesgue)
- si $\tilde{f}(x) = A_0 + \sum_{n=1}^{\infty}\left[A_n^0 \Phi_n(x) + B_n^0 \Psi_n(x)\right]$ est le développement de la fonction $f(x)$ en une base des fonctions orthogonales non sinusoïdales périodiques, alors

$$A_0^2 + \frac{1}{2}\sum_{n=1}^{\infty}\left[(A_n^0)^2 + (B_n^0)^2\right] = \frac{1}{2L}\int_{-L}^{L}|f(x)|^2 dx \quad \text{(le théorème de Parceval)}$$

- si les fonctions $f(x)$, $g(x)$, et $h(x)$, sont $2L$-périodiques, et dans l'intervalle $[-L, L]$ ils sont continus par morceaux et dérivables à gauche et à droite en tous les point de l'intervalle, la série $\hat{f}(x) = f_0 + \sum_{n=1}^{\infty} A_n g_n(x) + \sum_{n=1}^{\infty} B_n h_n(x)$ converge vers $\hat{f}(x_i) = \frac{1}{2}\left[\lim_{x \to x_i^+} f(x) + \lim_{x \to x_i^-} f(x)\right]$ dans tous les points $x_i$.

- la série de Fourier résulté à la suite d'un développement en une quelconque base *2L*-périodique non sinusoïdale d'une fonction *f(x)* *2L*-périodique, continûment dérivable par morceaux et continue, converge uniformément sur **R** vers cette fonction.
- la série non sinusoïdale de Fourier d'une fonction *2L*-périodique de carré intégrable sur une période, converge en norme $L^2$ vers la fonction considérée
- la série non sinusoïdale de Fourier d'une fonction de carré sommable converge presque partout vers cette fonction (le théorème de Carleson)
- deux fonctions *2L*-périodiques, ayant les mêmes coefficients de leurs développement en la même base non sinusoïdale Fourier, sont égales presque partout. Notamment, dans le cas continu par morceaux, elles coïncident en tous les points de [−*L, L*] sauf un nombre fini
- soit *f(x)* une fonction *2L*-périodique, continue dans l'intervalle [−*L, L*]. Son développement de Fourier $f(x) = f_0 + \sum_{n=1}^{\infty} A_n g_n(x) + \sum_{n=1}^{\infty} B_n h_n(x)$, sinusoïdal ou non, convergent ou non, peut être intégré terme par terme, entre toutes limites d'intégration:

$\int f(x)dx = d_0 + f_0 x + \sum_{n=1}^{\infty} A_n \int g_n(x)dx + \sum_{n=1}^{\infty} B_n \int h_n(x)dx$, où $d_0$ est une constante arbitraire.

- soit *f(x)* une fonction *2L*-périodique, continue dans l'intervalle [−*L, L*], avec *f(−L)=f(L)* et avec la dérivée *f'(x)* lisse par portions dans cet intervalle. Le développement de Fourier, sinusoïdal ou non, de la fonction *f'(x)*, peut être obtenu en dérivant terme par terme le développement de Fourier de la fonction *f(x)*. La série obtenue converge ponctuellement vers *f'(x)* en tous les points de continuité et vers [*f'(x)*+*f'(−x)*]/2 en ceux de discontinuité.

Si $f(x) = f_0 + \sum_{n=1}^{\infty} A_n g_n(x) + \sum_{n=1}^{\infty} B_n h_n(x)$, alors: $f'(x) = \sum_{n=1}^{\infty} A_n g_n'(x) + \sum_{n=1}^{\infty} B_n h_n'(x)$     (5.1)

La condition *f(−L)=f(L)* imposée dans cette affirmation est assez restrictive, ce qui réduit l'utilité du théorème. Nous pouvons contourner cette condition si nous prenons en compte le fait que le composant de type pair $f_e$ de la fonction *f(x)* vérifie toujours la condition de différentiabilité, et que le composant impair $f_o$ peut être écrit comme une somme de la fonction différenciable $f_{os}$ et de la fonction-rampe: $f_r = x \cdot f_o(L)/L$. Alors:

$$\frac{d}{dx} f_o(x) = \frac{d}{dx}\left[f_{os}(x) + \frac{f_o(L)}{L}x\right] = \frac{d}{dx} f_{os}(x) + \frac{f_o(L)}{L}$$

Par exemple, dans le cas du développement en série de Fourier sinusoïdale:

$$f'(x) = \frac{f_o(L)}{L} + \sum_{n=1}^{\infty}\left[b_n \omega_n \cos(\omega_n x) - \left(a_n \omega_n + 2(-1)^n \frac{f_o(L)}{L}\right)\sin(\omega_n x)\right]$$

ce qui nécessite la connaissance des conditions aux limites $f_o(-L)$ et $f_o(L)$. Cette relation permet de résoudre certaines équations différentielles en déterminant les coefficients du développement en série sinusoïdale de Fourier de la fonction inconnue (similaire au développement en série de Taylor).

### 6. Conclusions

Nous avons constaté dans les sections précédentes que toute fonction *f(x)*:[−*L, L*], *2L*-périodique, qui appartient au sous-espace $L^2$, peut être développée, d'une manière similaire à celle indiquée par Fourier il y a plus de 200 ans, dans une multitude de variantes, en des bases formées par la valeur moyenne de la fonction *f(x)* dans l'intervalle [−*L, L*] et deux ensembles de quasi-harmoniques: un ensemble de fonctions paires $g_n(x)$ et un ensemble de fonctions impaires $h_n(x)$, ($n=1, 2, 3, ..., \infty$), fonctions périodiques, avec la période $2L/n$. Dans le cas le plus général, les quasi-harmoniques fondamentales (pour $n = 1$) sont toutes fonctions qui satisfont les conditions de Dirichlet. Elles peuvent donc être des fonctions non trigonométriques et la base peut être non orthogonale.

Le développement en série sinusoïdale de Fourier n'est qu'un cas particulier de ce développement, à savoir le cas où les quasi-harmoniques fondamentales sont sinusoïdales: $g_{o1}(x)=sin(\omega_0 x)$ et $g_{e1}(x)=cos(\omega_0 x)$.

Ces résultats génèrent un large éventail de résultats théoriques. Tout d'abord, une nouvelle perspective extrêmement large s'ouvre dans l'analyse des espaces de fonctions, dans leur analyse spectrale, dans le développement de nouveaux types de transformations intégrales, dans la construction de systèmes de fonctions d'ondelettes, etc.

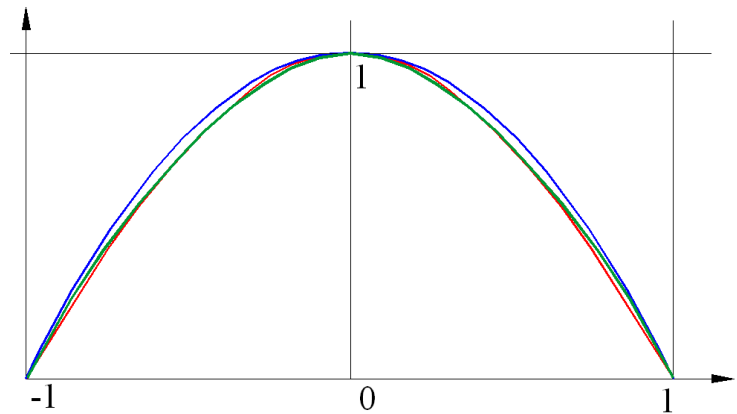

*Fig.11. Comparaison entre les courbes $cos(x\pi/2)$ ( rouge), $1-|x|^{1.75}$ (verte) et $1-x^2$ (bleue) dans l'intervalle $[-1, 1]$*

La comparaison des développements d'une fonction $f(x)$ en différentes bases complètes (le nombre de bases accessibles est devenu maintenant très important), permet de résoudre des nouveaux problèmes de convergence des séries numériques et des séries de fonctions et de trouver des nouvelles corrélations entre différents types de fonctions, etc. Deuxièmement, les conséquences pratiques sont également extrêmement importantes. Ce nouveau type de développement en série conduit à l'élaboration des nouvelles méthodes d'approximation de fonctions, dans lesquelles la précision de l'approximation peut être augmentée par la possibilité de choisir parmi un plus large éventail de possibilités. Par exemple, dans la figure 11, sont illustrés deux possibilités d'approximation d'une cosinusoïde.

Une perspective d'utilisation extrêmement prometteuse d'utilisation de ces types de développement en série est offerte par le domaine de la résolution numérique et analytique de larges catégories d'équations différentielles ordinaires et avec des dérivées partielles, linéaires et non linéaires.

## 7. Bibliographie